\input amstex
\documentstyle{amsppt}

\hsize=4.75in
\vsize=8in
\nologo

\def\I{{\Cal I}}
\def\T{{\Cal T}}
\def\A{{\Cal A}}
\def\K{{\Cal K}}
\def\B{{\Cal B}}

\def\O{{\Cal O}}
\def\C{{\Cal C}}
\def\D{{\Cal D}}

\def\S{{\Cal S}}
\def\L{{\Cal L}}

\def\H{{\Cal H}}
\def\Q{{\Cal Q}}
\rightheadtext {KMS states, entropy and the variational principle}
\leftheadtext{C.~Pinzari, Y.~Watatani, K.~Yonetani}
\topmatter
\title
KMS states, entropy and the variational principle\\
in full  $C^*$--dynamical systems
\endtitle
\author
C.~Pinzari\footnote{On leave of absence from Dipartimento di
Matematica, Universit\`a di Roma Tor Vergata, 00133 Roma. Supported by the
EU and
NATO--CNR.\hfil\break},
Y.~Watatani\footnote{Supported by the Grant--in--aid for Scientific 
Research of JSPS.\hfil\break},
K.~Yonetani$^2$
\endauthor
\affil
$^1$ Mathematics Department,\\
Massachusetts Institute of Technology, Cambridge, MA 02139,
USA. email: pinzari\@math.mit.edu\\
\\
$^2$ Graduate School of Mathematics,\\
Kyushu University,  Fukuoka 810-8560,
Japan. email:watatani\@rc.kyushu-u.ac.jp yonetani\@math.kyushu-u.ac.jp 
\endaffil

\abstract
To any periodic and full $C^*$--dynamical system $(\A, \alpha,
{\Bbb R})$, an invertible operator $s$ acting on the Banach space of
trace functionals of the fixed point algebra is canonically associated.
KMS states correspond to positive eigenvectors of $s$.
A Perron--Frobenius type theorem asserts the existence of
  KMS states at inverse temperatures equal the logarithms of  the inner
and outer
spectral radii of $s$ (extremal KMS states).
Examples arising from subshifts in symbolic dynamics, self--similar 
sets in fractal geometry and noncommutative metric spaces are discussed.
 
Certain subshifts are naturally associated to the system,
and criteria for the equality of their topological entropy and inverse 
temperatures of extremal KMS states are given. 

Unital completely positive maps $\sigma_{\{x_j\}}$
 implemented by  partitions of unity $\{x_j\}$ of grade $1$  are
considered, resembling the `canonical endomorphism' of the Cuntz
algebras.  
The relationship between the  Voiculescu  topological
entropy of $\sigma_{\{x_j\}}$ and the topological entropy of the
associated subshift is studied. Examples  where the equality holds
are discussed among 
Matsumoto algebras
associated to non finite type subshifts. 
In the general case $\text{ht}(\sigma_{\{x_j\}})$ is bounded by the sum
of the entropy of the subshift and a suitable entropic quantity
of the homogeneous subalgebra. Both summands are necessary. 

The measure--theoretic entropy
of $\sigma_{\{x_j\}}$, in the sense of Connes--Narnhofer--Thirring, is
compared to  the
classical measure--theoretic entropy of the subshift.

A noncommutative  analogue of  the classical variational principle for the
entropy  is obtained for the `canonical endomorphism' of
certain Matsumoto algebras. More generally, a necessary condition is
discussed. In the case of Cuntz--Krieger algebras an explicit
construction of the state with maximal entropy from the unique KMS state
 is done.
 \endabstract \endtopmatter \pagebreak
\document

\heading
Introduction
\endheading

Let $\A$ be a unital $C^*$--algebra endowed with a $2\pi$--periodic 
automorphic action $\alpha$ of ${\Bbb R}$. In algebraic statistical 
mechanics
elements of $\A$ represent kinematic observables of
an infinite quantum system, and $\alpha$ is the time evolution  of the 
system.
Equilibrium states of the system are  states on $\A$ which satisfy the 
KMS condition with respect to $\alpha$. 

We shall assume throughout the paper that $\alpha$ is $2\pi$--periodic, 
so it factors through
an action $\gamma$ of the circle ${\Bbb T}$, and also that $\gamma$ is 
{\it full\/}, in the
following sense.
Let $\A^k$ be the spectral subspace of elements
$a \in \A$ such that $\gamma_{z}(a) = z^ka$.
We assume that for all $k\in{\Bbb Z}$,
the closed linear span of $\{xy, x\in\A^k, y\in\A^{-k}\}$ is the fixed 
point
algebra $\A^0$. An example is given by a crossed product $C^*$--algebra 
$\A=\B\rtimes_{\beta}{\Bbb Z}$ by a single automorphism $\beta$, endowed
with the dual action $\gamma=\hat{\beta}$. All KMS states on $\A$ are
tracial, and are given by
 ${\Bbb T}$--invariant extensions to $\A$ of $\beta$--invariant tracial 
states on $\B$.

More generally, if ($\A$, $\gamma$) is not a dual
$C^*$--dynamical system, nontracial KMS states arise. 
Interesting examples are, in increasing generality,  the Cuntz--Krieger
algebras \cite{CK},  the Matsumoto algebras
associated with a subshift \cite{M}, and the
Pimsner algebras 
associated
with a full finite projective Hilbert $C^*$--bimodule \cite{P},
all endowed with the canonical gauge action.
While KMS states for the first two classes of $C^*$--algebras
are now well understood (see \cite{MFW}, \cite{E}, \cite{MYW}),
the third class of $C^*$--algebras is the main motivation of
the
present paper. (We shall see  that Pimsner $C^*$--algebras
are in fact a typical example, in the sense
that
 any unital, full and periodic  
$C^*$--dynamical system is isomorphic to a system constituted by
a Pimsner $C^*$--algebra associated
to a full, finite projective Hilbert bimodule, though not unique,
 and its canonical gauge action.)

We point out that, on the other hand, various  authors, regarding
the Cuntz--Krieger algebras as examples of noncommutative topological
dynamical systems, have computed the  Voiculescu topological entropy
\cite{V} of
the
so--called `canonical endomorphism' $\sigma$ (\cite{Ch},
\cite{BG}). However, among these, to the authors' knowledge,
it is only for the case of the Cuntz algebras $\O_d$ 
that a close relationship is known between KMS states,  
Voiculescu's topological entropy and measure--theoretic entropy in the
sense of
\cite{CNT}, see
\cite{Ch}. Choda's result states that the CNT entropy of $\sigma$ computed
with respect to the unique
KMS state equals the topological entropy of $\sigma$, which is, in turn,
$\log(d)$. This result can be regarded as a pivotal example
of noncommutative dynamical system for which a 
variational
principle for the entropy holds. Our
ultimate goal is that of
investigating the variational principle for the entropy and its
relationship with KMS states,  in
 full periodic $C^*$--dynamical systems.

KMS states for  $2\pi$--periodic actions have already been
 considered in the literature by several
authors. 
 Olesen and Pedersen   gave in \cite{OPI} an existence and
uniqueness theorem for
 KMS states of  the Cuntz algebras.   
 This was  generalized to the case of Cuntz--Krieger algebras
  by Enomoto, Fuji and the second--named author in \cite{EFW} 
and by Evans in \cite{E}. 

In 
\cite{BEH} Bratteli, Elliott and Herman construct, for any closed subset
$F$ of the extended real line, a simple $C^*$--algebra $\A_F$ endowed with
a  $2\pi$--periodic one--parameter group, for 
which
$F$ is precisely the set of inverse temperatures of KMS states,
 and such that for each $\beta\in F$, $\A_F$ has
a unique KMS state at inverse temperature $\beta$. 
We remark that if $F\subset{\Bbb R}$, 
then the ${\Bbb T}$--action is full, and if $F\subset(0,+\infty)$, 
$\A_F$ is purely infinite
(cf.  section 2). In a subsequent
paper \cite{BEK} Bratteli, Elliott and Kishimoto show that
even the set of KMS states with a specified inverse temperature
can be fairly arbitrary.

In a recent paper \cite{MWY} Matsumoto, Yoshida and the second--named
author 
study KMS states for the
Matsumoto $C^*$-algebras associated with a subshift in symbolic dynamics.
They develop a Perron--Frobenius type theorem
for a suitable positive operator naturally acting on a certain subalgebra, 
and show
that the 
logarithm of
its 
spectral radius 
arises as the inverse temperature of some KMS state. Furthermore
they  show a connection with the topological entropy of the 
underlying subshift. 

Our approach is close to that of \cite{MWY}, in that
we  emphasize the Perron--Frobenius  theory. The starting point is 
that to any periodic and full $C^*$--dynamical system we associate
certain completely positive maps on the underlying $C^*$--algebra, which
we 
interpret as being 
Perron--Frobenius type operators. 
KMS states correspond then to the positively scaled tracial states on the
fixed point
algebra.
We study the problem of existence of KMS states,  thus proving 
a Perron--Frobenius type theorem, and the relationship with the 
variational principle in ergodic theory.

The paper is organized  as
follows. 
In the first section  choose
 finite subsets $\{y_i\}$ and $\{x_j\}$ of $\A^1$
such that $\sum_j{y_j}^*y_j=I$ and $\sum_i{x_i}{x_i}^*=I$. Such multiplets 
exist because the group action is full. They 
can be regarded as
playing  the role of the canonical unitary  in 
 $\B\rtimes_\beta{\Bbb Z}$ implementing $\beta$.
We then consider two completely positive (cp) 
maps: $T_{\{y_i\}}: T\to\sum_i 
y_iT{y_i}^*$ and
$S_{\{x_j\}}: T\to\sum_j{x_j}^*Tx_j$ on $\A^0$, and also, by 
transposition,
 operators $t'$ and $s'$ which are inverses of one another
on the Banach space of trace functionals
on $\A^0$. 
These operators are  independent of the choice of the multiplets
$\{y_i\}$ and $\{x_j\}$.
 KMS states for the system at finite inverse temperatures
 then correspond to tracial states on $\A^0$ which are positively scaled 
by those cp maps, or, equivalently,
to tracial state eigenvectors of $s'$.

In the next section we show that, under the necessary
condition that the fixed point algebra has a tracial state, the inner
and outer spectral radii of $s'$  correspond to inverse temperatures
of  `minimal' and `maximal' KMS states (see Theorem 2.5 and  Corollary
2.6).
This can be regarded as a  Perron--Frobenius  theorem.
  The key point in the proof
 is that one needs to consider  
the trace functionals of the enveloping von Neumann algebra of $\A^0$,
endowed with its order structure.

In sections 3--5 we  discuss some examples.
In section 3 we apply our results to the Pimsner $C^*$--algebras generated
by finite 
projective
 Hilbert bimodules, and
 we thus deduce a criterion for existence of KMS states which applies, 
in 
particular,
to the case where the coefficient algebra is simple, unital and has a 
tracial state.

  In section 4 we construct Hilbert bimodules, and hence full
$C^*$--dynamical systems, via Pimsner's construction, naturally arising
from two different situations: subshifts of symbolic dynamics and 
self--similar sets in fractal
geometry.
In both cases the coefficient algebra is commutative, and the
corresponding Hilbert bimodules are described by a finite set of
endomorphisms. We show that
 in the former situation Pimsner's construction yields  the Matsumoto
$C^*$--algebras, while in the latter one gets a genuine Cuntz
algebra.  We also discuss a generalization of the latter example to
noncommutative metric spaces introduced by Connes \cite{Co}.
 It is interesting to compare our discussion with 
the papers  by
J\o rgensen-Pedersen \cite{JP} and  Bratteli-J\o rgensen \cite{BJ}, 
where 
the authors
 consider a relationship
between Cuntz algebras and multiresolutions in wavelet and fractal 
analysis.

In section 5 we look more closely at the subclass of the so called
Cuntz--Krieger bimodules (and the corresponding $C^*$--algebras). These
are bimodules for which the coefficient
algebra is a finite direct sum of unital simple $C^*$--algebras.
The leading and simplest example is, of course,  that of Cuntz--Krieger
algebras,
where each summand algebra is a copy of the complex numbers.
We show in particular that if each of the summands has a unique trace
and the defining $\{0, 1\}$--matrix $A$ is irreducible
then the associated Pimsner $C^*$--algebra has a unique KMS state
at inverse temperature $\log(r(A))$, where $r(A)$ is the spectral radius 
of
$A$.

In the next
section we associate to each pair ($\{y_i\}$, $\{x_j\}$)
of finite subsets of $\A^1$ as above, a pair of 
one--sided subshifts,
(${\ell}_{\{x_j\}}$, ${\ell'}_{\{y_i\}}$),
which roughly correspond to the operator $s'$ and its inverse $t'$.
We show that, under certain conditions,   the
topological entropies of these subshifts
 are precisely the minimal and maximal inverse temperatures
of KMS states. Furthermore we give a criterion for approximating 
such extremal temperatures with arbitrary (a priori non KMS)
tracial states satisfying suitable conditions.

In  section 7 we introduce a ucp  map $\sigma_{\{x_j\}}: 
T\to\sum_j{x_j}T{x_j}^*$
implemented by a multiplet $\{x_j\}$ of grade $1$
 as above, which should be compared with the map $S_{\{x_j\}}$. 
The main result of this section is the estimate 
$$\text{ht}(\sigma_{\{x_j\}})\leq
h_{\text{top}}(\ell_{\{x_j\}})+\text{ht}(\{\phi_{x_i, x_j}\})$$
where the l.h.s. is the Brown--Voiculescu topological entropy
 \cite{B}, \cite{V} 
of $\sigma_{\{x_j\}}$
and the second
summand at the r.h.s. is
the topological 
entropy, suitable defined, of the set of contractions $\phi_{x_i,
x_j}:T\to{x_i}^*Tx_j$ of
the homogeneous $C^*$--subalgebra $\A^0$. 
Both summands at the r.h.s. of  this inequality are necessary. Indeed
when $\A=\A^0\rtimes_\alpha{\Bbb Z}$ then the associated subshift
is trivial so its entropy is zero,
and the above inequality, when combined with monotonicity of
topological entropy
(\cite{B}, \cite{V}) 
leads to 
  Brown's result 
$$\text{ht}_\A(\text{Ad}(u))=\text{ht}_{\A^0}(\alpha),$$
where $u\in\A$ is a unitary implementing $\alpha$ \cite{B}.
Another extreme case is that of the Cuntz--Krieger algebras $\O_A$.
Now the second summand vanishes and the previous estimate 
yelds the result by Boca and Goldstein \cite{BG} that
$\text{ht}(\sigma_{\{x_j\}})=\text{ht}(\ell_{\{x_j\}})=\log(r(A))$,
see Corollary 7.8.

We then focus our attention on those algebras for which
$\text{ht}(\{\phi_{x_i, x_j}\})=0$ and we show that
if the $x_j$'s have pairwise orthogonal ranges
$$\text{ht}(\sigma_{\{x_j\}})=h_{\text{top}}(\ell_{\{x_j\}}).$$
We next  discuss new examples  of this occurrence among 
 Matsumoto algebras \cite{M} associated to a subshift. 

Our assumption of orthogonality introduces a certain trivialization
to the classical situation. In fact, in this case the algebra 
of continuous functions on the one--sided  subshift
$\ell_{\{x_j\}}$, together
with  the endomorphism induced by the left shift epimorphism,
 sits naturally inside the 
noncommutative 
dynamical
system ($\A$, $\sigma_{\{x_j\}}$). Therefore monotonicity of topological
entropy implies that  the  topological entropy of
$\ell_{\{x_j\}}$ is $\leq$ 
 the  topological entropy of the noncommutative  subshift
$\sigma_{\{x_j\}}$, thus leading 
to 
the equality, see Theorem 7.7. We discuss new examples of this occurrence
among the Matsumoto algebras associated to certain non finite type
subshifts.

In section 8 we investigate the 
CNT
dynamical entropy of $\sigma_{\{x_j\}}$. We show that, under the orthogonality assumption,
if $\phi$ is a $\sigma_{\{x_j\}}$--invariant state of $\A$ centralized by 
$\C(\ell_{\{x_j\}})$ then $h_\phi(\sigma_{\{x_j\}})\geq
h_\mu(\ell_{\{x_j\}})$,
where $\mu$ is the shift--invariant probability measure on
$\ell_{\{x_j\}}$
obtained restricting $\phi$. We also find a condition on $\sigma_{\{x_j\}}$
under which any such $\mu$ arises as the restriction of some $\phi$.
This enables us to obtain a variational priciple for certain systems
$(\A, \sigma_{\{x_j\}})$
for which $\text{ht}(\{\phi_{x_i,x_j}\})=0$.
More precisely, our variational principle asserts, for those
systems,
  the existence
of
$\sigma_{\{x_j\}}$--invariant states
of $\A$ with respect to which the CNT dynamical entropy equals the
Voiculescu  topological
entropy of $\sigma_{\{x_j\}}$.

In the last section we establish a closer relationship between 
KMS states and states with maximal entropy. The main point is that
a KMS state $\omega$ is to be understood as a  quasi--invariant measure
for the
noncommutative shift $\sigma_{\{x_j\}}$, as
$\omega\circ\sigma_{\{x_j\}}$ and $\omega$ are equivalent.
In classical ergodic theory, measures with this property are called conformal,
and play an important role, as they lead to measures with maximal entropy.
We thus show an explicit general way of constructing $\sigma_{\{x_j\}}$--invariant
measures from KMS states. We then consider 
basic examples, which we may think of as being  noncommutative
  Markov shifts:
systems $(\A, \gamma)$ containing some Cuntz--Krieger algebra $\O_A$
in a way that $\gamma$ restricts to the canonical gauge action on $\O_A$. 
We show that the $\sigma_{\{x_j\}}$--invariant
state $\phi$ previously derived from a KMS
state with maximal entropy restricts, 
on the algebra of continuous functions  on the classical Markov subshift
$\ell_A$,
 to the 
unique invariant measure $\mu$ with maximal entropy.
We thus conclude that if $\text{ht}(\{\phi_{x_i, x_j}\})=0$,  
$$h_\mu(\ell_A)=h_\phi(\sigma_{\{x_j\}})=h_{\text{top}}(\sigma_{\{x_j\}})=
h_{\text{top}}(\ell_A)=\log(r(A)).$$
This yields a generalization of Choda's result \cite{Ch} to the
Cuntz--Krieger
algebras, and Matsumoto algebras associated to certain non finite type
subshifts.

\heading
1. The scaling property
\endheading

 Recall that a state $\omega$
over a $C^*$--algebra $\A$ endowed with a one--parameter automorphism
group $\alpha$ is
called a KMS state at inverse temperature $\beta\in{\Bbb
R}$ if 
$$\omega(a\alpha_{i\beta}(b))=\omega(ba),\eqno(1.1)$$
for all $a, b$ in a dense $^*$--subalgebra of $\A_\alpha$, the set
 of entire
elements for $\alpha$ (which is in fact a dense $^*$--subalgebra).

We will only consider  $2\pi$--periodic one--parameter groups, i.e.
groups for which $\alpha$
comes from an action $\gamma$ of ${\Bbb T}$ by
$\alpha_t:=\gamma_{e^{it}}$. 
Furthermore, in view of applications to the algebras generated by 
Hilbert
bimodules, we will assume that $\A$ is unital and that
the group  action is full, in the sense explained in the introduction. 
Then we note that
the spectral subspace $\A^k$, for positive $k$, is in fact 
the linear span of the product set of $k$ copies of  $\A^1$.  
Moreover, by  definition of full $C^*$--dynamical system, there exist, for
all 
$n\in{\Bbb N}$,
finite subsets $\{y_i\}$ and $\{x_j\}$
of $\A^n$ such that $\sum_i{y_i}^*y_i=I$ and $\sum_j x_j{x_j}^*=I$.
We define correspondingly, for any tracial state $\tau$ on $\A^0$,
 $$\delta_n(\tau):=\tau(\sum_iy_i{y_i}^*)$$
$$\epsilon_n(\tau):=\tau(\sum_j {x_j}^*x_j).$$
We shall usually write $\delta(\tau)$ and $\epsilon(\tau)$
for $\delta_1(\tau)$ and $\epsilon_1(\tau)$ respectively.
\medskip

\noindent{\bf 1.1. Lemma} {\sl Let $(\A, \gamma, {\Bbb T})$ be a full 
$C^*$--dynamical
system, with $\A$ unital, and let $\tau$ be a tracial
state
on $\A^0$. Then $\delta_n(\tau)$ and $\epsilon_n(\tau)$ do not 
depend on the finite subsets $\{{y_i}\}$ and $\{{x_j}\}$ of 
$\A^n$
satisfying
the above relations. If in paricular $\{y_i\}$, $\{x_j\}\subset\A^1$,
 one has, for all $n\in{\Bbb N}$,
$$\|\sum {(x_{j_1}\dots x_{j_n})}^*x_{j_1}\dots x_{j_n}\|^{-1/n}
\leq\delta_n(\tau)^{1/n}\leq
\|\sum y_{i_1}\dots y_{i_n}({y_{i_1}\dots y_{i_n}})^*\|^{1/n}$$
$$\|\sum y_{i_1}\dots y_{i_n}({y_{i_1}\dots y_{i_n}})^*\|^{-1/n}
\leq\epsilon_n(\tau)^{1/n}\leq
\|\sum {(x_{j_1}\dots x_{j_n})}^*x_{j_1}\dots x_{j_n}\|^{1/n}$$
}\medskip

\noindent{\it Proof} We shall only prove the statements relative to
$\{y_i\}$, those relative to $\{x_j\}$ can be proved similarly.
Let $\{z_1,\dots, z_q\}\subset\A^n$ be another multiplet
satisfying $\sum_k{z_k}^*z_k=I$, and write $z_k=\sum_i a_{k,i}y_i$, 
where 
$a_{k,i}:=z_k{y_i}^*\in\A^0$. Then 
$$\tau(\sum_k{z_k}{z_k}^*)=
\tau(\sum_{k,i,j}a_{k,i}y_i{y_{j}}^*{a_{k,j}}^*)=
\tau(\sum_{i,j,k} y_i{y_{j}}^*{a_{k,j}}^*a_{k,i})=$$
$$\tau(\sum_{i,j,k} y_i{y_{j}}^*y_j{z_k}^*z_k{y_i}^*)=\tau(\sum_i
y_i{y_i}^*).$$
Note that $\delta_n(\tau)\leq\|\sum_i y_i{y_i}^*\|$. Furthermore
$$1=\tau(\sum_{i,j}x_j{y_i}^*y_i{x_j}^*)=\tau(\sum y_i{x_j}^*x_j{y_i}^*)
\leq
\|\sum {x_j}^*x_j\|\tau(\sum y_i{y_i}^*).$$ The conclusion follows
choosing subsets in $\A^n$ of the form 
$\{y_{i_1}\dots y_{i_n}\}$ and $\{x_{j_1}\dots x_{j_n}\}$, where 
$\{y_i\}$, $\{x_j\}\subset\A^1$ and satisfy $\sum_i{y_i}^*y_i=I$ and
$\sum_j x_j{x_j}^*=I$.\medskip
 
Let, for
$\lambda>0$, $\T\S_\lambda$ be the set of 
tracial states $\tau$ on $\A^0$ for which
$$\lambda\tau(xy^*)=\tau(y^*x),\quad x, y\in\A^1.\eqno(1.2)$$
We shall show that a state of $\A$ satisfies
the KMS condition  w.r.t. $\alpha$ 
if and only if its restriction to $\A^0$ is an element of some  
$\T\S_\lambda$.
We start with the following characterization of $\T\S_\lambda$.\medskip 

\noindent{\bf 1.2. Lemma} {\sl Let $(\A, \gamma, {\Bbb T})$ be a full 
$C^*$--dynamical 
system over 
a unital $C^*$--algebra.
For a tracial state $\tau$ on $\A^0$ and $\lambda>0$, the following 
conditions are equivalent:
\roster
\item $\tau\in\T\S_\lambda,$
\item 
$\tau(\sum_iy_ia{y_i}^*)=\lambda^{-1}\tau(a),\quad a\in\A^0,$
\item
$\tau(\sum_j {x_j}^*ax_j)=\lambda\tau(a),\quad a\in\A^0.$
\endroster
Here $\{y_i\}$ and $\{x_j\}$ are finite subsets of $\A^1$ satisfying
respectively
$$\sum_i {y_i}^*y_i=I,\quad \sum_j x_j{x_j}^*=I.$$
 $\lambda$ is uniquely determined by $\tau$:
$\lambda=\epsilon(\tau)=\delta(\tau)^{-1}.$}\medskip

\noindent{\it Proof} $(1)\to(2)$ and $(1)\to(3)$ are obvious.
We show that $(2)\to(1)$: for $x, y\in\A^1$, $y^*x\in\A^0$, so
$$\tau(y^*x)=\lambda\tau(\sum_i y_i y^*x{y_i}^*)=
\lambda\tau(\sum_i x{y_i}^*y_i y^*)=\lambda\tau(xy^*).$$ 
One similarly proves
 that $(3)\to(1)$.\medskip

The following result characterizes the set
of
tracial states on $\A^0$ which gives rise to KMS states for $(\A,
\gamma)$.
Let $F_0:\A\to \A^0$ denote the projection onto the fixed point algebra
obtained overaging over the circle
group action.
\medskip

\noindent{\bf 1.3. Proposition} {\sl 
The maps $\omega\to\omega\upharpoonright_{\A^0}$,
$\tau\to\tau\circ F_0$ set up a bijective correspondence between 
the set 
of KMS states $\omega$ for $(\A, \gamma)$ at inverse temperature 
$\beta$   and the set
$\T\S_{e^\beta}$.
 }\medskip

\noindent{\it Proof} If $\omega$ is a KMS state at inverse 
temperature $\beta$ then 
the KMS condition $(1.1)$ can be formulated, equivalently, for any pair
$x,y$ 
in the dense linear span of
 the $\A^k$'s. Therefore the 
 restriction $\tau$ of $\omega$ to $\A^0$ is a tracial state such that 
$\omega=\tau\circ F_0$
since $\omega$ is $\gamma$--invariant. Furthermore, if $x, y\in\A^1$,
$$\omega(y^*x)=\omega(x\alpha_{i\beta}(y^*))=e^{\beta}\omega(xy^*),$$
therefore $\tau\in\T\S_{e^\beta}$.
Conversely, if this condition is satisfied by some tracial state $\tau$ 
on $\A^0$, then
$\omega:=\tau\circ F_0$ is an extension of $\tau$ to a state on 
$\A$, and it is not difficult to check
that $\omega(y^*x)=\omega(x\alpha_{i\beta}(y^*))$ for $x,y\in\A^1$, and  
hence inductively for 
$x,y\in\A^1\dots\A^1=\A^k$.
If $x\in\A^k$, $y\in\A^h$, $h\neq k$ then 
$\omega(y^*x)=0=\omega(x\alpha_{i\beta}(y^*)),$ and the proof
 is complete. \medskip

\noindent{\it Remark\/}
A simple argument shows that the sequences 
$a_n:=\inf\{\epsilon_n(\tau),\tau\in\T\S(\A^0)\}$ and 
$b_n:=\sup\{\epsilon_n(\tau),\tau\in\T\S(\A^0)\}$ are respectively
supermultiplicative and
submultiplicative, therefore the sequences ${a_n}^{1/n}$, ${b_n}^{1/n}$
converge, and  $\lim_n{a_n}^{1/n}=\sup{a_n}^{1/n}$ 
and $\lim_n{b_n}^{1/n}=\inf {b_n}^{1/n}$.
\medskip

\noindent{\bf 1.4. Corollary} {\sl
If 
\roster
\item    $\lim_n\inf\{\epsilon_n(\tau),\tau\in\T\S(\A^0)\}^{1/n}>1$ 
(e.g. 
$\inf\{\epsilon(\tau), \tau\in\T\S(\A^0)\}>1$) 
 then every
 KMS state on $(\A, \gamma)$ has  positive inverse temperature,
\item  $\lim_n\sup\{\epsilon_n(\tau),\tau\in\T\S(\A^0)\}^{1/n}<1$ 
(e.g. $\sup\{\epsilon(\tau),\tau\in\T\S(\A^0)\}<1$)
 then every
 KMS state on $(\A, \gamma)$ has  negative inverse temperature,
\item  $\lim_n\inf\{\epsilon_n(\tau),\tau\in\T\S(\A^0)\}^{1/n}=
\lim_n\sup\{\epsilon_n(\tau),\tau\in\T\S(\A^0)\}^{1/n}=1$ (e.g.
$\epsilon(\tau)=1$  for all 
$\tau\in \T\S(\A^0)$)
then every KMS state on $(\A, \gamma)$ is tracial.
\endroster}

We next give  a criterion of faithfulness
for KMS states.\medskip

\noindent{\bf 1.5. Proposition} {\sl Let $(\A, \gamma)$ be a 
full periodic $C^*$--dynamical
system with $\A$ unital, and consider the $^*$--monomomorphism
$\alpha: \A^0\to M_p(\A^0)$ associated to a set $\{y_i\}_{i=1}^p $ of 
$\A^1$ such that 
$\sum_{i=1}^p {y_i}^*y_i=I$ and defined by $\alpha(a)=(y_ia{y_j}^*).$
 If $\A^0$ has no proper closed ideal $\I$
such that $\alpha(\A^0)\cap M_p(\I)=\alpha(\I)$ (e.g. $\A$ is simple), 
then any 
KMS state
of $(\A, \gamma)$ is faithful.}\medskip

\noindent{\it Proof} Let $\tau$ be the restriction of a KMS state 
$\omega$
to $\A^0$. Then $$\I:=\{a\in\A^0: \tau(a^*a)=0\}$$ is a closed ideal of 
$\A^0$. Since
for $x, y\in\A^1$, $a\in\A^0$,
$$(xay^*)^*(xay^*)\leq \|x\|^2ya^*ay^*,$$
we have, by $(1.2)$, that $xay^*\in\I$ if $a\in\I$. This  yields
$\alpha(\I)\subset 
M_p(\I)\cap\alpha(\A^0)$. We  show the reverse inclusion.
Let $a\in\A^0$ be such that $y_ia{y_j}^*\in\I$, $i,j=1,\dots,p$.
Then 
$\delta(\tau)\tau(a^*{y_i}^*y_ia{y_j}^*y_j)=\tau((y_ia{y_j}^*)^*y_ia{y_j
}^*)=0$.
Hence, summing up, we see that $a\in\I$, and this shows that 
$\alpha(\A^0)\cap M_p(\I)\subset\alpha(\I)$. It follows from our
assumption that $\I=\{0\}$, 
as, 
clearly, $\I\neq\A^0$.
Now the canonical conditional expectation $F_0:\A\to\A^0$ is faithful, 
so $\omega$ is faithful.
\medskip

We conclude this section recalling from \cite{GP}  a criterion for
 pure infinity of unital $C^*$--algebras, which we shall need in the 
sequel. We refrain from giving here the proof. We only point out  
that the arguments  essentially 
  go back to Cuntz' proof
of pure infinity of $\O_d$ \cite{C}. Also, the result
 is a 
generalization
of R\o rdam's result (cf. \cite{R}) about pure infinity of crossed 
products by proper corner
endomorphisms.

Following \cite{R}, we say that a $C^*$--algebra $\B$ has the 
comparability property if $\B$ has at least one tracial state,
and furthermore a projection $e\in\B$ is equivalent to a subprojection 
of $f$ if $\tau(e)<\tau(f)$ for all tracial states of 
$\B$. Assume that our $C^*$--algebra $\A$ has a nonunitary isometry 
$S$
in some $\A^n$, $n>0$, and also that the fixed point algebra $\A^0$ has
the comparability property. Then for all tracial states $\tau$ on $\A^0$
one has, by Lemma 1.1,
$$\delta_n(\tau)=\tau(SS^*)<1.\eqno(1.3)$$
(Note that, by Corollary 1.4 and the following Proposition 2.2,
all KMS states have positive inverse temperatures.)
It is then natural to ask under which conditions
$(1.3)$ guarantees the existence of a nonunitary isometry in $\A^n$,
or, better, pure infinity of $\A$. 
\medskip

\noindent{\bf 1.6. Theorem} {\sl \cite{GP} Let $(\A, \gamma,{\Bbb T})$ 
be
a full 
$C^*$--dynamical 
system,
 and assume
that $\A^0$  is unital, simple, separable, of real rank 
zero, and that  every $M_n(\A^0)$ has the comparability propery.
 If, for some $n>0$, 
$$\sup\{\delta_n(\tau), \tau\in\T\S(\A^0)\}<1,$$
 then $\A^n$ contains a nonunitary isometry. Furthermore, $\A$ is simple
and purely infinite.}\medskip

\heading
2. A Perron--Frobenius  theorem
\endheading

We  associate to each pair of finite subsets
 $\{y_i\}$, $\{x_j\}\subset\A^1$ satisfying $$\sum_i{y_i}^*y_i=I,$$
 $$\sum_j x_j{x_j}^*=I,$$ a corresponding  pair of completely
  positive maps  $T=T_{\{y_i\}}$ and $S=S_{\{x_j\}}$
on the homogeneous subalgebra $\A^0$:
$$T(a):=\sum_iy_ia{y_i}^*, \quad a\in\A^0,$$
$$S(a):=\sum_j{x_j}^*a{x_j}, \quad a\in\A^0.$$

Let ${T}'$, ${S}':{\A^0}^*\to{\A^0}^*$ denote the  Banach space adjoints 
of $T$ and $S$ respectively, and  let  $\T(\A^0)\subset {\A^0}^*$
the Banach subspace of trace  functionals.
Then one has the following result.\medskip

\noindent{\bf 2.1. Proposition} {\sl For any  finite
subset $\{y_i\}$ (resp. $\{x_j\}$) of $\A^1$ satisfying 
$\sum_i{y_i}^*y_i=I$ 
(resp.  $\sum_jx_j{x_j}^*=I$) the associated operator
 ${T}'$ (resp. ${S}'$)
leaves  $\T(\A^0)$ stable. Let $t'$ (resp. $s'$)
 be the  restriction of $T'$ to $\T(\A^0)$. Then $t'$ (resp. $s'$) does
 not depend on the set $\{y_i\}$ (resp. $\{x_j\}$) 
satisfying $\sum{y_i}^*y_i=I$ (resp. $\sum x_j{x_j}^*=I$). 
Furthermore $t'$ and 
$s'$ are inverses of one another.}\medskip

\noindent{\it Proof} It is easy to check that ${T}'$  
transforms trace functionals into trace functionals. Let
 $\{{y'}_k\}$ be another finite subset of $\A^1$
such that $\sum{{y'}_k}^*{y'}_k=I$, and write ${y'}_k=\sum a_{k,i}y_i$,
with $a_{k,i}={y'}_k{y_i}^*\in\A^0$. Then for any
 $\tau\in \T(\A^0)$, $a\in\A^0$, 
$$\tau(\sum_k{y'}_ka{{y'}_k}^*)=\sum_{k,i,j}\tau(a_{k,i}y_ia{y_j}^*{a_{k,j}}^*)=
\sum_{i,j,k}\tau(y_ia{y_j}^*{a_{k,j}}^*a_{k,i})=$$
$$\sum_{i,j}\tau(y_ia{y_j}^*y_j{y_i}^*)=\tau(\sum_iy_ia{y_i}^*).$$
Finally note that
 $$t's'(\tau)=s'(\tau)\circ T=\tau\circ ST=\tau,\quad \tau\in \T(\A^0),$$
by that trace property of $\tau$. Likewise,
$$s't'(\tau)=\tau,\quad\tau\in \T(\A^0).$$\medskip

Note that KMS states of $(\A, \gamma)$ correspond precisely to
the tracial state eigenvectors for $s'$ (or $t'$).
The following result, which has its own interest, explains 
why $\delta(\tau)=\epsilon(\tau)^{-1}$ when $\tau$ corresponds to a 
KMS state.   \medskip

\noindent{\bf 2.2. Proposition} {\sl  The map $h:\T\S(\A^0)\to\T\S(\A^0)
$
taking a tracial state $\tau$ to 
$$h(\tau)=\delta(\tau)^{-1}{t'}(\tau)$$
is a homeomorphism of $\T\S(\A^0)$ endowed with the weak$^*$--topology.
KMS states of $(\A, \gamma)$ correspond, as in Prop. 1.3, to fixed points
of $h$. 
The  inverse of $h$ is the map 
$$k(\tau)=\epsilon(\tau)^{-1}{s'}(\tau).$$
We have: 
$$\epsilon_n(h^n(\tau))=\delta_n(\tau)^{-1},$$
$$\delta_n(k^n(\tau))=\epsilon_n(\tau)^{-1},$$
thus if $\tau\in\T(\A^0)$ is the restriction of a  KMS state,
$\epsilon(\tau)=\frac{1}{\delta(\tau)}$.
 }\medskip

\noindent{\it Proof}
Clearly  $h(\tau)$ is a tracial state when
$\tau$ is. Furthermore $\tau\to\delta(\tau)$ is a positive valued
continuous
function on a compact set, therefore $h$  is continuous.
For the same reason $k$ is continuous on $\T\S(\A^0)$, and, by the trace
property,
 $h$ and $k$ are inverses of one another.\medskip

Our next aim
is to look more closely at 
 the spectrum $\sigma(s')$ of $s'$. 
The previous proposition shows  that
$$0\notin\sigma(s')=\sigma(t')^{-1}.$$
So we can define the inner and outer spectral radius of
$s':$
$$r_{\text{max}}(s'):=\max\{|\lambda|, \lambda\in\sigma(s')\}$$
$$r_{\text{min}}(s'):={r(t')_{\text{max}}}^{-1}=\min\{|\lambda|, 
\lambda\in\sigma(s')\}.$$
We give some estimates for $r_{\text{min}}(s')$ and 
$r_{\text{max}}(s')$.\medskip

\noindent{\bf 2.3. Proposition} {\sl Let $\{y_i\}$, $\{x_j\}\subset 
\A^1$
satisfy $\sum_i {y_i}^*y_i=I$ and $\sum_j{x_j}{x_j}^*=I$. Then one has
$$r_{\text{min}}(s')\geq\lim_n\|\sum_{i_1,\dots,i_n} y_{i_1}\dots 
y_{i_n}(y_1\dots y_{i_n})^*\|^{-1/n},$$
$$r_{\text{max}}(s')\leq
\lim_n\|\sum_{j_1,\dots,j_n}(x_{j_1}
\dots{x_{j_n}})^*x_{j_1}\dots{x_{j_n}}\|^{1/n}.$$}\medskip

\noindent{\it Proof}
Note that $$T(a)^*T(a)\leq 
\|\sum_i y_i{y_i}^*\|T(a^*a),\quad a\in\A^0,$$ which, together with
$T(I)=\sum_i{y_i}{y_i}^*$, implies 
$\|T\|=\|\sum_i y_i{y_i}^*\|$,
 hence inductively 
$$\|{T}^n\|=\|\sum_{i_1,\dots,i_n} y_{i_1}\dots y_{i_n}(y_{i_1}\dots
y_{i_n})^*\|.$$
Taking the $n$-th root and passing to the limit, 
one gets the spectral radius
of $T$:
$$r(T)=\lim_n\|\sum_{i_1,\dots,i_n} y_{i_1}\dots 
y_{i_n}(y_{i_1}\dots y_{i_n})^*\|^{1/n}$$
Similarly, one has
$$r(S)=
\lim_n\|\sum_{j_1,\dots,j_n}(x_{j_1}\dots x_{j_n})^*x_{j_1}\dots
x_{j_n}\|^{1/n}.$$
The proof is completed recalling that $s'$ is the restriction of $S'$ 
to a closed subspace, so 
$$r_{\text{max}}(s')\leq r(S')=r(S)$$
and similarly 
$$r_{\text{min}}(s')=r_{\text{max}}({t'}^{-1})\geq r(T)^{-1}.$$
\medskip

We next show  that
the inner and outer spectral radii of $s'$ correspond to inverse 
temperatures of KMS states, or, in other words, that they are
in the point spectrum of $s'$,  with corresponding  positive eigenvalues.
  The fact that the outer spectral
radius is in the point spectrum   was first proved in 
\cite{MWY} for the Matsumoto $C^*$--algebras associated with subshifts
\cite{M}. 
The key point in our situation  is that one needs to 
consider 
the trace functionals of $\A^0$
endowed with the  order structure which arises 
when we extend such traces to normal traces on the enveloping von Neumann
algebra. 

We anticipate the following,  possibly known,
 lemma.\medskip

\noindent{\bf 2.4. Lemma} {\sl Let $\sum_n\phi_n$ be a series
of   normal
 linear functionals 
on a von Neumann algebra $M$ weakly convergent to $\phi$.
If  each of the the absolute values 
$|\phi_n|$ is  tracial then 
$$|\phi|\leq\sum_n |\phi_n|.$$}\medskip

\noindent{\it Proof}
Let $\tau$ be a  positive tracial
linear functional, then 
$$|\tau (xy)| \leq ||y||\tau (|x|)$$
(see for example \cite{T}).
Consider the polar decompositions 
$$\phi_n(x)=|\phi_n|(xu_n),$$ 
$$\phi(x)= |\phi|(xu).$$
Then, for a positive $x$, we have
$$|\phi|(x) =\phi(xu^*)= \sum_n \phi_n(xu^*)
= \sum_n|\phi_n|(xu^*u_n)
\leq$$
$$\sum_n||\phi_n|(xu^*u_n)| \leq\sum_n |\phi_n|(x).$$\medskip
 
We are now in the position of proving our main result of this
section.\medskip

\noindent{\bf 2.5. Theorem} {\sl Let $(\A, \gamma,{\Bbb T})$ be a full 
$C^*$--dynamical system, and assume that
$\A$ is unital, and that $\A^0$ has a tracial state. Then
$r_{\text{min}}(s')$ and $r_{\text{max}}(s')$ 
 are eigenvalues of $s'$  with corresponding tracial state eigenvectors.}
\medskip

\noindent{\it Proof} We first show that $r_{\text{max}}(s')$ 
is a spectral value for $s'$ and then that it is in fact an eigenvalue
with a tracial state eigenvector. A similar argument will prove
that $r(t')$ is an eigenvalue for $t'={s'}^{-1}$ with a tracial
 state eigenvector.

By the uniform boundedness theorem,
there exists a sequence $\{z_n\}$ of complex numbers
such that $|z_n| \to r_{\text{max}}(s')^+$ and 
$\|R(z_n)\tau_0 \| \rightarrow \infty$ for some $\tau_0 \in \T(\A^0)$, 
 where $R(z)$ is the resolvent of 
$s'$ in $z$. Since $\T(\A^0)$ is linearly spanned by its 
tracial states, we may assume
that $\tau_0$ is a tracial state. 
 Consider, for $|z|>r_{\text{max}}(s')$, the Neumann series:
$$
             R(z) = \sum_{k=0}^{\infty} z^{-(k+1)}{s'}^k.\eqno(2.1)
$$
By the previous lemma, on the enveloping von Neumann algebra of $\A^0$,
$$
    |R(z)\tau_0| \leq \sum _{k=0}^{\infty} |z|^{-(k+1)}{s'}^k(\tau_0) = 
R(|z|)\tau_0,
$$
so  $\|R(|z_n|)\tau_0\|\to\infty$, and this 
shows that $r_{\text{max}}(s')\in\sigma(s')$.
$(2.1)$ also shows   that $R(\lambda)\tau_0$
is a nonzero positive functional for $\lambda>r_{\text{max}}(s')$,
 hence,  arguments similar to those of Lemma 3.1 in \cite{MWY} prove that 
 $$\tau_n:=
\frac{1}{\|R(|z_n|)\tau_0\|}R(|z_n|)\tau_0$$ is a sequence of tracial
states such that every weak$^*$--limit point of it is a  tracial state
eigenvector
with eigenvalue $r_{\text{max}}(s')$.\medskip

The previous theorem can be considered as an analogue of the
Perron--Frobenius theorem for matrices with nonnegative entries.\medskip

\noindent{\bf 2.6. Corollary} {\sl Let $(\alpha, {\Bbb R})$ be a 
$2\pi$--periodic
one--parameter automorphism group
of a unital $C^*$--algebra $\A$, such that the induced ${\Bbb
T}$--action $\gamma$  is full.  If $s'$ is defined as above, 
relatively to $\gamma$,
 then
the set of inverse temperatures of KMS states is a closed subset of
the interval 
$[\log(r_{\text{min}}(s')), \log(r_{\text{max}}(s'))]$ containing the
extreme points.}\medskip

\noindent{\it Proof} The subset of $\T\S(\A^0)$ corresponding to KMS
states
is weakly$^*$--compact by Prop. 2.2, furthermore the map
$\epsilon:\T\S(\A^0)\to{\Bbb
R}^+$ defined at the
beginning of section 1 is weakly$^*$--continuous. It follows that
the set of elements of the form $\log(\epsilon(\tau))$, when $\tau$
ranges over all
tracial states on $\A^0$ corresponding to KMS states, is compact. Now this
set
is precisely the set of possible inverse temperatures by Prop. 1.3. The
rest follows from the previous Theorem.\medskip

A KMS state of $(\A, \alpha)$ at inverse temperature 
$$\beta_{\text{min}}:=\log(r_{\text{min}}(s'))$$ or 
$$\beta_{\text{max}}:=\log(r_{\text{max}}(s'))$$ will be called
extremal.
Let $\beta$ be the inverse temperature of a KMS state, 
and set, as in the previous section, $$a_n=\inf\{\epsilon_n(\tau),
\tau\in\T\S(\A^0)\},$$
$$b_n=\sup\{\epsilon_n(\tau),\tau\in\T\S(\A^0)\}.$$
Then for all $n$,
${a_n}^{1/n}\leq e^\beta\leq{b_n}^{1/n}$, so 
$$\lim_n 1/n\log(a_n)\leq\beta\leq\lim_n1/n\log(b_n).$$
It is then natural to ask for which tracial states $\tau$, the sequence
$1/n\log(\epsilon_n(\tau))$
 approximates  the maximal
or the minimal inverse temperature. In   section 5 we shall 
give a sufficient condition.
\medskip

We conclude this section with the discussion of two examples
known in the literature. The first example, arising from ergodic theory,
 shows that in general, at a fixed
inverse temperature, there may be more than
one KMS state. \medskip

\noindent{\bf 2.7.  Example} 
Let  $(X, T)$ be a topological dynamical system:
$X$ is a  compact metric space endowed with a  
homeomorphism $T$. We  suppose that $X$ is not a finite set.  
Then it is well known that the $C^*$--algebra 
$\A=\C(X)\rtimes_{\alpha_T}{\Bbb Z}$ is simple if and only if $T$ is 
minimal,
i.e. there is no nontrivial closed 
subset $F\subset X$ such that $T(F)=F$. 
Here $\alpha_T$ is 
the automorphism of $\C(X)$ defined by $\alpha_T(f)=f\circ T^{-1}$.
 Tracial states
on $\C(X)\rtimes_{\alpha_T}{\Bbb Z}$ are in one--to--one correspondence
with $T$--invariant probability measures on $X$,
while there is no nontracial KMS state on $(\A, \gamma)$. 
The operator $s'$ therefore has spectrum contained in the unit circle.
However, $s'$ is the Banach space adjoint of $\alpha_T$, so its spectrum 
is the same as that of $\alpha_T$ which must be equal to ${\Bbb T}$ by 
simplicity
of $\A$ \cite{OP}.
There is an important example, due to Furstenberg, of a 
minimal
 analytic
diffeomorphism $T$  of ${\Bbb T}^2$ with nonunique 
invariant measures (see, e.g, \cite{Ma}), which thus leads to an example 
of nonuniqueness of tracial states on the simple crossed product
$C^*$--algebra $\C({\Bbb T}^2)\rtimes_{\alpha_T}{\Bbb Z}$.
\medskip

The next example shows that the set of inverse
temperatures can in general be an arbitrary closed subset
of ${\Bbb R}$.\medskip

\noindent{\bf 2.8. Example} In \cite{BEH} Bratteli, Elliott and Herman
construct an example of a simple $C^*$--algebra $\B$ endowed with ${\Bbb 
T}$--action for which the set of possible inverse
temperatures can be any arbitrary closed subset $F$ of ${\Bbb 
R}\cup\{+\infty, -\infty\}$. For each temperature
the corresponding state is unique. More in detail, $\B$ is obtained by 
cutting down the crossed product 
$\A\rtimes_{\alpha}{\Bbb Z}$
of an AF--algebra 
by some projection $P$ in $\A$. If neither $+\infty$ nor $-\infty$ belongs
to 
$F$, $\A$ itself is simple,
and this implies that the ${\Bbb T}$--action is full since $P$ is a full 
projection. If moreover
$F\subset(0, +\infty)$, one can choose $\alpha$ so that $\alpha(P)<P$,
see \cite{BEH}, 
hence $\alpha(\B^0)\subset\B^0$.
Then $\rho:=\alpha\upharpoonright_{\B^0}$ is a proper corner 
endomorphism  of $\B^0$, and one has
$\B=\B^0\rtimes_{\rho}{\Bbb N}$. Now by a result of R\o rdam 
\cite{R}, $\B$ is purely infinite.\medskip

\heading 
3. KMS states of the Pimsner algebras
\endheading 

In this section we discuss an application of the results of the previous 
section
 to the $C^*$--algebra $\O_X$
 associated to a 
 Hilbert $C^*$--bimodule $X$ over a  $C^*$--algebra  $\B$.
We refer the reader to \cite{P} for the construction of $\O_X$.
We just recall that both $X$ and $\B$ embed isometrically respectively as a 
Hilbert bimodule in and a $C^*$--subalgebra of $\O_X$.
We shall always assume that $X$ is finite projective and full, and that 
$\B$ is
unital. Therefore  any finite basis $\{x_j\}$ of $X$ yields, in $\O_X$,
the
relation 
$$\sum_{j} x_j{x_j}^*=I.$$
Furthermore any finite subset $\{y_i\}$ of $X$ such that $\sum_i <y_i,
y_i>=I$ yields in $\O_X$:
$$\sum_i {y_i}^*y_i=I.$$
$\O_X$ is endowed with a canonical gauge action $\gamma$ such that
$\gamma_z(x)=zx$, $z\in{\Bbb T}$, $x\in X$. Therefore we can conclude that 
($\O_X$, $\gamma$) is a full periodic $C^*$--dynamical system.  

We start proving that  systems of this form are typical
examples, in the sense that we can always easily associate
to  any unital, full, periodic  
$C^*$--dynamical system $(\A, \gamma)$, a  finite projective Hilbert
$C^*$--bimodule 
$X$ such that $\A=\O_X$ and $\gamma$ is the canonical ${\Bbb T}$--action.
  We should note, however, that  the Hilbert bimodule $X$  and its
coefficient 
$C^*$-algebra 
are, in general, not unique. In fact our construction leads to
a  maximal Hilbert bimodule.
In applications, it may be more convenient to start with smaller Hilbert bimodules.
It is well known the case of Cuntz--Krieger algebras discussed in
\cite{P}, where
the coefficient
algebra is finite--dimensional. 
In section 4 we shall discuss  the more general situation of Matsumoto algebras,
and we will construct
 natural minimal generating Hilbert bimodules.
\medskip

\noindent{\bf 3.1. Theorem} {\sl Let $(\A, \gamma)$ be a full 
$C^*$--dynamical system over ${\Bbb T}$ and assume that
$\A$ is unital.  Then there exists a full finite projective Hilbert 
$C^*$--bimodule $X$ over a 
unital $C^*$--algebra  $\B$ such that $(\A, \gamma)$ can be identified 
with 
$\O_X$, endowed with its canonical gauge action.}
\medskip

\noindent{\it Proof} Choose finite subsets 
$\{y_i\}$ and $\{x_j\}$
of $\A^1$ such that $\sum_i{y_i}^*y_i=I$ and $\sum_j x_j{x_j}^*=I$.
Let $\B$ be the fixed point algebra $\A^0$.  Set $X = \sum_j x_j\B$.
Then $X$ is a right Hilbert $\B$--module in $\A$ with  $\B$-valued inner
product: $<x,y>_\B = x^*y$.  The condition $\sum_j x_j{x_j}^*=I$ shows 
that
$\{x_j\}$ is a finite basis of $X$,  
and $X$ is full by  $\sum_i{y_i}^*y_i=I$.   
Left $\B$-action is given by $\phi (b)x = bx$ for $b \in \B$ and 
$x \in X$. For $x = \sum_j x_jb_j$, we have 
$ bx =\sum_j \sum_i x_ix_i^*bx_jb_j \in X$.
Thus $\phi$ is well defined.
Let $\hat{x}$ the image of an element $x\in X$ in $\O_X$.
By the universality of the Pimsner algebras, 
there exists a surjective $*$-homomorphism
$\varphi : \O_X \rightarrow \A$ such that $\varphi(\hat{x}_j) = x_j$.  Let  
$F_0 : \A \rightarrow \A^0$ and $E : \O_X \rightarrow \O_X^0$ be the 
natural
conditional expectations.  Since $\varphi E = F_0 \varphi$ and 
$E$ is faithful, 
$\varphi$ is injective by a well known argument.  Thus 
$\varphi$ is the 
 desired 
isomorphism.
\medskip

A KMS state on $(\O_X, \gamma)$ at some inverse 
temperature $\log(\delta)\in{\Bbb R}$ restricts to a tracial state 
$\tau$ on the coefficient algebra $\B$
satisfying, for $a\in \B$, $\tau(\sum_i<x_i, ax_i>)=\delta\tau(a)$, with 
$\{x_i\}$ a right basis of $X$. We show that conversely any such 
 trace extends to a KMS state.
\medskip

\noindent{\bf 3.2. Lemma} {\sl Let $\{x_1,\dots,x_d\}$ be a basis of 
$X$.
Any tracial state $\tau$ on the coefficient algebra  $\B$ satisfying 
$$\tau(\sum_i<x_i, ax_i>)=\delta\tau(a), \quad a\in\B,$$
with $\delta>0$, extends uniquely to a KMS state on $\O_X$ at inverse 
temperature $\log(\delta)$.
This extension is faithful if $\tau$ is faithful.}
\medskip

\noindent{\it Proof} We first prove uniqueness. A KMS state for $(\O_X, 
\gamma)$ is determined
by its restriction to the homogeneous $C^*$--subalgebra, and, by the 
trace--scaling property of the operator
$S_{\{x_j\}}$ on that subalgebra, it is in fact determined
by its restriction to the coefficient algebra $\B$. Conversely, if one is 
given a tracial state
$\tau_0$ on $\B$ as required then it is easy to check that
$$\tau_n(a):=\frac{1}{\delta}\tau_{n-1}(\sum {x_i}^*ax_i),\quad 
a\in\L(X^n), n\geq 1,$$
is a sequence of tracial states (faithful if $\tau_0$ is faithful) such 
that $\tau_{n+1}
\upharpoonright_{\L(X^n)}=\tau_n$, which thus gives rise
to a tracial state $\tau$ on the homogeneous subalgebra positively  scaled
by
$S_{\{x_j\}}$. 
Therefore $\tau$ extends
to a KMS state on $\O_X$ at inverse temperature $\log(\delta)$.\medskip

We now apply the results of the previous section to the Pimsner 
$C^*$--algebras.
\medskip

\noindent{\bf 3.3. Theorem} {\sl Let $\B$ be a unital 
$C^*$--algebra  with a tracial state, and let $X$ be a 
full finite projective Hilbert $C^*$--bimodule over  $\B$. Assume
 that
 for every tracial state $\tau$ on $\B$, and any basis $\{x_i\}$ of $X$,
$\tau(\sum_i <x_i, x_i>)>0$.
Then $\O_X$ has a KMS state.
\roster
\item Let $s': \T(\O_X^0)\to \T(\O_X^0)$ be the (invertible) operator
obtained restricting to $\T({\O_X}^0)$
 the Banach space adjoint of $S(a)=\sum_i {x_i}^*a{x_i}$, where $\{x_i\}$ 
is a
 basis of $X$. 
Then the set of  possible inverse 
temperatures is a closed subset of 
$[\log(r_{\text{min}}(s')), \log(r_{\text{max}}(s'))]$ 
containing   the extreme points.
\item  If $\O_X$ is ${\Bbb T}$--simple, every KMS state is faithful,
\item if $\tau(\sum_i<x_i, x_i>)>1$ for all $\tau\in\T\S(\B)$ then every 
KMS state has positive inverse
temperature,
\item if $\tau(\sum_i<x_i, x_i>)=1,$ for all $\tau\in\T\S(\B)$ 
then every KMS state is a tracial state,
\item if $\tau(\sum_i<x_i, x_i>)<1$ for all $\tau\in\T\S(\B)$ 
then every KMS state has a negative inverse temperature,
\item if $\B$ has a unique trace then ${\O_X}^0$ has a unique trace,
so $\O_X$ has a unique KMS state.
\endroster}\medskip

\noindent{\it Proof}   The function 
taking  a tracial state
$\tau$ to the tracial state 
$$a\in\B\to(\tau(\sum<x_i, x_i>))^{-1}\tau(\sum <x_i, ax_i>)$$
is weakly$^*$--continuous, so, by  the Schauder--Tychonov 
fixed point theorem, there is  a tracial state
$\tau$ on $\B$ such that $$\tau(\sum_i<x_i, ax_i>)=
\tau(\sum_i<x_i, x_i>)\tau(a),\quad a\in\B.$$ 
We can now extend such a  $\tau$ to a KMS state on $(\O_X, \gamma)$, by 
the
previous lemma. $(1)$ follows from Corollary 2.6. Since 
${\O_X}$ is ${\Bbb T}$--simple, $\O_X^0$ has no ideal
of the kind described  in our faithfulness
criterion, Prop. 1.5,  hence $(2)$ follows.
 $(3)$--$(5)$
 follow from Corollary 1.4.
Finally, if $\B$ has a unique tracial state then so does ${\O_X}^0$ since 
it is 
an inductive limit of $C^*$--algebras stably isomorphic to
$\B$ itself. Therefore $\O_X$ has a unique KMS state.
\medskip

\heading 4. Examples of full dynamical systems arising from
 subshifts,  self-similar sets and noncommutative metric spaces
\endheading

In this section we continue our discussion of  examples of  full
$C^*$--dynamical
systems obtained via   
Pimsner's construction. 
We start considering two different examples of Hilbert bimodules
both described
by families
of 
${}^*$-endo\-mor\-phi\-sms on  
commutative $C^*$-algebras, arising respectively from symbolic dynamics
and fractal
geometry. We shall also discuss a generalization of the latter example to
noncommutative metric spaces.
\medskip

\noindent{\bf 4.1. Subshifts in symbolic dynamics and Matsumoto algebras}

We recall the construction of the Matsumoto algebra $\O _{\Lambda}$ 
  associated with a two--sided subshift $\Lambda$, \cite{M}. 
Fix a finite discrete  set 
$\Sigma = \{ 1,2,...,d \}$, and let $\Sigma^{\Bbb{Z}}$  
 be the infinite product space endowed with 
the product topology. We denote by $\sigma$ the shift homeomorphism on 
$\Sigma^{\Bbb{Z}}$ defined by $(\sigma(x))_i = x_{i+1}$. For a
shift--invariant closed subset 
$\Lambda$ of $\Sigma^{\Bbb{Z}}$, the topological dynamical system 
$(\Lambda, \sigma\upharpoonright_{\Lambda})$ is called a subshift.  
We denote by $\Lambda_+$ the set of one--sided sequences 
$x \in \Sigma^{\Bbb{N}}$ such that $x$ appears in $\Lambda$.
For example, $\Sigma^{\Bbb N}={\Sigma^{\Bbb Z}}_+$. We shall still denote
by $\sigma$ the left shift epimorphism of $\Sigma^{\Bbb N}$.  
The dynamical system $({\Lambda_+}, \sigma\upharpoonright_{\Lambda_+})$ is
called the one--sided subshift
associated to $\Lambda$.  
A finite sequence $\mu = (\mu_1, \dots, \mu_k) $ of elements 
$\mu_j \in \Sigma$ is called a word. We denote by $|\mu|$ the length $k$ 
of $\mu$. 
For $k \in \Bbb{N}$, let 
$
\Lambda^k 
= \{ \mu | \mu \text{ is a word with length } k \text{ appearing in some 
} x \in \Lambda \}, 
$
\ $\Lambda_l = \bigcup_{k=0}^l \Lambda ^k$
and $\Lambda ^* = \bigcup_{k=0}^\infty \Lambda ^k$, where  $\Lambda^0$ 
denotes the set constituted by  
the empty word. 
\par

Let $\{e_1,...,e_d \}$ be an orthonormal basis of a $d$-dimensional 
Hilbert space 
$H ={\Bbb C}^d$. Let $ F^0$ be the one dimensional space ${\Bbb 
C}\Omega$ spanned by 
a normalized vector $\Omega$, called the vacuum vector, and let $F^k$ be
the
Hilbert 
space spanned by 
the vectors $e_\mu = e_{\mu_1} \otimes \cdots \otimes e_{\mu_k}$ for 
$\mu =(\mu_1,\dots \mu_k) \in \Lambda^k$. Consider the subspace 
$F_{\Lambda} = \oplus^{\infty}_{k=0} F^k$ of the full Fock space of $H$.

The creation operator $T_{\nu}$ by $e_{\nu}$ on $F_{\Lambda}$, for $\nu 
\in \Lambda ^{*}$, is 
defined by  
$$
T_{\nu}\Omega = e_{\nu}
\quad 
\text{ and }
\quad
T_{\nu}e_{\mu} =
\cases
e_{\nu} \otimes e_\mu, & \text{ if } \nu\mu \in \Lambda ^{*}\\
0                      & \text{ otherwise}
\endcases
$$
The unital $C^*$-subalgebra $\T_{\Lambda}$ of the algebra of bounded 
linear operators on $F_\Lambda$ 
generated by $\{T_i | i = 1,...,d \}$ is called the Toeplitz algebra
associated with $\Lambda$, and 
contains the algebra 
$\K (F_{\Lambda})$ of compact operators on $F_{\Lambda}$. The 
Matsumoto 
algebra $\O _{\Lambda}$ associated with the subshift $\Lambda$ is the 
quotient algebra
$\T _{\Lambda} / \K (F_{\Lambda})$. It is  generated by the quotient image
$\{S_i | i =
1,...,d 
\}$ of
$\{T_i | i = 1,...,d \}$.
The unitary representation of ${\Bbb T}$ on $F_\Lambda$ defining the
grading implements an automorphic action of ${\Bbb T}$ on $\T_\Lambda$
leaving $\K(F_\Lambda)$ stable. We thus obtain an automorphic ${\Bbb
T}$--action $\gamma$ on $\O_\Lambda$ such that 
$$\gamma_z(S_i)=zS_i,\quad z\in{\Bbb T}, i=1,\dots,d.$$
As $\sum _i S_iS_i^* = I$ and 
$\sum _i S_i^*S_i \geq I$, ($\O_\Lambda$, $\gamma$) is a full periodic
$C^*$--dynamical system.
We set $S_\mu = S_{\mu_1} \dots S_{\mu_k}$,  
for  $\mu = (\mu_1, \dots, \mu_k) \in \Lambda^* $.

  For each $i= 1, \dots , d$, we define a (not necessarily unital) 
  ${}^*$-endomorphism $\rho_i$ on 
$\ell ^{\infty}({\Lambda_+})$ by  
$$ 
(\rho_i(f))(x) =
\cases
f(i,x_1, \dots ) & \text{ if } (i,x_1, \dots ) \in  {\Lambda_+}\\
0                      & \text{ otherwise}
\endcases
$$
for $f \in  \ell ^{\infty}({\Lambda_+})$ and $x \in  {\Lambda_+}$.
Then we  have $\cap _i \text{ker} \rho_i = 0$.  

Consider the functions $q_{\mu} \in \ell ^{\infty}({\Lambda_+})$,  for
 $\mu =(\mu_1,\dots \mu_k) \in \Lambda^k$, defined by
 $$  
 q_{\mu}(x) =
\cases
1  & \text{ if } \mu x \in {\Lambda_+}\\
0                      & \text{ otherwise}
\endcases 
$$
Thus  $\rho_i(I) = q_i$. 
 (We should note that  $\rho_i(I) = q_i$  
is not a continuous function, in general.)

 Let $\Lambda_A$ be  the Markov subshift
defined by a $d\times d$ matrix $A=(a_{i,j})$ with entries in $\{0,1\}$
and with no
zero rows or
columns: 
$$\Lambda_A=\{x\in\Sigma^{\Bbb Z}: a_{x_i, x_{i+1}}=1, i\in{\Bbb Z}\}.$$ 
Then
each $\rho_i$ preserves $\C({{\Lambda_A}_+})$. 
For Markov subshifts, Matsumoto's construction yields the Cuntz--Krieger
algebras, see \cite{M},:
$$\O_{\Lambda_A}\simeq\O_A.$$
\medskip
 
 \noindent{\bf 4.2. Proposition} {\sl Let $\Lambda_A$ be  the Markov
subshift 
defined by a  matrix $A=(a_{i,j})\in M_d(\{0,1\})$ and set $\B =
\C({\Lambda_A}_+)$.  Consider 
the right Hilbert
 $\B$-module $X = \oplus _{i=1}^d q_i\B$ and   the $^*$-homomorphism 
$\phi : \B \rightarrow \L_{\B}(X)$ 
 given by the diagonal matrix
$\phi (a) = \text{diag} (\rho_i(a))_i$.  Then the Pimsner algebra $\O_X$ 
is isomorphic to the Cuntz-Krieger algebra $\O_A.$}
\medskip

\noindent{\it  Proof}  The commutative $C^*$-algebra $\D_A$ generated by 
$\{S_{\mu}S_{\mu}^* ; \mu \in \Lambda^* \}$ is isomorphic to
$\C({\Lambda_A}_+)$
via an isomorphism which identifies $S_{\mu}S_{\mu}^*$ with the
characteristic function
 $p_{\mu} = \chi_{[\mu]}$ of the
cylinder set $[\mu] = \{ x \in{\Lambda_+} ; 
(x_1, \dots, x_k)= (\mu_1, \dots, \mu_k) \} $
for $\mu \in \Lambda^k$. 
 Then the  ${}^*$-endomorphism $\gamma_i$ on $\D_A$
defined by $\gamma_i(T) = S_i^*TS_i$ for $T \in \D_A$  
corresponds 
to the  ${}^*$-endomorphism $\rho_i$ on $\C( {\Lambda_+})$. 
We have 
$$\rho_i(p_r) = \delta_{i,r}\sum_{\{j: a_{i j} =1\}} p_j \  \text{ and } 
\rho_i(p_{\mu}) = \delta_{i,\mu_1}p_{(\mu_2, \dots, \mu_k)}
$$
for $\mu \in \Lambda^k$ with $k \geq 2$. The corresponding formulae  for 
$\gamma_i(S_{\mu}S_{\mu}^*)$ hold.  Consider the right  Hilbert 
 $\D_A$-module $Y = \oplus _{i=1}^d S_i\D_A$ and   the $^*$-homomorphism 
$\phi : \D_A \rightarrow \L_{\D_A}(Y)$ 
 given by the diagonal matrix
$\phi (a) = \text{diag} (\gamma_i(a))_i$.
Since $(S_{\mu}S_{\mu}^*)S_i = S_i\gamma_i(S_{\mu}S_{\mu}^*)$, the 
Hilbert
bimodules $X$ and $Y$ are isomorphic.
Now a standard argument shows that $O_{A} \cong \O_Y \cong \O_X$. 
\medskip

If $\Lambda$ is a general subshift, the endomorphisms $\rho_i$,
$i=1,\dots, d$,
do not leave 
 $\C({\Lambda_+})$ stable.
Thus we should  replace $\C({\Lambda_+})$ by some unital 
$C^*$-subalgebra of 
$\ell ^{\infty}({\Lambda_+})$ which is invariant under $\rho_i$, 
$i=1, \dots d$. We shall choose, to this aim, the smallest such
$C^*$--algebra, which is related to the  Krieger left cover 
of a  sofic subshift and the past equivalence relation considered by
Matsumoto.

   Let $A({\Lambda_+})$ be the unital $C^*$-subalgebra
of $\ell^\infty({\Lambda_+})$ generated by 
$\{q_{\mu} ; \mu \in \Lambda^* \}$.
Since $\rho_i(I) = q_i$  and $\rho_i(q_{\mu}) = q_{\mu i}$, 
it is clear that $A({\Lambda_+})$
is the smallest unital $C^*$-subalgebra of 
$\ell ^{\infty}({\Lambda_+})$ which is invariant under $\rho_i$, 
$i=1, \dots d$. 
 \medskip
 
 \noindent{\bf 4.3. Theorem} {\sl Let 
 $\B = A({\Lambda_+})$ be the commutative $C^*$--algebra associated to a
subshift $\Lambda$ as above. Consider the Hilbert right
 $\B$-module $X = \oplus _{i=1}^d q_i\B$ and   the $^*$-homomorphism 
$\phi : \B \rightarrow \L_{\B}(X)$ 
 given by the diagonal matrix
$\phi (a) = \text{diag} (\rho_i(a))_i$.  Then the Pimsner algebra $\O_X$ 
is isomorphic to the Matsumoto algebra $\O_{\Lambda}.$}
\medskip

\noindent{\it Proof} For $l \in {\Bbb N}$, let  
 $ A_l$ be the $C^*$-subalgebra of  $\O _{\Lambda}$ generated by  
               $\{S_{\mu}^*S_{\mu}, \mu \in \Lambda_l \} $
and $A_{\Lambda}$ be the $C^*$-subalgebra  of  $\O _{\Lambda}$
 generated by elements  
               $\{S_{\mu}^*S_{\mu}, \mu \in \Lambda^* \}.$
Then  
   $(A_l)_l$ is an increasing sequence of commutative 
finite--dimensional algebras and   $A_{\Lambda} = \varinjlim A_l$. 
Similarly, let $ A_l({\Lambda_+})$ be the $C^*$-subalgebra of 
$\ell ^{\infty}({\Lambda_+})$  generated by  
               $\{q_{\mu}, \mu \in \Lambda_l \} $.  
Then  $(A_l({\Lambda_+}))_l$ is an increasing sequence of commutative 
finite--dimensional algebras and   
$A({\Lambda_+}) = \varinjlim A_l({\Lambda_+})$.  
For $x \in {\Lambda_+},$ let 
$\Lambda_l(x) = \{\mu \in \Lambda_l ;  \mu x \in  {\Lambda_+} \}$.  
 Matsumoto introduced in \cite{M2} the following notion of 
{\it past equivalence\/} relation.
 Two points
$x$ and $y \in  {\Lambda_+}$ are called $l$-past equivalent,  
$x \sim_l y$, if $\Lambda_l(x) = \Lambda_l(y)$. 
 The corresponding set of  equivalent 
classes is denoted by $\Omega _l :=  {\Lambda_+} / \sim_l $.  
For $\mu \in \Lambda_l$, if $x \sim_l y$, then $q_{\mu}(x) = 
q_{\mu}(y)$.  
Thus $q_{\mu}$ defines a function $\hat{q}_{\mu} \in \C(\Omega _l)$.  
The set
$\{\hat{q}_{\mu}~ \in \C(\Omega _l) ; \mu \in \Lambda_l\}$ separates the 
points 
in $\Omega _l$, thus it generates $\C(\Omega _l)$.   $A_l({\Lambda_+})$
is precisely the set of functions in $\ell ^{\infty}({\Lambda_+})$ which 
have the same value on each  $l$-past equivalent class and we have an 
isomorphism between  $ A_l({\Lambda_+})$ and $\C(\Omega _l)$.  
We see directly  that the commutative $C^*$-algebra $A_{\Lambda}$ 
 is isomorphic to $\B = A({\Lambda_+})$
via an isomorphism which identifies   $S_{\mu}^*S_{\mu}$ with $q_{\mu}$, 
 $\mu \in \Lambda^*$. 
The  $^*$-endomorphism $\gamma_i$ on $A_{\Lambda}$
defined by $\gamma_i(T) = S_i^*TS_i$ for $T \in A_{\Lambda}$
  corresponds 
to the  $^*$-endomorphism $\rho_i$ on $\B= A({\Lambda_+})$. 
We have $\rho_i(q_{\mu}) = q_{\mu i}$ and 
$\gamma_i(S_{\mu}^*S_{\mu}) = S_{\mu i}^*S_{\mu i}$
for $\mu \in \Lambda^k$.  
Consider the right  Hilbert 
 $A_{\Lambda}$-module $Y = \oplus _{i=1}^d S_iA_{\Lambda}$ and 
   the $^*$-homomorphism 
$\phi : A_{\Lambda} \rightarrow \L_{A_{\Lambda}}(Y)$ 
 given by the diagonal matrix
$\phi (a) = \text{diag} (\gamma_i(a))_i$.
Since $(S_{\mu}^*S_{\mu})S_i = S_i\gamma_i(S_{\mu}^*S_{\mu})$, the 
Hilbert
bimodules $X$ and $Y$ are isomorphic by the identification
of   $\B$ with 
$A_{\Lambda}$.
The universality of the Matsumoto algebra and the Pimsner algebra 
immediately shows that $O_{\Lambda} \cong \O_Y \cong \O_X$. 
\medskip

\noindent{\bf 4.4. Contractions of compact metric spaces}

We next dicuss an example associated with a self-similar set 
in fractal geometry.   
Let $\Omega$ be a (separable) complete metric space and 
let $\{\gamma_1,...,\gamma_d \}$
be a finite family of nonzero proper contractions of 
$\Omega$ with  Lipschitz
constants $c_i = \text{Lip}(\gamma_i) < 1$.  We assume that $d \geq 2$. 
Then there
exists a unique  nonempty compact set $K \subset \Omega$ satisfying
the (exact) invariance condition
$$
K = \cup_i \gamma_i (K).
$$
The above invariance condition shows that the compact set $K$ is 
self-similar
in a weak sense.  For example the Cantor set,  the Koch curve and the
 Sierpinski 
gasket are typical examples of self-similar sets.  We refer the reader to the book of
 Hutchinson \cite{H}  for more information on fractal geometry.
The topological dimension of $K$ is 
dominated by the Hausdorff dimension of $K$, and the Hausdorff dimension 
of $K$
is dominated by the similaritiy dimension of $K$,  which is a finite
number $D$ satisfying $\sum_i c_i^D = 1$.   Thus $K$ has a 
finite topological dimension. 

   Consider the $C^*$-algebra $\B = \C(K)$ and the canonical Hilbert 
right $\B$-module $X = \B^d$.  
For each $i$, we define an endomorphim $\phi_i$ on $\B$ by 
$$
   (\phi_i(a))(z) = a(\gamma_i(z)), \quad a \in \B, \quad z \in K.
$$
Left $\B$-action $\phi : \B \rightarrow \L_{\B}(X)$ 
is defined by the diagonal matrix
$\phi (a) = \text{diag} (\phi_i(a))_i$.  
We see that the (exact) invariance condition
 $K = \cup_i \gamma_i (K)$ is equivalent to the fact that  
$\phi$ is injective.  The bimodule $X$ generates the Pimsner 
$C^*$-algebra
$\O_X$.  Let $\{x_1, \dots, x_n \}$ be the canonical basis of $X$.  
Then the corresponding elements $\{S_1,\dots, S_d\}$ of $\O_X$ generate a
copy of the Cuntz
algebra
$\O_d\subset\O_X$.  By construction, the Pimsner $C^*$-algebra 
$\O_X$ is isomorphic to the universal $C^*$-algebra generated by 
$\B = \C(K)$ and $\O_d$ satisfying the relations 
$aS_i = S_i\phi_i(a)$  for $a \in \B$ and $i = 1,\dots,d$.

In \cite{H} Hutchinson shows that there exists a unique regular Borel
probability measure $\mu$ on $K$ satisfying, 
for any measurable set $F$,
$$
   \mu (F) = \sum _{i=1}^d \frac{1}{d} \mu (\gamma _i^{-1} (F)).
$$
Consider the trace $\tau_0$
on $\B$ corresponding to the  probability measure $\mu$.
 Then  $\tau_0$
satisfies 
$$
    \tau_0 (\sum _i <x_i,ax_i>) = d\tau_0(a), \quad    a \in \B.
$$      
Since $<x_i,x_i> = 1$ for $i = 1,\dots ,d$  , $\O_X$ has a KMS
state at the inverse temperature $\beta$ if and only if
$\beta = \log d$.  Moreover the uniqueness of the  probability
measure implies that the corresponding KMS state is also unique.
  We shall show that the algebra $\O_X$ is in fact the Cuntz algebra 
$\O_d$.   Before proving this, we study  a  more general situation
to include  standard $d$-times around embeddings.

Let $K$ be a compact metric space.  Consider the $C^*$-algebra $\B = 
\C(K)$ and  the state space $\S$ of $\B$.  
Let $\text{Lip}(K)$ be the space of Lipschitz functions, and let 
 $\text{Lip}(f)$ denote the Lipschitz constant of $f\in \text{Lip}(K)$.  
   In \cite{H} Hutchinson considers the following metric $L$ on $\S$:
$$
  L(\varphi_1, \varphi_2) = \sup\{|\varphi_1(f) -  \varphi_2(f)| ;
 f \in \text{Lip}(K), \text{Lip}(f) \leq 1  \}.
$$
Then $(\S, L)$ is a complete metric space, and the topology
defined by $L$ is precisely   the weak$ ^*$-topology of $\S$.

Consider the canonical Hilbert right $\B$-module $X = \B^d$     
and any injective  unital $^*$-homomorphism 
$\phi : \B \rightarrow \L_{\B}(X)$.  We identify  $\phi (a)$
with the matrix $(\phi_{ij} (a))_{ij} \in \B \otimes M_d({\Bbb C})$
  for 
$a \in \B$. 
   
Then the  Pimsner $C^*$-algebra 
$\O_X$ is isomorphic to the universal $C^*$-algebra generated by 
$\B = \C(K)$ and the Cuntz algebra $\O_d$ satisfying the relations 
$aS_j = \sum _i  S_i\phi_{ij}(a)$  for $a \in \B$ and $i,j = 1,\dots,d$.       

\medskip

 \noindent{\bf 4.5. Proposition} {\sl  In the above situation, let 
$\Psi$ be the unital  positive map on 
 $\B = \C(K)$ defined by 
$\Psi (a) = \frac{1}{d} \sum_i \phi_{ii}(a)$.
 If the Banach 
space adjoint $\Psi^*$  induces a proper contraction on $\S$ with 
respect to 
the metric $L$, then $\O_X^0$ has a unique
tracial state.  Moreover  $\O_X$ has a KMS
state at the inverse temperature $\beta$ if and only if
$\beta = \log d$ and the corresponding KMS state is also unique.}
\medskip

\noindent{\it Proof}
We identify $\L(X^n)$ with $\B \otimes M_{d^n}({\Bbb C})$.  Using the 
commutation relation, it is easy to see that  the 
inclusion map $\Phi_n : \L(X^n) \rightarrow \L(X^{n+1})$ 
is described by  matrices 
$$
\Phi_n((a_{\alpha, \beta})_{\alpha, \beta})
    = (\phi_{ij}(a_{\alpha, \beta}))_{(\alpha,i),(\beta,j)},
$$
where  $ \alpha$ and $\beta$ run the set $\{1,\dots ,d\}^n$ of words 
with
length $n$. Thus $ \Phi_n = \phi \otimes id$ on 
$\B \otimes M_{d^n}({\Bbb C})$.  
 A tracial state on ${\O_X}^0 = \varinjlim _n \L(X^n)$
is described by a sequence of tracial states $\{\tau_n\}$, 
$n\geq0$ on $\L(X^n)$ such that 
$\tau_{n+1}\upharpoonright_{\L(X^n)}=\tau_n$. 
Therefore one needs to assign a 
sequence of 
states $(\varphi_n) \in \S$ with 
$\tau_n = \varphi_n  \otimes tr$ on $
\L(X^n) \cong  \C(K) \otimes  M_{d^n}({\Bbb C})$.  
The coherence relations require that 
$$\Psi^*(\varphi_{n+1})=\varphi_n,\quad n\geq0.$$
Since $K$ is compact, the diameter of $(\S, L)$ is bounded.  Hence
$\cap_{n=1}^{\infty}  \Psi^{*n} (\S)$ 
consists of a single point $\omega_0$.
The constant sequence of
states $(\varphi_n) \in \S$ with 
$\varphi_n =  \omega_0$ gives a tracial state $\tau$ on ${\O_X}^0$,
because $\omega_0$ is the unique fixed point of $\Psi^*$ in $\S$.
Choose another tracial state.  
The coherence relations  
$\Psi^*(\varphi_{n+1})=\varphi_n$ shows that any 
$\varphi_n$ belongs to $\cap_{r=1}^{\infty}  \Psi^{*r} (\S)$.
Therefore $\varphi_n = \omega_0$.  Thus  $\O_X$ has a unique KMS state 
at 
inverse temperature $\log d$.    
\medskip

We remark that the present situation is  similar to that of  the 
Cuntz-Krieger algebras 
associated to aperiodic matrices.
In fact,  for any  state $\omega \in \S$,  $\Psi^{*n} (\omega)$ 
converges
to the unique $\omega_0$  in $\S$ with respect to $L$.  This resembles 
 the  Perron--Frobenius
Theorem for aperiodic matrices.

\medskip

\noindent{\bf 4.6. Example}  In the fractal case,  each proper 
contraction
$\gamma_i$ induces an endomorphism $\phi_i$ on  $\B = \C(K)$ satisfying 
$\text{Lip}(\phi_i(f)) \leq c_i\text{Lip}(f)$ for $f\in \text{Lip}(K)$.
Hence $\Psi^* = \frac{1}{d} \sum_i \phi_i^*$ 
is a proper contraction on $\S$ with respect to 
the metric $L$.  
\medskip

\noindent{\bf 4.7. Example}   We next study the example  of  standard
$d$-times around embeddings.  
 Let  $\B = \C({\Bbb T})$ with 
 ${\Bbb T} = {\Bbb R}/ {\Bbb Z}$ 
  and $X = \B^d$ be  the natural right Hilbert 
 $\B$-module.  Then a standard $d$-times around embedding
is defined by a map
 $\phi : \C({\Bbb T}) \rightarrow \C({\Bbb T}, M_d({\Bbb C}))$ 
of the form
 $$
 (\phi(f))(t) = u_t \text{diag} (f(\frac{t}{d}), f(\frac{t+1}{d}), 
 \dots, f(\frac{t+d-1}{d}))u_t^*,
 $$
 where $(u_t)_t$ is a continuous path of unitaries in $M_d({\Bbb C})$ 
 such that $u_0 = I$ and such that $u_1$ is the unitary matrix
corresponding 
to the operator taking vectors  $e_1,\dots,e_d$ of the canonical
basis of ${\Bbb C}^d$ to $e_d, e_1,\dots, e_{d-1}$ respectively.
 We regard $\phi$ as a map 
 $\phi : \B \rightarrow \L_{\B}(X) = \B \otimes M_d({\Bbb C})$.
Then $\Psi^*=\frac{1}{d}\sum_i{\phi_i}^*$  induces a proper contraction
on $\S$ with respect to 
the metric $L$.  
In fact 
$$
(\Psi (f))(t) = \frac{1}{d} (f(\frac{t}{d}) + f(\frac{t+1}{d}) +  
 \dots + f(\frac{t+d-1}{d})).
 $$
We assume that $d = 2$ for the simplicity of  notation.  For 
 $f\in \text{Lip}({\Bbb T})$ and $x, y \in {\Bbb T}$, choosing 
carefully the 
 nearest pairs between  $\{ \frac{x}{2}, \frac{x+1}{2} \}$ and 
  $\{ \frac{y}{2}, \frac{y+1}{2} \}$, we have 
  $|(\Psi (f))(x) - (\Psi (f))(y)| 
  \leq \frac{1}{2}\text{Lip}(f)d(x,y)$.  
Hence  $\text{Lip}(\Psi (f)) \leq \frac{1}{2} \text{Lip}(f)$.  
Therefore for $\varphi_1, \varphi_2 \in \S$,  we have 
$$
 |(\Psi^* (\varphi_1))(f) - (\Psi^* (\varphi_2))(f)| 
 = |\varphi_1(\Psi(f)) -\varphi_2(\Psi(f))| 
 \leq L(\varphi_1, \varphi_2)\frac{1}{2} \text{Lip}(f).
$$
Thus  $L(\Psi^* (\varphi_1), \Psi^* (\varphi_2)) 
\leq \frac{1}{2} L(\varphi_1, \varphi_2).
$
\medskip
 
\noindent{\it 4.8. Remark}
One can easily show, using known results, that under suitable
circumstances $\O_X$ is simple and purely infinite. Indeed, assume that
 $K$ is totally disconected
or connected  
and
has  a finite  topological dimension. If $\O_X^0 = \varinjlim _n 
\L(X^n)$
is simple, then $\O_X^0$ is of real rank zero by \cite{BDR},  since  
$\O_X^0$ has a unique trace.    
By a result of
Martin and Pasnicu
\cite{MP}, $\O_X^0$
has the comparability property on every matrix algebra.
Thus we can  apply  a result by \cite{GP} (Theorem 1.6)
 and conclude that the  $\O_X$ 
is simple and purely infinite.  For example, in the case of a standard 
$d$-times around embeddings all the assumptions are
satisfied.  In fact 
$\O_X^0 = \varinjlim _n \L(X^n)$ is a Bunce-Deddens algebra. 
Note that we have naturally embedded an AT algebra into a purely infinite
simple $C^*$--algebra.       
Again, in the fractal case , it is easy to show that $\O_X^0$ is simple.  
So we can  
apply the preceding argument.  However,  it is not difficult 
to show that in this case 
 $O_X$ is canonically isomorphic to the Cuntz algebra $\O_d$
(a fact which will be later generalized to noncommutative metric spaces):  
We identify $\L(X^n)$ with $\C(K, M_{d^n}({\Bbb C}))$.  Then the 
inclusion map $\Phi_n : \L(X^n) \rightarrow \L(X^{n+1})$ 
is described by  block diagonal matrices 
$$
(\Phi_n(f))(t) = \text{diag} (f(\gamma_1(t)), \dots , 
f(\gamma_d(t))).
$$
Let $\omega = (\omega_1,\dots, \omega_k) \in 
\{1, \dots , d\}^k$ be a finite word and 
$\gamma_{\omega} = \gamma_{\omega_1} \dots\gamma_{\omega_k}$.
Then the inclusion map $\Phi_{n+k,n}$ of $\L(X^n)$
into $\L(X^{n+k})$ is given by   
$$
   (\Phi_{n+k,n}(f))(t) =  \text{diag}  
(f(\gamma_{\omega}(t)))_{\omega}.
$$ 
By the uniform continuity of $f$,  $\Phi_{n+k,n}(f)$ is approximated 
by
a constant matrix  up to $\varepsilon$ for a sufficient large $k$.  
Thus $\O_X^0$ is a UHF algebra $M_{d^{\infty}}$ and $O_X$ is exactly the 
 Cuntz algebra $\O_d$ generated by the original operators  
 $\{S_1, \dots , S_d\}$. 
\medskip

\noindent{\bf 4.9. Contractions of noncommutative metric spaces}

The preceding argument suggests  a generalization to
noncommutative metric spaces introduced by Connes in \cite{Co}.  
Our setting will be the following. Let $\A$ and $\B$ be unital
$C^*$-algebras.
Suppose that $\B$ is a Banach bimodule over $\A$.  
Let $\delta : \A \supset \text{Dom}(\delta) \rightarrow \B$ be a densely 
defined $^*$-derivation with $\text{ker}\delta = {\Bbb C}I$. 
 Let $\S$ be the state space of $\A$.  
 Consider the following metric $L$ on $\S$:
$$
  L(\varphi_1, \varphi_2) = \sup\{|\varphi_1(a) -  \varphi_2(a)| \  ; \ 
 a \in \text{Dom} (\delta) , \| \delta (a) \| \leq 1  \}.
$$  
The metric is allowed to take the value $\infty$.

In \cite{RiI},  Rieffel considers  the question of whether  the
metric topology agrees with the underlying weak$^*$ topology on the state 
space.  His setting  is, however, more
general, as he works with normed vector  spaces endowed with seminorms not 
necessarily arising from $^*$--derivations.

 We assume that 
$$\{ a \in \text{Dom}(\delta) \ ; \  \| \delta (a) \| \leq 1 \}/ {\Bbb
C}I\quad\text{is bounded in\ }\A/{\Bbb C}I.\eqno(4.1)$$  
By Proposition 1.6 in \cite{RiI} this condition is equivalent to the fact
that 
the metric $L$ on $\S$ is bounded.

Let $\{\phi_1, \dots , \phi_d\}$ be a finite family of unital 
$^*$-endomorphisms on $\A$, with $d \geq 2$.  Recall that 
the  crossed product $C^*$-algebra 
$C^*(\A;\phi_1, \dots , \phi_d)$ of $\A$ by $\{\phi_1,\dots,\phi_d\}$
 is the universal 
$C^*$-algebra generated by the image of a $C^*$--homomorphism
$\pi : \A \rightarrow C^*(\A;\phi_1, \dots , \phi_d)$
 and the Cuntz algebra 
$\O_d$ with the generators $S_1, \dots, S_d$ 
satisfying the relations 
$\pi(a)S_i = S_i\pi(\phi_i(a))$  for $a \in \A$ and $i = 1,\dots,d$.  
We note that $\pi$ is isometric if and only if 
$\cap _i \text{ker} \phi_i = 0$.  In this case the  crossed 
product
 $C^*$-algebra $C^*(\A;\phi_1, \dots , \phi_d)$ is isomorphic to $\O_X$,
where $X$ is 
the trivial Hilbert right $\A$-module $X = \A^d$    
endowed with the diagonal left $\A$--action:
 $\phi : \A \rightarrow \L_{\A}(X)$ 
$\phi (a) = \text{diag} (\phi_1(a),\dots,\phi_d(a))$. 
\medskip

\noindent{\bf 4.10. Proposition} {\sl In the above setting, 
 assume that the restrictions $\gamma _i$ of the Banach space adjoint 
$\phi _i^*:\A^*\to\A^*$  to 
the state space $\S$ of $\A$ are proper contractions with respect to $L$.
Then the endomorphism crossed product $C^*$-algebra 
$C^*(\A;\phi_1, \dots , \phi_d)$  is canonically isomorphic to the Cuntz
algebra  
$\O_d$ and has  a unique KMS state
at inverse temperature $\log d$. }  
\medskip

\noindent{\it Proof} Let $c$ be the maximum of the Lipschitz norms 
$c_i = \text{Lip} (\gamma_i)$,
$i = 1, \dots , d$.
For any $a \in \text{Dom} (\delta)$ and $\varphi$ , $\psi \in \S$, we 
have 
$$
 | \varphi (a) - \psi (a) | \leq  L(\varphi , \psi)\| \delta (a) \|.
$$
For a finite word $\alpha = (\alpha_1, \dots, \alpha_n) \in \{1, \dots ,
d\}^n$, 
we use the multi-index notation $\phi_{\alpha}$ and
$\gamma_{\alpha}$.
Then we have that 
$$
L(\gamma_{\alpha}(\varphi), \gamma_{\alpha}(\psi)) \leq c^n L(\varphi, 
\psi) 
\leq c^n \text{diam} (\S,L),
$$
for any pair of states $\varphi$ , $\psi \in \S$.

We shall show that $\pi(A)$ is included in  the canonical UHF subalgebra 
 $M_{d^{\infty}}$ of the Cuntz algebra $\O_d$. 
   For any $a \in \text{Dom} (\delta)$ and 
 $\varepsilon > 0$, there exists $n\in{\Bbb N}$ such that
 $$
 c^n \text{diam} (\S,L) \| \delta (a) \| \leq \varepsilon.
$$
Fix a state $\omega_0 \in \S$.  
Consider the diagonal matrix $t = \text{diag} (\omega_0(\phi_{\alpha}(a)))
_{\alpha}
 \in M_{d^n}({\Bbb C})$.  Then for any state $\omega \in \S$, we have
 $$
 | \omega (t_{\alpha \alpha} - \phi_{\alpha}(a) )|
 = | \omega_0(\phi_{\alpha}(a)) - \omega ( \phi_{\alpha}(a)) |
 \leq  L(\gamma_{\alpha}(\omega_0), 
\gamma_{\alpha}(\omega))\|\delta(a)\| 
 \leq \varepsilon .
 $$
 Hence $\| t_{\alpha \alpha} - \phi_{\alpha}(a) \| \leq 2\varepsilon
 $.  
 Since 
 $$
 \pi(a) = \pi(a)\sum_{\alpha}S_{\alpha}S_{\alpha}^* 
          = \sum_{\alpha} S_{\alpha}\pi(\phi_{\alpha}(a))S_{\alpha}^*,
$$          
we  have that 
$$
\| \sum_{\alpha} t_{\alpha \alpha} S_{\alpha}S_{\alpha}^* - \pi(a) \|
= \| \sum_{\alpha} S_{\alpha}\pi(t_{\alpha \alpha}-\phi_{\alpha}(a))
S_{\alpha}^* \|
\leq 2\varepsilon.
$$
Thus   $C^*(\A;\phi_1, \dots , \phi_d)$  is 
precisely the Cuntz algebra $\O_d$ generated by the original  
 $\{S_1, \dots , S_d\}$.
\medskip

\noindent{\it Remark} M. Rieffel has kindly pointed out to us that
in the proof
of the previous Proposition
we never use the fact that $L$ comes from a $^*$--derivation, 
but rather only that it is a seminorm satisfying the boundedness condition
$(4.1)$. Seminorms of this kind were studied in \cite{RiII}.
\medskip

\noindent{\bf 4.11. Example}
  Let $\D$ be a noncommutative  
unital 
$C^*$-algebra.  Let $\A = \{ a \in \C([0,1],\D) ; a(0) \in {\Bbb C}I \  
\}$.  
Set $Y = \{(x,y) \in [0,1] \times [0,1] ; x \not= y \}.$  Let 
$\B = C^b(Y,\D)$ be the set of $\D$-valued bounded continuous functions 
on $Y$.  Then $\B$ is a Banach bimodule over $\A$ by 
$$
(a_1fa_2)(x,y) = a_1(x)f(x,y)a_2(y).
$$
  
Let $\delta : \A \supset \text{Dom}(\delta) \rightarrow \B$ be the densely 
defined $^*$-derivation of De Leeuw, given by 
$$
(\delta (a))(x,y) = \frac{a(x)-a(y)}{|x-y|},
$$ 
where $\text{Dom}(\delta)$ is the set of Lipschitz functions in
$\A$.  Then
 $\text{ker}\delta = {\Bbb C}I$.  Let $\alpha$ be the ${}^*$-endomorphism
on 
 $\A$ defined by $\alpha (f)(x) = f(\frac {x}{2})$.  
 Then the restriction $\gamma$ of the Banach space adjoint 
$\alpha^*:\A^*\to\A^*$  to 
the state space $\S$ is a  proper contraction with respect to $L$.
Therefore  the endomorphism crossed product $C^*$-algebra 
$C^*(\A;\alpha, \alpha)$  is isomorphic to the Cuntz algebra  
$\O_2$.

\medskip

\heading
5. KMS states of Pimsner $C^*$--algebras \\ associated to  Cuntz--Krieger
bimodules
\endheading

In this section we illustrate Theorem 2.5 by examples. We shall discuss
some
situations where there is a unique  
KMS state, or more generally,
where the set of KMS states can be easily characterized. The inspiring 
example is that of the Cuntz--Krieger
algebras that we discuss here below.  Some of the following facts 
 are well known in terms of path algebras 
in subfactor theory.
\medskip

\noindent{\bf 5.1. KMS states of Cuntz--Krieger algebras}

Let $\O_A$ be the Cuntz--Krieger algebras associated to a 
 matrix 
$A=(a_{i j})\in M_d(\{0,1\})$.
A KMS state $\omega$ for the canonical 
circle action restricts to a tracial state $\tau$ on the f.d. 
commutative subalgebra
$\A$ generated by the ranges $P_1,\dots, P_d$ of the generating partial
isometries
$S_1,\dots,S_d$. Let
$\lambda=(\lambda_1,\dots,\lambda_d)$ $\in{\Bbb R^+}^n$ be defined by 
$\lambda_i=\omega(P_i)$.
 Since $\omega$ is normalized, we have   $\sum \lambda_i=1$.
The scaling property
$s'(\tau)=\epsilon(\tau)\tau$  says, when checked on $\A$, that
$\lambda$ is a nonnegative eigenvector of $A$, and hence,
 when $A$ is irreducible,
 it is the unique normalized
Perron  eigenvector for $A$ by the Perron--Frobenius Theorem \cite{G}.
Let us analyse more in detail the structure of the Banach space 
$\T({\O_A}^0)$ and the spectrum of the operator $s'$.
Let $L_r$, $r\geq1$, denote the unital finite--dimensional 
$C^*$--subalgebra of $\O_A$
generated by elements of the form $S_{i_1}\dots S_{i_r}P_k(S_{j_1}\dots 
S_{j_r})^*$.
Set $L_0={\Bbb C}P_1+\dots+{\Bbb C}P_d$. Then the set of 
minimal central projections for $L_r$ is
$\{\sigma^r(P_k), k=1,\dots,d\}$, where $\sigma$ is the canonical 
endomorphism of ${L_0}'\cap\O_A$
defined by $\sigma(T)=\sum_i S_iT{S_i}^*$. 
So $L_r\sigma^r(P_k)\simeq M_{d_{r,k}}({\Bbb C})$
where 
$$d_{r, k}=\text{Card}\{(i_1,\dots,i_r): S_{i_1}\dots
S_{i_r}P_k\neq0\}=$$
$$\text{Card}\{(i_1,\dots,i_r): 
a_{i_1,i_2}a_{i_2, i_3}\dots a_{i_r, k}\neq0\}.$$
Therefore $d_{r,k}$ is the  sum of the entries in the 
$k$--th column of 
 $A^r$:
$$d_{r,k}=\sum_{i_1,\dots,i_r}a_{i_1,i_2}a_{i_2, i_3}\dots a_{i_r, k}.$$
The $(j,i)$--entry of inclusion matrix 
of  $L_r\subset L_{r+1}$ can be computed by looking at
the projection $\sigma^r(P_i)\sigma^{r+1}(P_j)=\sigma^r(S_iP_j{S_i}^*)$
 which $0$ when $a_{i,j}=0$,
otherwise it is the sum $d_{r,i}$ minimal projections of 
$L_{r+1}\sigma^{r+1}(P_j)$. Thus 
the inclusion matrix of $L_r\subset L_{r+1}$ is $A^t$. A tracial state 
on ${\O_A}^0$
is described by a sequence of positive traces $\{\tau_r\}$, 
$r\geq0$ on $L_r$ such that $\tau_0$
is normalized and
$\tau_{r+1}\upharpoonright_{L_r}=\tau_r$. Therefore one needs to assign 
a 
sequence of 
nonnegative column vectors $(t_r)\in{\Bbb R_+}^d$
 which will be the values that $\tau_r$ takes on the minimal
projections of $L_r$. The coherence relations require that 
$$A(t_{r+1})=t_r,\quad r\geq0$$
while the positivity and normalization properties translates into:
$$t_r(k)\geq0,\quad r\geq0, k=1,\dots,d,$$
$$\sum_k t_0(k)=1.$$
The latter  implies, as expected, normalization of each $\tau_r$:
$$\sum_k t_r(k)d_{r,k}=\sum_i (A^rt_r)(i)=\sum_it_0(i)=1, r\geq0.$$
Removing  positivity and normalization, but requiring
instead a norm bound for the sequence $(t_r)$, one finds that
$\T({\O_A}^0)$ is
described by 
$$\{(t_r)_{r\geq0}: t_r\in{\Bbb C}^d, t_r=At_{r+1}, r\geq0, 
\sup_r\sum_i (A^r|t_r|)(i)<\infty\},$$
with the Banach space norm 
$$\|(t_r)_{r\geq0}\|=\sup_r\sum_i (A^r|t_r|)(i).$$
The operator $s'$ acts as:
$$s'((t_r)_{r\geq0})=(At_0, t_0, t_1,\dots)$$
while its inverse is
$$t'((t_r)_{r\geq0})=(t_1,t_2,\dots).$$\medskip

\noindent{\bf 5.2. Proposition} {\sl If $A=(a_{i,j})\in M_d(\{0,1\})$ 
is an irreducible symmetric  matrix then 
$\T({\O_A}^0)$ is linearly spanned by elements of the form
$(\lambda^{-r}t_0)_{r\geq0}$ where $t_0$ is an eigenvector of $A$ with
eigenvalue $\lambda$ and $|\lambda|=r(A).$ Furthermore
$$\sigma(s')=\{\lambda\in\sigma(A): |\lambda|=r(A)\}.$$}
\medskip

\noindent{\it Proof}
Let $t_0$ be an eigenvector for
$A$ with eigenvalue $\lambda$ such that $|\lambda|=r(A)$. 
Set $t_r:=\lambda^{-r}t_0$, so that $At_{r+1}=t_r$. We have:
$$\|A^r|t_r|\|_2=r(A)^{-r}\|A^r|t_0|\|_2\to\|E_0|t_0|\|_2,$$ 
where $E_0$ is the rank one orthogonal projection onto the span of the
Perron eigenvector. It follows that $(t_r)\in \T({\O_A}^0)$.
Furthermore $(t_r)$ is also an eigenvector of $s'$ with the same 
eigenvalue.
The same argument shows that if $t_0\in \T({\O_A}^0)$ were
an eigenvector with eigenvalue
$\lambda$ such that $|\lambda|<r(A)$, and $t_r$ is defined as above,
 then $\|A^r|t_r|\|_2$ would be unbounded, so that
$(t_r)$ does not define an element of $\T({\O_A}^0)$.
 With similar arguments 
one sees that $\T({\O_A}^0)$ is linearly spanned by vectors of the form
$(\lambda^{-r}t_0)$ where $|\lambda|=r(A).$ 
We  show that if either $\lambda\notin\sigma(A)$ or
$\lambda\in\sigma(A)$ but $|\lambda|<r(A)$
then  
 $s'-\lambda $ is invertible. 
We start assuming that $\lambda\notin\sigma(A)$.
Then $s'-\lambda$  is clearly injective.
We show that it is also surjective. Given $(v_r)\in \T({\O_A}^0)$
set $t_r=(A-\lambda)^{-1}v_r$. Then 
$$At_r=A(A-\lambda)^{-1}v_r=(A-\lambda)^{-1}Av_r=
(A-\lambda)^{-1}v_{r-1}=t_{r-1}.$$
For any matrix $B$ with complex entries let $B^+$ stand for the 
matrix with entries the absolute values of the corresponding elements
of $B$, and let $M_B$ be the 
maximum
of the absolute values of its entries.
Then
$$A^r|t_r|\leq A^r{(A-\lambda)^{-1}}^+|v_r|\leq 
M_{(A-\lambda)^{-1}}A^r|v_r|$$
which shows that $(t_r)\in \T({\O_A}^0)$ and $(s'-\lambda)(t_r)=v_r$.
 Assume now that $\lambda\in\sigma(A)-\{0\}$ but $|\lambda|<r(A)$. 
We have already noted that $s'-\lambda$ is injective. 
Furthermore since for each $(t_r)\in \T({\O_A}^0)$, any $t_r$ 
belongs to the range of
$A-\lambda$, with similar arguments one shows that $s'-\lambda$
is surjective.
\medskip

\noindent{\it  Remark}
If $A$ is aperiodic, then the homogeneous subalgebra 
of $\O_A$ has a unique trace. More generally,
 if one drops the assumption that $A$ is symmetric,
then, with a more extensive use of Perron--Frobenius theory, one can 
still
show that eigenvectors corresponding to nonmaximal eigenvalues  do
not appear in the point spectrum of  $s'$, hence our result shows
that $\sigma(s')\subset\{\lambda: |\lambda|=r(A)\}.$ 
Note also that if we more generally start with a reducible matrix $A$,
then
we are in a situation of  nonuniqueness of 
KMS states for $\O_A$ (corresponding to the minimal and maximal
Perron eigenvalues of $A$).\medskip

\noindent{\bf 5.3. KMS states as Markov traces arising from inclusions
of finite algebras with finite Jones index}

Let $N\subset M$ be an inclusion of 
$II_1$--factors with finite index or of 
finite--dimensional $C^*$--algebras such that $Z(N)\cap Z(M)={\Bbb C}$.
Let, in the latter case,  $A$ be the inclusion matrix. 
 Let $\tau$ be a faithful tracial state
on $M$, and consider the unique $\tau$--preserving conditional 
expectation
$$E_\tau: M\to N.$$ Endow $X=M$ with the $C^*$--bimodule structure
over $N$
as follows. The structure of $N$--bimodule is 
defined by left and right multiplication, while the $N$--valued 
inner product is 
$$<x,y>_N:=E_\tau(x^*y).$$
Then $X$ is full and  finite projective as a right  $N$--module 
(\cite{GHJ}). 
It is not difficult to check that 
$$\L_N(X^r)=M_{2r-1},$$
where $$M_{-1}=N\subset M_0=M\subset M_1\subset\dots$$ 
is the the Jones tower. 
A KMS state at inverse temperature $\beta$
for $\O_X$ corresponds precisely 
to a Markov trace for the tower, which is unique, and one has
$$\beta=\log([M:N]).$$
See \cite{K}.   If $N\subset M$ are finite factors, each term of the 
tower is a finite factor,
hence its trace space is one dimensional, and spanned by the Markov 
trace.
So $\text{dim}\T({\O_X}^0)=1$, and $s'$ acts multiplying by $[M:N]$. 
If $N\subset M$ are finite--dimensional $C^*$--algebras, 
the inclusion matrix of
$\L_N(X^r)\subset \L_N(X^{r+1})$ is $A^tA$, which 
is symmetric and irreducible (\cite{GHJ}). It is not difficult to show, 
with arguments similar
to those of the previous example, that $\T({\O_X}^0)$ is again linearly 
spanned
by traces corresponding to eigenvectors of $A^tA$ with 
maximal eigenvalue. \medskip

\noindent{\bf 5.4. KMS states of Pimsner algebras associated with
Cuntz--Krieger bimodules}

After these motivating examples,
we  consider, more generally,   systems  of the form $(\O_X,\gamma)$, 
where
  $X$ is what we call a Cuntz--Krieger Hilbert $C^*$--bimodule and
$\gamma$
is the canonical gauge 
action. Such Hilbert bimodules, and simplicity of the corresponding
$C^*$--algebras $\O_X$, have been considered in \cite{KPW}.
Consider $d\geq2$ unital simple $C^*$--algebras 
$\A_1,\dots\A_d$ 
and a matrix $A=(a_{i,j})\in M_d(\{0, 1\})$ with no row and no column 
identically zero.
Let, for any pair of indices $i,j$ such that $a_{i,j}=1$, $X_{i,j}$ be a 
full, finite
projective  $\A_i$--$\A_j$ Hilbert bimodule, and let 
$X=\oplus_{i,j:a_{i,j}=1} X_{i,j}$ be endowed
with the natural structure of Hilbert bimodule over 
$\A:=\A_1\oplus\dots\oplus\A_d$. Then since no row of $A$
is zero, left $\A$--action is faithful, and since no column of $A$ is 
zero, $X$ is full.
Clearly $X$ is finite projective as a right module. 
We assume that there is a system of tracial states $\tau_1,\dots,\tau_n$ 
on
 $\A_1,\dots\A_n$ respectively
satisfying, for each pair of indices for which $a_{j,k}=1$, 
$$\tau_k(\sum_r<{x_{r}}^{k,j} 
a{x_{r}}^{k,j}>)=\lambda_{j,k}\tau_j(a),\quad a\in \A_j,$$
for some $\lambda_{j,k}>0$. Here $\{{x_{r}}^{k,j}\}_r$ if a basis of 
$X_{j,k}$.
We set $\lambda_{j,k}=0$ if $a_{j,k}=0.$ 
We will call $\{\tau_i\}$ a coherent set of traces. Note that 
$(\lambda_{j,k})$ is irreducible precisely
when $A$ is. If each $\A_j$ has a unique tracial state $\tau_j,$
the set $\{\tau_j\}$ 
is coherent. This is indeed the case of   Cuntz--Krieger algebras.
Let $P_j$ be the identity of $\A_j$. Recall from \cite{KPW} that
$\A'\cap\O_X$ has a unique unital endomorphism $\sigma$ such that 
$$\sigma(a)x=xa,\quad x\in X, a\in\A'\cap\O_X.$$ 
Recall also that if $A$ is irreducible, $\O_X$ is simple \cite{KPW}, and
if $A$ is aperiodic, $\A$ is separable and has real rank zero and all 
$M_n(\A_j)$ have the comparability
property, then  $\O_X$ is simple and purely infinite (cf. Theorem
1.6). For all $r\geq1$ 
$\L_\A(X^r)$ is a finite
direct sum of unital simple $C^*$--algebras,
and its minimal central projections are 
$\sigma^r(P_1),\dots\sigma^r(P_d).$
One has $\L_\A(X^r)\sigma^r(P_j)=\L_{\A_j}(X^rP_j),$ where $X^rP_j$ is 
regarded as an
$\A$--$\A_j$ Hilbert bimodule.\medskip

Let $\tau_1,\dots,\tau_d$ be a coherent choice of tracial states on
 $\A_1,\dots,\A_d$ respectively.
Let
 $\{{u_i}^{(r),j}\}$ be a basis of $X^rP_j$. 
Then the positive functional 
$$a\in\L_{\A_j}(X^rP_j)\to
\tau_j(\sum_i<{u_i}^{(r),j}, a{u_i}^{(r),j}>)$$ is nonzero, tracial and 
independent of the basis. Let
$\epsilon_r(\tau_j)$ denote its norm. After
 normalization, we get a
tracial state ${T_j}^{(r)}$ on $\L(X^rP_j).$ If $\A_j$ has a unique 
tracial state,
${T_j}^{(r)}$ is the unique tracial state of $\L(X^rP_j).$
Consider a tracial state $\tau_r$ on $\L(X^r)$ which restricts to 
a multiple 
of ${T_j}^{(r)}$ on $\L(X^rP_j)$, and let 
${t_j}^{(r)}=\tau(\sigma_r(P_j)),$
so ${t_j}^{(r)}\geq0$ and $\sum_j{t_j}^{(r)}=1.$ Then $\tau_{r}$
restricts to $\tau_{r-1}$ if and only if
for all $a\in\L(X^{r-1}P_j)$, and all $j$,
$${t_j}^{(r-1)}{T_j}^{(r-1)}(a)=\sum_k{t_k}^{(r)} 
{T_k}^{(r)}(a\sigma^r(P_k)).$$
Now if $\{x_j\}$ is a basis of $X$, $\{x_{j_1}\dots x_{j_r}P_k\}$ is a 
basis of $X^rP_k$, so
$${T_k}^{(r)}(a\sigma^r(P_k))=
\epsilon_r(\tau_k)^{-1}\tau_k(\sum_{j_1 \dots j_r} 
P_k{x_{j_r}^*}(x_{j_1}\dots x_{j_{r-1}})^*ax_{j_1}\dots 
x_{j_{r-1}}x_{j_r}P_k)=$$
$$\lambda_{j,k} \epsilon_r(\tau_k)^{-1}\epsilon_{r-1}(\tau_j) 
{T_j}^{(r-1)}(a).$$
So $\tau_r$ restricts to $\tau_{r-1}$ if and only if
$$\sum_k \lambda_{j,k}{t_k}^{(r)} \epsilon_r(\tau_k)^{-1}=
{t_j}^{(r-1)}\epsilon_{r-1}(\tau_j)^{-1}.$$
Set ${v_k}^{(r)}:=\epsilon_r(\tau_k)^{-1}{t_k}^{(r)}.$ Then a  solution 
is obtained
choosing for $v^0$ the Perron--Frobenius eigenvector of  the nonnegative 
matrix $(\lambda_{j,k})$
with the normalization $\sum{v_k}^0=1$, and iteratively 
$v^r=\lambda^{-1}v^{r-1}$,
where $\lambda$ is a positive eigenvalue of $(\lambda_{j,k}).$
Note that if $A$ is aperiodic, then $(\lambda_{j,k})$ is aperiodic
as well, so such a $v^0$ is the only possible
solution.
One can easily check that the tracial  state $\tau$ thus obtained on 
${\O_X}^0$ satisfies
$$\tau(\sum_j{x_j}^*ax_j)=\lambda\tau(a),\quad a\in{\O_X}^0,$$
and therefore gives rise to a KMS state of $\O_X$.

Summarizing, we have proved the following result.
\medskip

\noindent{\bf 5.5. Theorem} {\sl Let $A=(a_{i,j})\in M_d(0,1)$ be an 
irreducible matrix, $\A_1,\dots,\A_d$
unital simple $C^*$--algebras with a nonzero trace, and let,
for any pair of indices for which $a_{i,j}=1$, $X_{i,j}$ be a full, 
finite projective,
Hilbert $\A_i$--$\A_j$ bimodule. Consider $X=\oplus_{i,j: a_{i,j}=1} 
X_{i,j}$ as a 
Hilbert bimodule over $\A_1\oplus\dots\oplus\A_d$. Then 
 any normalized Perron--Frobenius eigenvector 
of the irreducible  matrix $(\lambda_{i,j})$
associated to a coherent system of tracial states on 
$\A_1,\dots,\A_d$ defines, as described above, a KMS state
of $\O_X$ at inverse temperature $\beta=\log r((\lambda_{i,j}))$. 
In particular, if
each $\A_j$ has a unique tracial state, then there is
a unique system of coherent tracial states, and the corresponding
KMS state is the unique KMS state of $\O_X$.}\medskip

Note that if each $\A_j$ has a unique trace and $A$ is aperiodic,
then ${\O_X}^0$ has a unique trace, which correponds necessarily to
the unique KMS state.\medskip

\heading
6. Inverse temperatures and topological entropy
\endheading

The aim of this section is to establish a relationship between inverse 
temperatures
of extremal
KMS states and the topological entropy of certain  subshifts naturally 
associated to
$(\A, \gamma).$

Let $\{y_i\}$ and $\{x_j\}$ be finite subsets of $\A^1\backslash\{0\}$ 
satisfying 
$\sum_i{y_i}^*y_i=I$ and $\sum_jx_j{x_j}^*=I,$ and let $T=T_{\{y_i\}}$ and 
$S=S_{\{x_j\}}$
be the completely positive maps of $\A^0$ defined in  
section 2.
We define
$$h(T_{\{y_i\}}):=
\lim_n \frac{1}{n}\log(\sum_{i_1,\dots, i_n}\|{y_{i_1}\dots
y_{i_n}}({y_{i_1}\dots
y_{i_n}})^*\|)
$$
$$h(S_{\{x_j\}}):=
\lim_n\frac{1}{n}\log(\sum_{j_1,\dots,j_n}\|({x_{j_1}\dots
x_{j_n}})^*x_{j_1}\dots
x_{j_n}\|).$$
In the  case  of the Matsumoto $C^*$--algebras $\O_\Lambda$,
 $h(S)$ is the {\it topological entropy} of the  shift homeomorphism
$\sigma\upharpoonright_\Lambda$, and 
 it was shown to
coincide with the maximal inverse temperture of certain KMS states in  
\cite{MWY}.

Note that since
 $(\sum\|y_{i_1}\dots y_{i_n}\|^2)^{1/n}\geq 1$ for all $n$ and
$\|y_i\|\leq1$, for all $i$, then
$$0\leq h(T_{\{y_i\}})\leq\log(\text{Card}\{y_i\}),$$
and similarly, 
$$0\leq h(S_{\{x_j\}})\leq\log(\text{Card}\{x_j\}).$$
Note however that if $\A$ is a crossed product
$C^*$--algebra by a single automorphism $\alpha$ 
then $h(\alpha)=h(\alpha^{-1})=0.$ More generally, since the sequences
 defining $h(T)$ and $h(S)$ converge to their greatest lower bounds,
we see that these are positive if and only if one has repectively
$$(\sum\|y_{i_1}\dots y_{i_n}\|^2)^{1/n}\geq 1+\varepsilon$$
$$(\sum\|x_{j_1}\dots x_{i_n}\|^2)^{1/n}\geq 1+\varepsilon$$
for all $n$ and some $\varepsilon>0$.

We now associate to a fixed finite set of nonzero elements 
$\{x_j\}\subset\A^1$
such that $\sum{x_j}{x_j}^*=I$, 
a {\it one--sided subshift} $\ell_{\{x_j\}}$, defined as 
follows. Set
 $\Sigma:=\{1,\dots, d\}$, where 
$d=\text{Card}\{x_j\}$. Then
$$\ell_{\{x_j\}}=\{\lambda\in 
\Sigma^{\Bbb N}: x_{\lambda_1}\dots x_{\lambda_r}\neq0, r\in{\Bbb N}\}.$$
Clearly  $\ell_{\{x_j\}}$ is a closed subset of $\Sigma^{\Bbb N}$
mapped onto itself by the left shift homomorphism:
$$\sigma(\lambda)_i=\lambda_{i+1}.$$

Notice that if $(\A, \gamma)$ does not result from a crossed product by 
an automorphism,
or a proper corner endomorphism, $d\geq2$.

Note also that, thank to the realation $\sum_i{x_j}{x_j}^*=I$,
any $n$--tuple $(\lambda_1,\dots,\lambda_r)$ $\in\Sigma^r$ 
such that $x_{\lambda_1}\dots x_{\lambda_r}\neq 0$ extends to an
element of
$\ell_{\{x_j\}}$. In particular, $\ell_{\{x_j\}}$ is nonempty.

Replacing the ${\Bbb T}$--action $\gamma$ by the action
$z\in{\Bbb T}\to {\gamma'}_z:=\gamma_{z^{-1}}$,
we see that we also have, for any finite
subset $\{y_i\}\subset \A^1\backslash\{0\}$ such that
$\sum_i{y_i}^*y_i=I$, a one--sided subshift:
$${\ell'}_{\{y_j\}}=\{\lambda\in{\Sigma'}^{\Bbb N}: y_{\lambda_r}\dots
y_{\lambda_1}\neq 0, r\in{\Bbb N}\},$$ 
where $\Sigma'$ is the set of the first $d':=\text{Card}\{y_i\}$ positive
integers.

We also introduce the following two--sided subshifts:
$$\Lambda_{\{x_j\}}=\{\lambda\in\Sigma^{\Bbb Z}: x_{\lambda_r}\dots
x_{\lambda{r+s}}\neq 0, r\in{\Bbb Z}, s\in{\Bbb N}\},$$
and
$${\Lambda'}_{\{y_i\}}=\{\lambda\in{\Sigma'}^{\Bbb
Z}: y_{\lambda_{r+s}}\dots
y_{\lambda{r}}\neq 0, r\in{\Bbb Z}, s\in{\Bbb N}\}.$$

\noindent{\it Remark}
Even though it would seem more convenient to work with two--sided
subshifts,
we should point out that these may be rather small, in the sense
that 
a finite word $(\lambda_1,\dots,\lambda_r)$
occurring, e.g., in $\ell_{\{x_j\}}$
does not necessarily extend to a word in $\Lambda_{\{x_j\}}$. 
The following simple example well describes the situation.
Consider the $C^*$--algebra $\A=M_2({\Bbb C})\otimes \C({\Bbb T})$,
and define the following $2\pi$--periodic automorphic action $\alpha$
of ${\Bbb R}$:
$$\alpha_t=\text{ad}v_t\otimes \beta_t$$
where $$v_t=\text{diag}(e^{it}, 1)$$ and
$$\beta_t(f)(e^{i\tau})=f(e^{i(\tau+t)}).$$
Let $e_{i,j}$, $i,j=1,2$, be a system of matrix units for $M_2({\Bbb C})$,
and define
$x_1=e_{1,2}\otimes I$, $x_{2}=e_{2,2}\otimes u$, $y_1=e_{1,2}\otimes I$,
$y_2=e_{1,1}\otimes u$,
where $u(z)=z$, $z\in{\Bbb T}.$
Then all the above elements are in $\A^1$, and satisfy
$x_1{x_1}^*+x_2{x_2}^*=I$, ${y_1}^*y_1+{y_2}^*y_2=I$,
so $(\A, \alpha)$ is a full $C^*$--dynamical system.
Note that ${x_1}^2=0$, $x_1x_2\neq0$, $x_2x_1=0$, ${x_2}^2\neq0$,
so
$$\Lambda_{\{x_j\}}=\{(\dots,2,2,2,\dots)\}$$
while
$$\ell_{\{x_j\}}=\{(1,2,2,\dots),(2,2,2,\dots)\}.$$
Thus $(1,2)$ is a $2$--word appearing in $\ell_{\{x_1,x_2\}}$ which can
not be extended to any two--sided sequence of $\Lambda_{\{x_1,x_2\}}$.

We  give a condition ensuring that $\ell_{\{x_j\}}$ and
${\ell'}_{\{y_i\}}$ are the positive parts of $\Lambda_{\{x_j\}}$ and 
${\Lambda'}_{\{y_i\}}$ respectively.
\medskip

\noindent{\bf 6.1. Proposition} {\sl  Let $\{x_j\}$,  
be 
a finite
subset of $\A^1\backslash\{0\}$ such that 
$$\sum_j{x_j}{x_j}^*=I$$ and 
let
$\ell_{\{x_j\}}$ and $\Lambda_{\{x_j\}}$  be the associated 
one--sided and two--sided subshifts respectively.
If $\sum_j {x_j}^*{x_j}$  is invertible then
$\Lambda_{\{x_j\}}\neq\emptyset$. Furthermore
$$\ell_{\{x_j\}}={\Lambda_{\{x_j\}}}_+.$$
An analogous statement holds for ${\ell'}_{\{y_i\}}$ and
${\Lambda'}_{\{y_i\}}$.}\medskip

\noindent{\it Proof} Since $\sum_j{x_j}{x_j}^*=I$, for any $(i_{1},\dots, 
i_{r})$ such that
$x_{i_{1}}\dots x_{i_{r}}\neq0$ there is an $i_{r+1}$ such that 
$x_{i_1}\dots x_{i_r}x_{i_{r+1}}\neq0$, and therefore there is a 
sequence $(i_n)_{n\geq 1}$ such that $x_{i_1}\dots x_{i_n}\neq0$ for 
all $n\in{\Bbb N}$.
To complete the proof relative to the set $\{x_j\}$ it is now
clear that it suffices to show, for any such $(i_n)$, 
the existence of $i_0\in\Sigma$ such that $x_{i_0}x_{i_1}\dots
x_{i_n}\neq0$ for all 
$n\geq0$.
If this were not the case, for all $k\in \Sigma$ there would exist
$n_k$ 
such that
$x_{{k}}x_{i_1}\dots x_{i_{n_k}}=0$. 
Letting $n=\max\{n_k, k\in
\Sigma\}$,
we must have
$x_{{k}}x_{i_1}\dots x_{i_{n}}=0$ 
 for all $k\in \Sigma$, and 
therefore,  $\sum_k{x_k}^*x_k$ being invertible,
$x_{i_1}\dots x_{i_{n}}=0$. This is 
now a contradiction.
The statement relative to the set $\{y_i\}$ can be proved 
similarly.\medskip

In particular, if  $A=(a_{ij})$ is a 
 $\{0,1\}$--matrix with no zero row or column,
then the generating partial isometries $\{S_i\}$ of the
Cuntz--Krieger 
algebra
$\O_A$
satisfy both   
$$\sum_{i} S_i{S_i}^*=I,$$
and 
$${S_i}^*S_i=\sum_j a_{ij}S_j{S_j}^*,$$
thus
$\sum_i {S_i}^*S_i$ is invertible. 
One has $\Lambda_{\{S_i\}}=\Lambda_A$ and $\ell_{\{S_i\}}={\Lambda_A}_+$.
More generally, if $\Lambda$ is a nonempty subshift of $\Sigma^{\Bbb Z}$
then the canonical set of generating partial isometries
$\{S_i\}$ of the Matsumoto $C^*$--algebra $\O_\Lambda$
still satisfy the above conditions
(see 4.1) and we have also in this case 
$\Lambda_{\{S_i\}}=\Lambda$ and $\ell_{\{S_i\}}=\Lambda_+$.
Another example is provided by the 
algebras generated by certain Cuntz--Krieger bimodules
$X=\oplus_{(i,j): a_{i,j}=1} X_{i,j}$ as described in section 5. 
More precisely, if $X_{i, j}$ is the Hilbert bimodule defined by 
a unital $^*$--isomorphism $\phi_{i,j}: \A_i\to\A_j$, then one can 
define $S_{i,j}$ to be the identity of $\A_j$ regarded as an element of
$X_{i,j}$. So $S_i=\sum_{j: a_{ij}=1} S_{i,j}$ are partial isometries of 
${\O_X}^1$
satisfying the Cuntz--Krieger relations with respect to $A=(a_{ij})$.
Therefore $\Lambda_{\{S_i\}}$ is again the two--sided Markov subshift
defined by the matrix $A=(a_{i,j})$.
The example arising from fractal geometry discussed
in section 4 is in the same spirit, 
in that the natural basis of the generating module are generators of a 
Cuntz
algebras, so the associated one or two--sided subshifts are full.

Note that all the examples above discussed have in common the fact
that there is   a multiplet $\{x_i\}\subset\A^1$ such that
$\sum_ix_i{x_i}^*=I$ consisting 
of elements with pairwise orthogonal ranges (and therefore they are 
necessarily
partial isometries).

We start by establishing  general estimates for the extremal inverse 
temperatures
using the topological entropies of the associated subshifts.  
Recall \cite{DGS} that for a one--sided (or two--sided) subshift $(\ell, 
\sigma\upharpoonright_{\ell})$
the topological entropy can be computed as
$$h_{\text{top}}(\sigma\upharpoonright_{\ell})=
\lim_n \frac{1}{n}\log(\theta_n(\ell)),$$
where $\theta_n(\ell)$ is the cardinality of the set $\ell^r$ of distinct
words 
of length $n$ occurring in $\ell$.
\medskip

\noindent{\bf 6.2. Proposition} {\sl Let $\{y_i\}$,
$\{x_j\}\subset\A^1\backslash\{0\}$ 
be 
finite
subsets such that $\sum_i{y_i}^*y_i$ $=I$ and $\sum_j x_j{x_j}^*=I$, and 
let
${\ell'}_{\{y_i\}}$ and $\ell_{\{x_j\}}$ be the corresponding
one--sided subshifts, defined 
as above.
Then
$$\beta_{\text{min}}\geq-h(T_{\{y_i\}})\geq
-h_{\text{top}}(\sigma\upharpoonright_{{\ell'}_{\{y_i\}}}),$$
$$\beta_{\text{max}}\leq h(S_{\{x_j\}})\leq
h_{\text{top}}(\sigma\upharpoonright_{\ell_{\{x_j\}}}).$$}\medskip

\noindent{\it Proof} By Prop. 2.3 and the triangle inequality 
$\beta_{\text{min}}\geq-h(T_{y_i})$
and $\beta_{\text{max}}\leq h(S_{\{x_j\}})$. The rest follows from
$\sum_{i_1,\dots, i_n}\|x_{i_1}\dots x_{i_n}\|^2\leq
\theta_n(\ell_{\{x_j\}}),$
and its analogue for $\{y_i\}$.\medskip

We shall see that in the general situation if it is possible to choose 
the multiplets $\{y_i\}$ and $\{x_j\}$ carefully, then the 
topological entropies of the corresponding subshifts
lead  to the extremal inverse temperatures of KMS states.
We first present an intermediate result, which gives
a sufficient condition for $h(T)$ and $h(S)$ to coincide with 
the topological entropy of the associated subshifts.\medskip

\noindent{\bf 6.3. Proposition} {\sl Set $$l_n:=
\min\{\|y_{i_1}\dots y_{i_n}\|^2: y_{i_1}\dots y_{i_n}\neq0\},$$ 
$$m_n:=\min\{\|x_{j_1}\dots x_{j_n}\|^2: x_{j_1}\dots x_{j_n}\neq0\}.$$
 If 
$$
\limsup_n l_n^{1/n}=1$$
then $$h(T_{\{y_i\}})=
h_{\text{top}}(\sigma\upharpoonright_{{\ell'}_{\{y_i\}}}).$$
Similarly, if
$$
\limsup_n m_n^{1/n}=1$$ then $$h(S_{\{x_j\}})=
h_{\text{top}}(\sigma\upharpoonright_{\ell_{\{x_j\}}}).$$}\medskip

\noindent{\it Proof} We shall prove only the first statement.
Let $\theta_n$ denote the number of words of length
$n$ occurring in ${\ell'}_{\{y_i\}}$, i.e. the number on $n$-tuples 
$(i_1,\dots,i_n)\in {\Sigma'}^n$ such that $y_{i_1}\dots y_{i_n}\neq0.$
Then given $\varepsilon>0$, for infinitely many indices $n$,
$$(1-\varepsilon){\theta_n}^{1/n}\leq{\theta_n}^{1/n}{l_n}^{1/n}
\leq (\sum \|y_{i_1}\dots y_{i_n}\|^2)^{1/n}\leq{\theta_n}^{1/n},$$
hence taking the logarithm of the limit over $n$,
$$h(T)=\lim_n \frac{1}{n}\log\theta_n=h_{\text{top}}
(\sigma\upharpoonright_{{\ell'}_{\{y_i\}}}).$$\medskip

The previous result applies  whenever  one is working
with a multiplet consisting of partial isometries with 
mutually orthogonal ranges.

We next show that, strengthening   the hypotheses of the previous 
result,
all the inequalities of Proposition 6.2 become equalities.
More precisely, if the positive evaluations of a tracial state
 $\tau$
on the iterated basic monomials, e.g. 
$(x_{i_1}\dots{x_{i_n}})^*x_{i_1}\dots{x_{i_n}}$,
 do not get too small when $n$ increases, then 
the maximal inverse temperature
$\beta_{\text{max}}$  can be approximated iterating the operator
$s'$ on $\tau$. The proof is inspired by an 
analogous 
 result in \cite{MWY} for the Matsumoto algebras associated to 
subshifts.\medskip

\noindent{\bf 6.4. Theorem} {\sl Let $\{y_i\}$, $\{x_j\}$ be finite 
subsets
of $\A^1$ such that $\sum{y_i}^*y_i=I$ and $\sum x_j{x_j}^*=I$.
If $\tau$ is a tracial state on $\A^0$, set
 $$\mu_n(\tau):=\min\{\tau(y_{i_1}\dots 
y_{i_n}(y_{i_1}\dots y_{i_n})^*): y_{i_1}\dots y_{i_n}\neq0\},$$ 
$$\nu_n(\tau):=\min\{\tau((x_{j_1}\dots 
x_{j_n})^*x_{j_1}\dots y_{j_n}): x_{j_1}\dots x_{j_n}\neq0\}.$$
If
$$\limsup_n {\mu_n}^{1/n}(\tau)=1$$ then 
$$\beta_{\text{min}}=-h(T_{\{y_i\}})=
-h_{\text{top}}(\sigma\upharpoonright_{{\ell'}_{\{y_i\}}})=$$
 $$=-\limsup_n \frac{1}{n}\log(\delta_n(\tau)).$$
In particular, if $\tau$ 
is the restriction of a  KMS state $\omega$, then $\omega$ has minimal
inverse temperature.

If instead
$$\limsup_n {\nu_n}^{1/n}(\tau)=1$$
then
$$\beta_{\text{max}}=h(S_{\{x_j\}})=
h_{\text{top}}(\sigma\upharpoonright_{\ell_{\{x_j\}}})=$$
$$=\limsup_n \frac{1}{n}\log(\epsilon_n(\tau)).$$
If $\tau$ is the restriction of a   KMS state $\omega$, 
then $\omega$ has maximal inverse temperature.}\medskip

\noindent{\it Proof} We shall prove only the first statement. For 
infinitely many
indices $n$,
$$(1-\varepsilon)^n\theta_n\leq\theta_n\mu_n(\tau)\leq \delta_n(\tau)=
\tau(\sum y_{i_1}\dots{y_{i_n}}
(y_{i_1}\dots{y_{i_n}})^*)\leq \theta_n,$$
hence taking the $n$-th root and then the logarithm of the limit over a
subsequence,
 by the arbitrariness of $\varepsilon$, we get 
$$\limsup_n\frac{1}{n}\log(\delta_n(\tau))=
h_{\text{top}}(\sigma\upharpoonright_{{\ell'}_{\{y_i\}}})=h(T_{\{y_i\}}).$$
The last equality follows from Prop. 6.3. Now
$\delta_n(\tau)^{1/n}=\|{t'}^n(\tau)\|^{1/n}\leq \|{t'}^n\|^{1/n}\to
r(t')$, so
$$h_{\text{top}}(\sigma\upharpoonright_{{\ell}'_{\{y_i\}}})\leq
-\beta_{\text{min}},$$
which, together with Proposition 6.2, proves the first statement. If in
particular
$\tau$ arises from a KMS state $\omega$ at inverse temperature 
$\beta$, then $\delta_n(\tau)=e^{-n\beta}$, so 
$\beta=\beta_{\text{min}}.$
\medskip

The previous result can be regarded as an analogue of the well known
fact from Perron--Frobenius theory that for an
irreducible nonnegative
 matrix $A$ the maximal eigenvalue  can be approximated
by $r(A)=\limsup_n\|A^n(\tau)\|^{1/n}$, where $\tau$ is any vector with 
positive 
entries \cite{G}.\medskip

\noindent{\bf 6.5. Corollary} {\sl If there is a finite subset $\{y_i\}$ 
(resp. $\{x_j\}$)
of $\A^1$ 
such that
\roster
\item
$$\sum_i{y_i}^*{y_i}=I\quad\text{(resp.}\  \sum_j x_j{x_j}^*=I\text{)},$$
\item
$$y_i{y_h}^*=0,\quad i\neq h\quad\text{(resp.}\  {x_j}^*x_k=0, j\neq 
k\text{)},$$
\item The algebra $\C$ generated by all finite products of the form
$$y_{i_1}\dots y_{i_n}(y_{i_1}\dots y_{i_n})^*$$
$$\text{(resp.}\ (x_{j_1}\dots x_{j_n})^*x_{j_1}\dots x_{j_n})$$ is 
finite--dimensional,
\endroster
then the conclusions of the previous theorem hold for any faithful 
tracial
state on $\A^0$.}\medskip

\noindent{\it Proof} Just note that   under our assumptions any of the 
nonzero
basic monomials is a projection, and therefore it majorizes a minimal
projection in $\C$. It follows that $\lim_n \mu_n(\tau)^{1/n}=1$
(resp. $\lim_n \nu_n(\tau)^{1/n}=1$) for any faithful trace 
$\tau$ on $\A^0$, so the previous theorem applies.\medskip

In particular, this result applies to all the examples discussed at the 
beginning
of this section.

\heading
7. Topological entropy of canonical ucp maps
\endheading

Let
$\{x_j\}$ be a finite set of  a $C^*$--algebra $\A$ of grade $1$
 such that $\sum_j x_j{x_j}^*=I$.
In the previous section we have associated to this set a one--sided 
subshift 
$(\ell_{\{x_j\}}, \sigma\upharpoonright_{\ell_{\{x_j\}}})$
of the Bernoulli shift $(\Sigma^{\Bbb N}, \sigma)$, where $\Sigma$
is the state space of the first $d$ positive integers, and
$d=\text{Card}\{j: x_j\neq0\}$, 
in a way that, under suitable circumstances, its classical topological 
entropy
equals an extremal inverse temperature of KMS states.
One can also associate to the subset 
$\{x_j\}$  a unital completely positive map defined by
$$\sigma_{\{x_j\}}: T\in\A\to \sum_j x_jT{x_j}^*.$$
In view of the results of the previous section,
we ask whether there is a relationship between the
the Voiculescu topological entropy of this map and the classical
topological entropy of the subshift $\ell_{\{x_j\}}$.
If $\A=\O_n$ is the Cuntz algebra with generators
$S_1,\dots, S_n$, Choda shows in
\cite{Ch} that
the topological entropy of the canonical endomorphism $\sigma_{\{S_i\}}$
is $\log(n)$, i.e. the topological entropy of the associated full shift.
In the more general case where
 $\A=\O_A$ is a Cuntz--Krieger algebra, 
and $\{S_j\}$ is the canonical set of generating partial isometries,
Boca and Goldstein \cite{BG} have recently computed the Voiculescu
entropy 
\cite{V}
of this map, and they have shown   that it equals the 
logarithm of the 
spectral radius of $A$, or, in other words, the classical topological
 entropy of the underlying
finite type subshift $\ell_A$ \cite{DGS}. However, special cases, although
extreme from a certain point of view, of full periodic $C^*$--dynamical
systems are the crossed products by an automorphism $\alpha$.  Brown
showed in 
\cite{B} that $\text{ht}_{\A^0\rtimes{\Bbb
Z}}(\text{Ad}(u))=\text{ht}_{\A^0}(\alpha)$, where
$u$ is a unitary implementing $\alpha$. It is obvious that the associated
subshift is in this case trivial, so its entropy is zero.

In Theorem 7.4 we give, for full periodic $C^*$--dynamical systems,
an  upper bound for $\text{ht}(\sigma_{\{x_j\}})$
which allows to recover the above discussed results as special cases.
 We
will then apply this result to find
new examples, among the Matsumoto algebras associated to non finite type
subshifts, where
$$\text{ht}(\sigma_{\{x_j\}})=h_{\text{top}}(\ell_{\{x_j\}})$$
still holds.

In the beginning of this section the automorphic action of the circle 
plays no role, therefore we
shall not assume that the $x_j$'s are of grade $1$. We define the 
associated one--sided subshift 
$\ell=\ell_{\{x_j\}}$ as in the previous section.

We now show that the ucp map $\sigma_{\{x_j\}}$ can be understood as a
noncommutative subshift.
 Let
  $\sigma\upharpoonright_{\ell}$ be the restriction of $\sigma$ 
to $\ell$ and  
$T_{\ell}$ the $^*$--monomorphism of $\C(\ell)$ obtained by 
transposing 
$\sigma\upharpoonright_{\ell}$, i.e.
$$T_{\ell} f(x)=f(\sigma(x)),\quad x\in\ell.$$\medskip
Also, we will consider a natural basis of neighborhoods for  
$\ell$.
For each $(i_1,\dots, i_r)\in{\ell}^{r}$, consider the cylinder 
set
$$[i_1\dots i_r]=\{(x_j)_j\in\ell: x_1=i_1,\dots, x_r=i_r\}.$$
For a fixed $r\in{\Bbb N}$, these constitute an open and closed cover
of $\ell$ with cardinality 
$\theta_r=\text{Card}{\ \ell}^{r}$.\medskip

\noindent{\bf 7.1. Proposition} {\sl Let $\{x_j\}$ be a finite subset of 
$\A$ such that
$\sum_j x_j{x_j}^*=I,$ and let 
 ($\ell, \sigma\upharpoonright_{\ell}$) be the associated 
one--sided
subshift.  Then
there 
is a unique unital completely positive map $\Phi: \C(\ell)\to\A$ 
taking the characteristic function of $[i_1\dots i_r]$ to
 $x_{i_1}\dots x_{i_r}(x_{i_1}\dots x_{i_r})^*$. One has
$\sigma_{\{x_j\}}\circ \Phi=\Phi\circ T_{\ell}.$ Moreover, if the 
sequence 
$(m_n)_n$ defined in Proposition 6.3 does not converge to $0$,
then $\Phi$ is faithful.}\medskip

\noindent{\it Proof} 
We first notice that, for each $r\in{\Bbb N}$, the map 
 $\Phi_r: \Sigma^{\Bbb N}\to \Sigma^r$ projecting
onto the first $r$ coordinates takes  $\ell$ onto the subset
${\ell}^{r}$ of $\subset \Sigma^r$
consisting of $\theta_r$ elements. Therefore there is a natural
$^*$--monomorphism $\phi_r: {\Bbb C}^{\theta_r}\to \C(\ell)$ 
taking  a $\theta_r$--tuple assuming value $1$ on $(i_1,\dots, 
i_r)$
and zero elsewhere to the characteristic function of $[i_1\dots i_r]$.
 Similarly,
there are, for $r\leq s$, natural $^*$--monomorphisms 
$\phi_{r,s}:{\Bbb C}^{\theta_r}\to {\Bbb C}^{\theta_s}$ such that 
$\phi_s\phi_{r,s}=\phi_r$.
Since the cylinder sets $\{[i_i\dots i_r], r\in{\Bbb N}\}$ form a basis
of closed and open sets for $\ell$, we see that the image of all 
the $\phi_r$ is dense.
It follows that $\C(\ell)$ is the inductive limit of the ${\Bbb 
C}^{\theta_r}$'s
under the maps $\phi_{r,s}$.
We define the ucp map $\Phi_r:{\Bbb C}^{\theta_r}\to\A$
which takes
 the characteristic function of
$[i_i,\dots, i_r]$ to the element $x_{i_1}\dots x_{i_r}(x_{i_1}\dots 
x_{i_r})^*$. Since 
$\Phi_s\phi_{r,s}=\Phi_r$
$s\geq r$, thanks to  $\sum_j x_j{x_j}^*=I$, we get
a ucp map $^0\Phi:\cup_r{\Bbb C}^{\theta_r}\to\A$. Since 
$\|\Phi_r(f)\|\leq\|f\|$,
$^0\Phi$ extends to a ucp map on $\C(\ell)$, which is the desired
$\Phi$. 
The relation $\Phi\circ T_{\ell}=\sigma_{\{x_j\}}\circ \Phi$ can be 
easily checked
on the total set of characteristic functions of cylinder sets.

We now construct a conditional expactation $E_r:\C(\ell)\to {\Bbb 
C}^{\theta_r}$.
Choose a faithful normalized Borel measure $\mu$ on $\ell$, and associate
to a 
function 
$f\in\C(\ell)$ the $\theta_r$--tuple with coordinates 
$$E_r(f)(i_1,\dots,i_r)= \frac{1}{\mu([i_1\dots 
i_r])}\int_{[i_1,\dots,i_r]}f(x)d\mu(x),$$
for each $(i_1,\dots, i_r)\in{\ell}^{r}$. One can easily check
that $E_r(I)=I$ and that $E_r(fa)=E_r(f)a$, $a\in{\Bbb C}^{\theta_r}$. 
Clearly $(E_r)_r$
converges pointwise in norm to the identity.
Assume now that $f\in\C(\ell)$ is a positive element such that 
$\Phi(f)=0$.
Then $\Phi E_r(f)$ converges to $0$. On the other hand 
$$\|\Phi E_r(f)\|=\|\sum_{i_1,\dots, i_r}(\frac{1}{\mu([i_1\dots i_r])}
\int_{[i_1\dots i_r]} f(x)d\mu(x))x_{i_1}\dots 
x_{i_r}(x_{i_1}\dots x_{i_r})^*\|\geq $$
$$(\frac{1}{\mu([i_1\dots i_r])}\int_{[i_1\dots i_r]} f(x)d\mu(x)) 
m_r,$$
therefore if $m_r$ does not converge to $0$, a subsequence of  $E_r(f)$ 
converges to $0$, 
so $f=0$.\medskip

Note that if $(m_n)_n$ does not converge to $0$ then $\limsup_n 
m_n^{1/n}=1$ thus 
we are in the position of applying Proposition 6.3.
A particularly important case is when the $x_j$' have pairwise 
orthogonal ranges.\medskip

\noindent{\bf 7.2. Corollary} {\sl If there is a finite  subset $\{x_j\}
\subset\A$ such that
$$\sum_j x_j{x_j}^*=I,$$
$${x_i}^*x_j=0,\quad i\neq j,$$
then the ucp map $\Phi:\C(\ell)\to\A$ constructed in the previous 
proposition
is in fact a $^*$--monomorphism. The restriction of $\sigma_{\{x_j\}}$ 
to $\C(\ell)$
corresponds to the one--sided subshift 
$T_{\ell}$.}\medskip

For any subshift $\Lambda$ the Matsumoto $C^*$--algebra $\O_\Lambda$
satisfies 
the requirements of the previous result \cite{M}.

Let now ($\A$, $\gamma$) be a full $C^*$--dynamical system over ${\Bbb 
T}$.
Our next aim is to compare  the topological entropy of the ucp map 
$\sigma_{\{x_j\}}$,
when $\{x_1,\dots,x_d\}$ is a finite subset of $\A^1\backslash\{0\}$ such
that
$\sum_jx_j{x_j}^*=I$, with entropic properties of the canonical
homogeneus $C^*$--algebra. 
We refer the reader respectively to \cite{V}  for  the notion
 of
topological entropy for nuclear $C^*$--algebras and to \cite{B} for its
generalization to
exact $C^*$--algebras, and
to
\cite{BG} for the generalization of the topological entropy 
$\text{ht}(P)$ of a ucp map
$P$ on a unital exact $C^*$--algebra.

We start defining an entropic quantity for the homogeneous subalgebra
$\A^0$ which, in the case where $\A=\A^0\rtimes_\alpha{\Bbb Z}$,
reduces to the topological entropy of $\alpha$. We shall assume that
$\A^0$ is an exact $C^*$--algebra. We start fixing a choice of nonzero
elements $\{x_i\}$ of $\A^1$ such that $\sum_i x_i{x_i}^*=I$. 
Let us define, for $\mu=(i_1,\dots,i_r)\in{\ell}^{r}$,
$x_\mu:=x_{i_1}\dots x_{i_r}$. We set $x_{\emptyset}=I$ and
$|\emptyset|=0$. We shall also consider the operators 
$q_{\alpha,\beta}={x_{\alpha}}^*x_{\beta}\in\A^0$ for 
$|\alpha|=|\beta|\geq0$. Note that $q_{\emptyset,\emptyset}=I$.

Let be given $\pi:\A^0\to\B(\H)$  a faithful
$^*$--representation,        $\omega\subset\A^0$  a finite subset
and $\delta>0$.
Set, for $n\in{\Bbb N}$,
$$\omega^{(n)}=\{x_{\mu}^*q_{\delta,\delta'} T q_{\epsilon,
\epsilon'}x_\nu,$$
$$T\in\omega, |\mu|=|\nu|\leq n-1, |\delta|=|\delta'|\leq n-1,
|\epsilon|=|\epsilon'|\leq n-1\}.$$
Note that $\omega^{(n)}$ depends on the  contractions 
$\phi_{x_i, x_j}: T\in\A^0\to{x_i}^*T x_j$ rather then on the elements
$\{x_i\}$. 
We define, for $a>0$,
$$\text{ht}_a(\pi,\{\phi_{x_i, x_j}\},
\omega,\delta)=
\limsup_n\frac{1}{n}\log\text{rcp}(\pi, 
\omega^{(n)},\frac{\delta}{{\theta_{n-1}}^a}),$$
$$\text{ht}_a(\pi,\{\phi_{x_i, x_j}\},\omega)=\sup_{\delta>0}
\text{ht}_a(\pi, \{\phi_{x_i, x_j}\}, \omega,\delta),$$
$$\text{ht}_a(\pi, \{\phi_{x_i, x_j}\})=
\sup_{\omega\subset\A^0\text{finite}}\text{ht}_a(\pi, \{\phi_{x_i, x_j}\},
\omega).$$

We will use the same notation as \cite{B}. Thus in particular
for a finite set $\Omega\subset\A^0$, $\text{rcp}(\pi, \Omega, \delta)$
is computed with respect to factorizations of completely positive
contractions, not necessarily unital, from $\A^0$ to
$\B(\H)$ via finite
dimensional
$C^*$--algebras. 

Brown proves in \cite{B} that $\text{rcp}(\pi,\Omega,\delta)$ is
independent of the choice of $\pi$. We will regard $\A$ faithfully
represented on a Hilbert space $\H$, and we take $\pi$ to be the
inclusion $\iota_{\A^0}$ of $\A^0$ in $\B(\H)$.
Moreover
we will avoid indicating $\pi$ in the above definitions.

We anticipate, for later use, the following
 immediate consequence of the definition. \medskip

\noindent{\bf 7.3. Lemma}
{\sl 
\roster
\item
If $\alpha$ is an automorphism of a unital $C^*$--algebra $\A^0$,
$\A=\A^0\rtimes_\alpha{\Bbb Z}$, and $u$ is any unitary of $\A^1$
implementing $\alpha$ on $\A^0$, 
$$\text{ht}_a(\phi_{u,u})=\text{ht}(\alpha^{-1}).$$

\item
If $\omega$ is such that for some finite
dimensional $C^*$--subalgebra $\D\subset\A^0$ which is the range of a
conditional expectation, $\omega^{(n)}\subset
\D$ 
except for finitely many $n$, then $\text{ht}_a(\{\phi_{x_i, x_j}\},
\omega)=0$.
\endroster}\medskip

\noindent{\it Proof} (1) for $n\in{\Bbb N}$,
 $\theta_n=1$ and
$\omega^{(n)}=\omega\cup\dots\alpha^{-n+1}(\omega)$.
(2) Let $E:\A^0\to\D$ be a conditional expectation, and set $\phi:=E$,
$\psi:=\iota_{\D}$, so that for $T\in\omega^{(n)}$ and infinitely many
indices $n$,
$\psi\circ\phi(T)=T$. This implies that
$\text{rcp}(\omega^{(n)},\frac{\delta}{{\theta_{n-1}}^a})\leq\text{rank}(\D)$,

Let now $\omega\subset\A^0$ be a finite subset and $n_0\in{\Bbb N}$.
A typical finite subset of $\A$ has the form
$$F(\omega, n_0)=\{x_\gamma T, T\in\omega, |\gamma|\leq n_0\}.$$
Our aim is to show the following result.\medskip

\noindent{\bf 7.4.  Theorem} {\sl Let ($\A$, $\gamma$, ${\Bbb T}$) be 
a full $C^*$--dynamical system, with $\A^0$ exact.
 Let $\sigma_{\{x_i\}}$ and
$\ell_{\{x_i\}}$ be the ucp map and the one--sided subshift
associated to a set $\{x_i\}\subset\A^1$ satisfying $\sum_i
{x_i}{x_i}^*=I$.
Then for any finite subset $\omega\subset\A^0$ and $n_0\in{\Bbb N}_0$,
$$\text{ht}(\sigma_{\{x_i\}}, F(\omega, n_0))\leq
h_{\text{top}}(\sigma\upharpoonright_{\ell_{\{x_i\}}})+
\text{ht}_2(\{\phi_{x_i, x_j}\},\omega).$$}\medskip

We will prove this theorem combining appropriate analogues of arguments
of Brown \cite{B} for crossed product $C^*$--algebras and Boca and
Goldstein \cite{BG} for Cuntz--Krieger algebras.

Motivated by \cite{B}, we define certain cp maps.
Let $F\subset {\Bbb N}_0$ be a finite subset. Set
$$S_F: T\in\A\to({x_\alpha}^*m_{|\alpha|-|\beta|}(T)x_\beta)_{\alpha,
\beta\in I_F}\in M_{\theta_F}(\A^0).$$
Here $I_F:=\cup_{r\in F}\ell^r$,
$\theta_F=\sum_{r\in F}\theta_r$, and, for $k\in{\Bbb Z}$,
$m_k:\A\to \A^k$ is the natural projection.
Note that $S_F$ is contractive and cp.
For a contractive cp map $\phi: \A^0\to\B$ set
$$\phi_F:=\iota\otimes\phi\circ S_F:\A\to M_{\theta_F}(\B)$$
which is contractive and cp.

Let  $f\in\ell^2({\Bbb N}_0)$ have  support in $F$,
and define
the cp map $$\tilde{S}_{F, f}: T=(T_{\alpha,\beta})\in
M_{\theta_F}(\B(\H))\to\sum_{\alpha,\beta\in I_F}
f(|\alpha|)\overline{f(|\beta|)} x_\alpha T_{\alpha,\beta}x_{\beta}^*\in
\B(\H).$$ 
Note that $\tilde{S}_{F, f}(I)={\|f\|_2}^2 I$.
Again, for a contractive cp map 
$\psi:\B\to\B(\H)$, define
$$\psi_{F, f}:=\tilde{S}_{F,
f}\circ\iota\otimes\psi: M_{\theta_F}(\B)\to\B(\H).$$
So $\psi_{F, f}$ is cp contractive if  $\|f\|_2\leq1$.

Finally, for a contractive cp map $\Lambda:\A^0\to\B(\H)$
set 
$$\Phi_{\Lambda, F, f}:=\tilde{S}_{F, f}\circ\iota\otimes\Lambda\circ
S_F: \A\to\B(\H).$$
Note that in particular, if $\phi$ and $\psi$ are as above,
$$\Phi_{\psi\circ\phi, F, f}=\psi_{F, f}\circ \phi_{F}$$
which factors through the algebra $M_{\theta_F}(\B)$.
One can easily show that for an element of fixed degree
$X\in\A^k$,
$$\Phi_{\iota_{\A^0}, F, f}(X)=\sum_{p\in
F\cap(F+k)}f(p)\overline{f(p-k)} X.$$  
The following lemma is our analogue of Lemma 3.4 in \cite{B}.\medskip

\noindent{\bf 7.5.  Lemma} {\sl Let $\omega$ be a finite subset
of the unit ball of $\A^0$, $n_0\in{\Bbb N}$ and $\delta>0$. Consider the
set
$F'(\omega, n_0):=\{Tx_\gamma, T\in\omega, |\gamma|\leq n_0\}.$
Then there is a finite set $F\subset{\Bbb N}_0$
 which depends only on $n_0$ and $\delta$ and not on $\omega$
such that
$$\text{rcp}(F'(\omega, n_0),
\delta)\leq$$
$$\theta_F\text{rcp}(\cup_{|\gamma|\leq
n_0}\cup_{|\alpha|,|\beta|\in
F, |\alpha|=|\beta|+|\gamma|}{x_\alpha}^*\omega x_\gamma x_\beta,
\frac{\delta}{2\max_{p\in F}\theta_p}).$$}\medskip

\noindent{\bf Proof} We proceed as in the proof of Lemma 3.4 in
\cite{B}. Let $f\in\ell^2({\Bbb Z})$ be a function with finite support
$E$,
$\|f\|_2\leq1$ such that $|f*\tilde f(n)-1|<\delta/2$, $n=0,1,\dots,n_0$. 
Here $\tilde f(p)=\overline{f(-p)}$.
Replacing $f$ with a suitable translate if necessary, we may assume
$E\subset\{n\in{\Bbb Z}: n\geq n_0\}$. 
Set $F=E\cup(E-1)\cup\dots\cup (E-n_0)$, which is a subset of ${\Bbb
N}_0$.
Note that $F$ depends on $\delta$ and $n_0$ but not on $\omega$. Consider
contractive cp maps
$\phi:\A^0\to\B$ and $\psi: \B\to\B(\H)$ with $\B$ finite dimensional
such that 
$$\|\psi\circ\phi(a)-a\|<\frac{\delta}{2\max_{p\in
F}\theta_p},\quad a\in\cup_{|\gamma|\leq n_0}\cup_{|\alpha|,
|\beta|\in F, |\alpha|=|\beta|+|\gamma|}{x_\alpha}^*\omega x_\gamma
x_\beta.$$ 
Let us choose $\B$ with minimal rank.
Then the cp contractive map $\Phi_{\psi\circ\phi, F, f}$
factors
through the finite dimensional algebra $M_{\theta_F}(\B)$, which has
rank $\theta_F\text{rank}(\B)$. We are thus left to show that
for $a\in F'(\omega, n_0)$,
$$\|\Phi_{\psi\circ\phi, F, f}(a)-a\|<\delta.$$
We write $a=Tx_\gamma$ with $T\in\omega$, $|\gamma|\leq n_0$.
Then the l.h.s. is bounded by
$$\|\Phi_{\psi\circ\phi, F, f}(Tx_\gamma)-\Phi_{\iota_{\A^0}, F,
f}(Tx_\gamma)\|+\|\Phi_{\iota_{\A^0}, F, f}(Tx_\gamma)-Tx_\gamma\|.$$
Now the computation of $\Phi_{\iota_{\A^0}, F, f}$ given before
on elements with fixed degree shows that the second summand
is bounded by $$|\sum_{p\in F\cap (F+|\gamma|)}
f(p)\overline{f(p-|\gamma|)}-1|\ \|Tx_\gamma\|=|f*\tilde
f(|\gamma|)-1|\|Tx_\gamma\|<\delta/2.$$
We now evaluate the first summand.
$$\|\Phi_{\psi\circ\phi, F, f}(Tx_\gamma)-\Phi_{\iota_{\A^0}, F,
f}(Tx_\gamma)\|\leq
\|\iota_{M_{\theta_F}}\otimes(\psi\circ\phi-\iota_{\A^0})\circ
S_F(Tx_\gamma)\|=$$
$$\|\sum_{p\in F\cap(F+|\gamma|)}
(\sum_{|\alpha|=p, |\beta|=p-|\gamma|} e_{\alpha,
\beta}\otimes(\psi\circ\phi-\iota_{\A^0})({x_\alpha}^*Tx_\gamma
x_\beta))\|
=$$
$$\max_{p\in F\cap(F+|\gamma|)}
\|\sum_{|\alpha|=p, |\beta|=p-|\gamma|} e_{\alpha,
\beta}\otimes(\psi\circ\phi-\iota_{\A^0})({x_\alpha}^*Tx_\gamma
x_\beta)\|\leq\delta/2.$$
The last inequality follows from our choice of $\psi$ and $\phi$ and from 
the fact that if a matrix $A\in
M_{h,k}(\A^0)$ has entries of norm bounded by $c$ then
$\|A\|\leq(hk)^{1/2}c$. \medskip

 Consider
 the contractive cp map
$$\rho_r:\A\to M_{\theta_r}(\A)$$ taking $T\in\A$ to the matrix 
$({x_\mu}^*Tx_{\nu})_{\mu,\nu\in{\ell}^{r}}.$
(One can easily check that $\rho_r$ is a unital $^*$--monomorphism
with image  the corner algebra $P_rM_{\theta_r}(\A)P_r$, where
$P_r=({x_\mu}^*x_{\nu})$ 
is an orthogonal
projection.)

 For $n, n_0\in{\Bbb N}$, $m\geq n+n_0-1$, $l=0,\dots,n-1$,
$|\alpha|\leq 
n_0$, $T\in\A^0$, we compute
$$\rho_m(\sigma^l(x_{\alpha}T))=\rho_m(\sum_{|\eta|=l} 
x_{\eta\alpha}T{x_{\eta}}^*)=
(\sum_{|\eta|=l}x_{\mu}
^*x_{\eta\alpha}T{x_{\eta}}^*x_{\nu})_{\mu,\nu\in\ell^m}.$$
Writing $\mu=\gamma\mu'$, $|\gamma|=|\eta|+|\alpha|=l+|\alpha|$, and 
$\nu=\delta\nu'$, 
$|\delta|=|\eta|=l$
we have:
$$\rho_m(\sigma^l(x_{\alpha}T))=
(\sum_{|\eta|=m}{x_{\mu'}}^*q_{\gamma,\eta\alpha}Tq
_{\eta,\delta}x_{\nu'})_{\mu,\nu\in\ell^m}.$$
Setting, again, $\nu'=\epsilon\nu''$, with 
$|\epsilon|=|\mu'|=m-l-|\alpha|$, we obtain that
$$\rho_m(\sigma^l(x_{\alpha}T))=
(\sum_{|\eta|=l}
({x_{\mu'}}^*q_{\gamma,\eta\alpha}Tq_{\eta,\delta}x_{\epsilon})x_{\nu''
})_{\mu,\nu\in\ell^m}.$$
If now $T$ ranges over a finite subset $\omega\subset\A^0$,  we see 
that
the image of  
$$F(\omega, n_0)=\{x_{\alpha}T, T\in\omega, |\alpha|\leq 
n_0\}\subset\A$$
under the maps $\rho_{n+n_0-1}\circ\sigma^l$, $l=0,\dots, n-1$, is 
constituted by matrices of size $\theta_{n+n_0-1}$
with entries sums of at most $\theta_{n-1}$ elements in 
$F'(\omega^{(n+n_0)}, n_0)$.
\medskip

\noindent{\it Proof of Theorem 7.4.} 
We apply the previous Lemma to the sets $F'(\omega^{(n+n_0)}, n_0)$
for fixed $n_0$ and $\omega$ and all $n\in{\Bbb
N}$. 
Note that the corresponding set $F$ can be chosen independent of $n$.
 We can thus find for each $n\in{\Bbb N}$ 
a contractive cp map $\Lambda_n:\A\to\B(\H)$ factoring through a
finite
dimensional
algebra $\B$ of rank 
$$\theta_F\text{rcp}(\cup_{|\gamma|\leq n_0}\cup_{|\alpha|, |\beta|\in F,
|\alpha|=|\beta|+|\gamma|}{x_\alpha}^*\omega^{(n+n_0)}x_{\gamma}x_{\beta},
,\frac{\delta}{2{\theta_{n+n_0-1}}^2\max_{p\in F}\theta_p})\leq$$
$$\theta_F\text{rcp}(\omega^{(n+n_0+\max
F)},\frac{\delta}{2{\theta_{n+n_0-1+\max F}}^2\max_F\theta_p})$$
such that
$$\|\Lambda_n(Tx_\gamma)-\pi(Tx_\gamma)\|<\frac{\delta}{{\theta_{n+n_0-1}}^2},
\quad
T\in\omega^{(n+n_0)},\quad |\gamma|\leq n_0.$$

Consider the ucp map $\Psi_m: M_{\theta_m}(\B(\H))\to\B(\H)$ taking the 
matrix $(t_{\mu, \nu})$
to the operator $$\sum_{|\mu|=|\nu|=m}\pi(x_\mu)t_{\mu,\nu}\pi(x_\nu)^*.$$
Then the map
$\Psi_{n+n_0-1}\circ\iota_{M_{\theta_{n+n_0-1}}}\otimes
\Lambda_n\circ\rho_{n+n_0-1}:\A\to\B(\H)$ factors through an algebra
of rank bounded by
$$\theta_{n+n_0-1}\theta_F\text{rcp}(\omega^{(n+n_0+\max
F)}, \frac{\delta}{2\theta_{n+n_0+\max F-1}^2\max_F\theta_p}).$$
Thus if we show that
$$\|\Psi_{n+n_0-1}\circ\iota\otimes\Lambda_n\circ\rho_{n+n_0-1}
(\sigma^l(x_\gamma T))-\sigma^l(x_\gamma T)\|<\delta$$
for $T\in\omega$, $|\gamma|\leq n_0$ and $l=0,\dots, n-1$, we will
deduce
that 
$$\text{rcp}(F(\omega, n_0)\cup\dots\cup\sigma^{n-1}(F(\omega,
n_0)),\delta)\leq$$
$$\theta_{n+n_0-1}\theta_F\text{rcp}(\omega^{(n+n_0+\max
F)}, \frac{\delta}{2{\theta_{n+n_0+\max F-1}^2\max_{F}\theta_p}})$$ and 
the
conclusion will follow.
Now as $\Psi_m\circ\rho_m=\iota_{\A}$, it
suffices
to show that
$$\|\iota_{M_{\theta_{n+n_0-1}}}\otimes(\Lambda_n-\iota_{\A})
\circ\rho_{n+n_0-1}(\sigma^l(x_\gamma
T))\|<\delta.$$ This follows from our choice of $\Lambda_n$
and from the fact that entries of
$\rho_{m}(\sigma^l(x_\gamma T))$, $m= n+n_0-1$ are sums of at most 
$\theta_{n-1}$ elements of $F'(\omega^{(n+n_0)}, n_0)$.\medskip

\noindent{\bf 7.6. Corollary} {\sl Consider the same situation as in
Theorem 7.4. Let $(\omega_\alpha)_{\alpha\in A}$ be a net of finite
subsets  of $\A^0$ 
with total union.
 Then
$$\text{ht}({\sigma_{\{x_j\}}})\leq 
h_{\text{top}}(\sigma\upharpoonright_{\ell_{\{x_j\}}})+
\lim_{\alpha}\text{ht}_2(\{\phi_{x_i, 
x_j}\}, \omega_\alpha).$$}\medskip

\noindent{\it Proof} This is a straightforward consequence of the fact 
that
$\cup_{\alpha, n_0} F(\omega_\alpha, n_0)\cup
F(\omega_\alpha, n_0)^*$ is
total
in 
$\A$ and of the Kolmogorov--Sinai property of the entropy of
a ucp map,
\cite{V}, \cite{B}, \cite{BG}.
\medskip

\noindent{\it Remark} If in particular $\A=\A^0\rtimes_\alpha{\Bbb Z}$ and
$u\in\A^1$ is a unitary implementing $\alpha$ on $\A^0$ then $\ell_{u}$ is
a single point space, so its entropy is zero. By
Lemma 7.3 and the previous Corollary, we recover
Brown's result that
$\text{ht}_{\A}(\text{Ad}(u))\leq\text{ht}_{\A^0}(\alpha)$ (and therefore
one deduces an equality by monotonicity of topological  entropy
\cite{B}.)\medskip


The case where we can choose the $x_j$'s with pairwise orthogonal ranges
is of course of special interest, the Cuntz algebras, Cuntz--Krieger
algebras and Matsumoto algebras belonging to this class.
The next result shows that the estimate of the entropy can be made more
precise in this case.\medskip

\noindent{\bf 7.7. Theorem} {\sl Let $(\A, \gamma, {\Bbb T})$ be a full
$C^*$--dynamical system, with $\A^0$ exact. Let $\{x_j\}\subset\A^1$ be a
finite subset such that
$$\sum_j x_j{x_j}^*=I,$$
$${x_i}^*x_j=0,\quad i\neq j,$$
$$\sum_j{x_j}^*x_j\quad\text{is invertible}.$$
 Then
$$h_{\text{top}}(\sigma\upharpoonright_{\ell_{\{x_j\}}})\leq
\text{ht}({\sigma_{\{x_j\}}})\leq h_{\text{top}}
(\sigma\upharpoonright_{\ell_{\{x_j\}}})+\lim_\alpha\text{ht}_1(\{\phi_{x_i, 
x_j}\},\omega_\alpha),$$
where $(\omega_\alpha)_{\alpha\in A}$ is any net of finite subsets
of $\A^0$ with
total union in $\A^0$. 
If in particular 
for some  net
$(\omega_\alpha)_{\alpha\in A}$ 
$\text{ht}_2(\{\phi_{x_i,
x_j}\},\omega_\alpha)=0$, $\alpha\in A$, then
$$\text{ht}(\sigma_{x_j})=h_{\text{top}}
(\sigma\upharpoonright_{\ell_{\{x_j\}}}).$$}\medskip

\noindent{\it Proof} The proof of the second inequality $\leq$ goes
exactly as
that of Theorem 7.4 with the only exception that 
entries of $\rho_m(\sigma^l(x_\gamma T))$ are now already elements of 
$F'(\omega^{(n+n_0)}, n_0)$. We show that
$$\text{ht}({\sigma_{\{x_j\}}})\geq 
h_{\text{top}}(\sigma\upharpoonright_{\ell_{\{x_j\}}}).$$
By monotonicity of  topological entropy \cite{B}, \cite{V} and
Corollary 7.2,
$$\text{ht}({\sigma_{\{x_j\}}})\geq \text{ht}(T_{\ell}),$$
where, as before, $T_{\ell}$ denotes the $^*$--monomorphism of 
$\C(\ell)$ implemented by
the one--sided shift. We are thus left to show that
$\text{ht}(T_{\ell})\geq 
h_{\text{top}}(\sigma\upharpoonright_{\ell_{\{x_j\}}}).$
The proof is similar to that of [BG], which in turn goes back to 
[V, Proposition 4.6]. Let $\mu$ be a $\sigma$--invariant probability 
Borel measure on the two--sided subshift
$\Lambda=\Lambda_{\{x_j\}}$ defined before Prop. 6.1,
 and let us restrict
it to a
$\sigma$--invariant probability 
measure
on $\ell=\Lambda_+$. For any ucp map $\gamma: M\to\C(\ell)$, with $M$ 
finite
dimensional, let $h_{\mu, T_{\ell}}(\gamma)$
be defined as in [CNT], by means of the function 
$H_\mu(\gamma, T_{\ell}\gamma,\dots,{T_{\ell}}^{n-1}\gamma)$.
Reasoning as in [V, Prop. 4,6] we see that
$h_{\mu, T_{\ell}}(\gamma)\leq \text{ht}(T_{\ell}).$
Choosing $M={\Bbb C}^{\theta_n}$ and $\gamma: {\Bbb 
C}^{\theta_n}\to\C(\ell)$
the natural inclusion, then one finds, thanks to [CNT, Remark III.5.2], 
that
the classical measurable entropy $H_\mu(\sigma\upharpoonright_\Lambda)$ 
is $\leq$
$\text{ht}(T_{\ell})$. Taking the supremum over all invariant 
measures we obtain the claim,
by the classical variational principle for topological entropy [DGS, 
Theorem 18.8].\medskip

\noindent{\it Remark} It is natural to ask whether the upper bound for
$\text{ht}(\sigma_{\{x_j\}})$ described in the previous result can be
further improved to
$h_{\text{top}}(\sigma\upharpoonright_{\ell_{\{x_j\}}})+
\text{ht}_0(\{\phi_{x_i, x_j}\}).$\medskip

We next show that the estimates 
above obtained are good enough to compute 
 $\text{ht}(\sigma_{\{x_j\}})$ in the case of Cuntz--Krieger algebras.
This was first done  by Boca and Goldstein \cite{BG}.\medskip

\noindent{\bf 7.8. Corollary} {\sl \cite{BG} Let $\A=\O_A$ be a
Cuntz--Krieger
algebra
defined by a $\{0,1\}$--matrix $A$, and let $\{S_i\}_1^d$ be the canonical
set of generating partial isometries. 
Then
$\text{ht}(\sigma_{\{S_i\}})=h_{\text{top}}
(\sigma\upharpoonright_{\ell_{\{S_i\}}})=\log(r(A)).$}\medskip

\noindent{\it Proof} The second equality  is the well known
computation of entropy of finite type subshifts of order two.
See, for example, Proposition 17.12 in \cite{DGS}.
The first
equality will follow from the previous
Theorem provided we show that
 there is an increasing sequence $\omega_p$, $p\in{\Bbb N}$
of finite subsets of ${\O_A}^0$ with total union such that
$$\text{ht}_2(\{\phi_{S_i, S_j}\},\omega_p)=0,\quad p\in{\Bbb N}.$$
Set $P_i:=S_i{S_i}^*$ and define
$\omega_p:=\{S_{\alpha}P_i{S_\beta}^*, i=1,\dots,d, |\alpha|=|\beta|\leq
p\}$. It is clear that $\cup_p\omega_p$ is total in the homogeneous
subalgebra. One easily checks that, for fixed $p$ and all $n$,
${\omega_p}^{(n)}$
is contained in  the linear span of $\omega_p$, which
is finite
dimensional $C^*$--algebra, so by Lemma
7.3 $\text{ht}_2(\phi_{\{S_i, S_j\}}, \omega_p)=0$.\medskip

We next look at the class of Matsumoto algebras $\O_\Lambda$
associated to a general
subshift $\Lambda\subset\{1,\dots, d\}^{\Bbb Z}$ \cite{M}. See also 
subsection 4.1.
Let $\{S_i\}_1^d$ be the canonical set of generating partial isometries.
We recall from \cite{M} a few properties of $\O_\Lambda$. First, the
relations
$${S_i}^*S_j=0,\quad i\neq j\eqno(7.1)$$ 
$$\sum_i{S_i}{S_i}^*=I,\eqno(7.2)$$
which easily imply that
for any pair of words with the same  
length,
$q_{\alpha,\beta}:={S_\alpha}^*S_\beta$=0 unless
$\alpha=\beta.$
We will write $q_{\alpha}$ for $q_{\alpha,\alpha}$,
$\alpha\in\cup_r\Lambda^r$. Note that these are projections.
Furthermore  the following commutation
relations hold in $\O_\Lambda$: for $\mu,\nu\in\cup_r\Lambda^r$:
$$q_\mu S_\nu=S_\nu q_{\mu\nu},\eqno(7.3)$$
$$q_\mu q_\nu=q_\nu q_\mu.\eqno(7.4)$$
By (7.4) the algebra $\Q_l$ generated by the projections
$\{q_{\alpha},\alpha\in\cup_{k=0}^l\Lambda^k\}$ is commutative and
therefore finite dimensional. 

These properties imply that the   finite sets
$$\omega_{k,l}:=\{S_\alpha E {S_{\beta}}^*,
|\alpha|=|\beta|\leq k, \ E\text{ minimal projection in }\Q^l \}$$
 has total union in ${\O_\Lambda}^0$.
It follows that ${\O_\Lambda}^0$ is AF \cite{M}, so
$\O_\Lambda$ is a nuclear 
$C^*$--algebra.

Using properties $(7.1)$--$(7.4)$ one can show with tedious computations
that for all $n\in{\Bbb N}$,
$${\omega_{k,l}}^{(n)}\subset\{S_{\alpha} q {S_\beta}^*,
|\alpha|=|\beta|\leq k,\  q \text{  projection in }
\Q_{2n+\max(k,l)}\}.$$ 
This computation is aimed to give an estimate for $\text{ht}_a(\phi_{S_i,
S_j}, \omega_{k,l}).$\medskip

\noindent{\bf 7.9. Lemma} {\sl If $\O_\Lambda$ is the Matsumoto
$C^*$--algebra
associted to a subshift $\Lambda$ we have, for $a>0$, and $k,l\in{\Bbb
N}_0$,
$$\text{ht}_a(\phi_{S_i, S_j}, \omega_{k,l})\leq
2\limsup_n\frac{1}{n}\log(\text{dim}(\Q_n)).$$
}\medskip

\noindent{\it Proof} Let $\phi:{\O_\Lambda}^0\to\B$ and
$\psi:\B\to{\O_\Lambda}^0$ be unital completely positive maps such that
$\|\psi\phi(a)-a\|<\delta/\theta_k$, when $a$ ranges the projections of
$\Q_{2n+2\max(k,l)}$, and assume that $\B$ has minimal rank.
Consider as before  the maps $\rho_k:{\O_\Lambda}^0\to
M_{\theta_k}({\O_\Lambda}^0)$,
$\Psi_k: M_{\theta_k}({\O_\Lambda}^0)\to\O_{\Lambda}^0$ which satisfy
$\Psi_k\circ\rho_k=\iota_{{\O_\Lambda}^0}$. Define
$\phi':=\iota_{M_{\theta_k}}\otimes\phi\circ\rho_k$ and
$\psi':=\Psi_k\circ\iota_{M_{\theta_k}}\otimes\psi$.
Then if $a$ is a projection of $\Q_{2n+\max(k,l)}$ and
$|\alpha|=|\beta|\leq k$,
$$\|\psi'\circ\phi'(S_{\alpha}a{S_\beta}^*)-S_\alpha a{S_\beta}^*\|\leq
\|\sum_{|\gamma|=|\gamma'|=k} e_{\gamma,
\gamma'}\otimes(\psi\circ\phi-\iota)({S_\gamma}^*S_\alpha
a{S_\beta}^*S_{\gamma'})\|<\delta$$
because ${S_\gamma}^*S_{\alpha}a{\S_\beta}^*S_{\gamma'}$ is a 
projection in $\Q_{2n+2\max(k,l)}$.
Any element of ${\omega_{k,l}}^{(n)}$ being
of the form $S_\alpha a{S_\beta}$, we deduce that for all $n\in{\Bbb N}$,
and $\delta>0$,
$$\text{rcp}({\omega_{k,l}}^{(n)},
\delta)\leq\theta_k\text{rcp}(\text{Proj}(\Q_{2n+2\max(k,l)}),
\delta/\theta_k)\leq\text{dim}(\Q_{2n+2\max(l,k)}).$$
The last inequality follows from the existence
of a conditional expectation onto
$\Q_{2n+2\max(k,l)}$. The rest follows choosing $\delta$ of the form
$\frac{\delta}{{\theta_{n-1}}^a}$, 
taking the logarthm, dividing
by $n$ and passing to the $\limsup$.
\medskip

We combine the previous result with Theorem 7.6.\medskip

\noindent{\bf 7.10. Theorem} {\sl If $\sigma_{\{S_i\}}$ is the ucp map
associated
to the canonical set of generators
of a Matsumoto $C^*$--algebra $\O_\Lambda$,
$$h_{\text{top}}(\Lambda)\leq \text{ht}(\sigma_{\{S_i\}})\leq
h_{\text{top}}(\Lambda)+2\limsup_n
\frac{1}{n}\log(\text{dim}(\Q_n)).$$}\medskip

We conclude this section discussing two examples of subshifts
 for which
$\text{ht}(\sigma_{\{S_i\}})=h_{\text{top}}(\Lambda)$.
The first example beyond Markov shifts is that of sofic subshifts, see
\cite{DGS} and therein
quoted references.
 By \cite{M} a subshift is sofic if and only if $\cup_n\Q_n$ is finite
dimensional. Then  Theorem 7.10 yields the desired
equality. Another example is that of $\beta$--shifts associated to 
$\beta$--expansion of real numbers \cite{Re}, \cite{Par}, \cite{Bl}.
In this case it is proved in \cite{KMW} that if the $\beta$--shift is not
sofic, $\text{dim}(\Q_n)=n+1$,
and this leads again to the same conclusion. 

 \heading{8. CNT dynamical entropy and variational principle}
\endheading

In fact, if $\A=\O_n$ Choda shows in \cite{Ch} not only that
$h_{\text{top}}(\sigma)=\log(n)$ but also that, if $\phi$ is the unique KMS state
of $\O_n$, then 
$h_\phi(\sigma)=h_{\text{top}}(\sigma)=\log(n)$, where the l.h.s.
denotes the Connes--Narnhofer--Thirring dynamical entropy of $\sigma$ \cite{CNT}.
This result has its own importance, as it exhibits
 a fundamental example where a noncommutative variational principle
for the entropy holds true. (We refer the reader to \cite{DGS} for a
formulation of the variational
principle for the entropy in ergodic theory for compact spaces.)

It is an open problem whether a noncommutative variational principle 
for the entropy of $C^*$--algebras holds. In this section we give a
 class of examples for which this is true, thus generalizing Choda's result. 
Examples will be  the  canonical ucp map 
of the Cuntz--Krieger algebras, or certain  
 Matsumoto algebras associated to non finite type subshifts.

In this section we show that
the CNT dynamical entropy of the ucp map $\sigma_{\{x_j\}}$,
defined in the previous section is $\geq$  the m.t. entropy
of the associated subshift $\Lambda_{\{x_j\}}$, as defined, e.g., in \cite{DGS}, see Theorems
8.5, 8.6.
This inequality looks similar to that of Theorem 7.7
relative to the topological entropy, but it goes in the  reverse order.
We start establishing  the setting of the CNT entropy.

Let $\A$ be a unital $C^*$--algebra, and let $\gamma_i:\A_i\to\A$,
$i=1,\dots,n$ be ucp maps from finite--dimensional $C^*$--algebras,
 and let $\phi$ be a state on $\A$. Let us recall from \cite{CNT} that an
Abelian model
for $(\A, \phi, \gamma_1,\dots,\gamma_n)$ is given by
an Abelian finite--dimensional $C^*$--algebra $\B$, a state $\mu$ on $\B$
and subalgebras $\B_1,\dots,\B_n$ of $\B$ for which there is a ucp map
$E:\A\to\B$ with $\phi=\mu\circ E$. Consider first the entropy 
of the Abelian model $(\B, \mu, \B_1,\dots,\B_n)$ as defined in 
[CNT,  III.3] and then the quantity $H_\phi(\gamma_1,\dots,\gamma_n)$,
defined as the supremum of the entropies of all the Abelian models
(see [CNT, Definition III.4]). The following result
is an  obvious consequence of the definition.\medskip

\noindent{\bf 8.1. Lemma} {\sl If, for $i=1,\dots, n$,
  $\gamma_i:\A_i\to\A$ is a $^*$--monomorphism,
$\vee_{i=1}^n\gamma_i(\A_i)$ is finite--dimensional and commutative and
if
there exists a conditional expectation $E:\A\to\vee_{i=1}^n\gamma_i(\A_i)$
such that $\phi\circ E=\phi$,
then
$$H_\phi(\gamma_1,\dots,\gamma_n)\geq S(\phi
\upharpoonright_{\vee_{i=1}^n\gamma_i(\A_i)}),$$
where the r.h.s. denotes the classical m.t. entropy of the restriction of $\phi$ 
to ${\vee_{i=1}^n\gamma_i(\A_i)}$.}
  \medskip

\noindent{\it Proof} Let us define  $\B={\vee_{i=1}^n\gamma_i(\A_i)}$, 
$\B_i=\gamma_i(\A_i)$, $\mu=\phi\upharpoonright_\B$. Then $$(\B, \mu, \B_1,\dots,\B_n)$$
is an Abelian model for $(\A,\phi,\gamma_1,\dots,\gamma_n)$. Let $E_i:\B\to\B_i$
denote the canonical conditional expectation associated to $\mu$. Then
$E_i\circ E\circ \gamma_i: \A_i\to \gamma_i(\A_i)$ conincides with $\gamma_i$, which is a $^*$--isomorphism,
thus its {\it entropy defect\/} is zero (see \cite{CNT}, section II).
It follows from [CNT, Definition III.4] that 
$$H_\phi(\gamma_1,\dots,\gamma_n)\geq S(\phi
\upharpoonright_{\vee_{i=1}^n\gamma_i(\A_i)}).$$
\medskip

Let now $\sigma$ be a ucp map of our $C^*$--algebra $\A$
such that $\phi\circ\sigma=\phi$,
and let $\gamma: M\to\A$ be a ucp map from a finite--dimensional
$C^*$--algebra $M$.
Define the m.t. dynamical entropy of $\gamma$ with respect to
$\phi$ to be
$$h_{\phi,\sigma}(\gamma)=\lim_n \frac{1}{n} 
H_\phi(\gamma,\sigma\circ\gamma,\dots, \sigma^{n-1}\circ\gamma),$$
and, finally, define the m.t. dynamical entropy of $\sigma$ as
$h_\phi(\sigma)=\sup_\gamma\{h_{\phi,\sigma}(\gamma)\}$, where the
supremum is taken over
all possible $\gamma: M\to\A$.
\medskip

\noindent{\bf 8.2. Corollary} {\sl 
Let $\A$ be a unital $C^*$--algebra,  $\phi$ a state of $\A$,
and let $\sigma$ be a ucp map of $\A$ such that $\phi\circ\sigma=\phi$.
Let $\gamma: M\to\A$ be a unital $^*$--monomorphism
from a commutative finite--dimensional $C^*$--algebra $M$.
Assume that the smallest $\sigma$--stable $C^*$--subalgebra  $\C$ of $\A$
containing $\gamma(M)$ 
is commutative and that $\sigma\upharpoonright_\C$ is a $^*$--monomorphism.
 If, for $n\in{\Bbb N}$, there exists a 
conditional expectation
$E_n:\A\to\C_n:=\vee_{i=0}^{n-1}\sigma^i\circ\gamma(M)$
such that $\phi\circ E_n=\phi$, then 
$$h_{\phi,\sigma}(\gamma)\geq
 h_{\phi\upharpoonright_\C}(\sigma\upharpoonright_\C, \gamma(M)),$$
where the r.h.s. denotes the classical m.t. dynamical entropy of the
partition
of the spectrum of $\C$
defined by $\gamma(M)$ with respect to $\sigma\upharpoonright_\C$ (see,
e.g., [DGS, Def. 10.8].
It follows that
$$h_\phi(\sigma)\geq h_{\phi\upharpoonright_\C}(\sigma\upharpoonright_\C),$$
where the r.h.s. denotes the classical m.t. entropy of the epimorphism
of the spectrum of $\C$ defined by the restriction of $\sigma$
([DGS, Def. 10.10].}\medskip

\noindent{\it Proof} Just apply the previous lemma to
$\gamma_1=\gamma,\dots,\gamma_n=\sigma^{n-1}\gamma$ and then pass to the limit.
The last assertion is a consequence of the classical Kolmogorov--Sinai
property of  the entropy.
\medskip

In order to apply the above result, one needs to know under which conditions 
on the system $(\A, \sigma, \gamma)$ as in Cor. 8.2 every invariant
measure $\mu$ on $\C$
extends to a $\sigma$--invariant state $\phi$ on $\A$ fulfilling all the 
requirements of the previous Corollary.
We start giving a well known  sufficient condition for the existence 
of invariant conditional expectations. \medskip

\noindent{\bf 8.3. Lemma} {\sl Let $M\subset\A$ be a unital inclusion of
$C^*$--algebras, with $M$ commutative and finite--dimensional, and let
$\phi$ be a  state on 
$\A$ faithful on  $M$ and such that
 $$\phi(am)=\phi(ma), \quad m\in M.$$ Then
there is a unique conditional expectation $E:\A\to M$ such that
$\phi\circ E=\phi$.}\medskip

\noindent{\it Proof} Set $E(a)=\sum \phi(e)^{-1}\phi(ae)e$, and check
that $E$ is the desired conditional expectation. Uniqueness follows easily from
faithfulness of $\phi$ on $M$.\medskip

We next give a condition on  $(\A, \sigma, \gamma)$ so that every invariant measure 
on the spectrum of  $\C$ extends to a $\sigma$--invariant state on $\A$ containing $\C$
in its centralizer. In view of the previous Lemma, this would imply
the existence of conditional expectations $E_n$ as in Corollary 8.2, 
satisfying all the
necessary
 requirements. \medskip

\noindent{\bf 8.4. Proposition} {\sl Let $\A$ be a unital
 $C^*$--algebra endowed with a ucp map $\sigma$,
and let $\C\subset\A$ be a unital 
$\sigma$--stable, AF, 
 commutative, $C^*$--subalgebra. If 
$$\sigma(ca)=\sigma(c)\sigma(a), \quad {c\in\C, a\in\A},$$
then every state $\mu$ on $\C$ satisfying $\mu\circ\sigma=\mu$ extends 
to a state $\phi$ on $\A$ such that
\roster
\item $$\phi\circ\sigma=\phi,$$
\item $$\phi(ca)=\phi(ac),\quad c\in\C, a\in\A.$$
\endroster
In particular, if $\mu$ is faithful, for every finite--dimensional
$C^*$--subalgebra $M\subset\C$ there exists a unique conditional
expectations $E_M: \A\to M$ such that
$$\phi\circ E_M=\phi.$$}\medskip

\noindent{\it Proof} Let $\C_1\subset\C_2\subset\dots$ be an increasing
sequence 
of unital finite--dimensional $C^*$--subalgebras of $\C$ with dense union,
and, for $n\in{\Bbb N}$,
let $F_n$ be the set of minimal projections of $\C_n$. 
Set
$$K_0=\{\phi\in\S(\A): \phi\upharpoonright_\C=\mu\},$$
which is a nonempty convex and compact subset of the state space
$\S(A)$ in
the weak$^*$--topology.
The function $f_0$ taking any state $\phi$ on $\A$ to 
the state $a\to\sum_{e\in F_1}\phi(eae)$ is weakly$^*$--continuous
and leaves $K_0$ invariant, thus, by the Schauder--Tychonov
fixed point theorem, the fixed point set
$$K_1=\{\phi\in K_0: f_0(\phi)=\phi\}$$
is nonempty. Note that $K_1$ is still compact and convex.
Define now the weakly$^*$--continuous function $f_1:\S(\A)\to\S(\A)$ taking $\phi$
to $a\in\A\to\phi(\sum_{e\in F_2} eae)$ and check again that 
$K_1$ is invariant under $f_1$, so that
$$K_2=\{\phi\in K_1: f_1(\phi)=\phi\}$$
is nonempty. We thus find iteratively a decreasing sequence $K_0\supset K_1\supset K_2\dots$
of nonempty compact convex subsets of $\S(\A)$. Consider
the compact convex set $K:=\cap_{n\in{\Bbb N}} K_n$.
A state $\phi$ is in $K$ if and only if
$$\phi\upharpoonright_\C=\mu,$$
$$\phi(ae)=\phi(ea),\quad e\in\cup_n F_n, a\in\A,$$
and therefore, being $\cup_n\C_n$ dense in $\C$,
$$\phi(ca)=\phi(ac), \quad c\in\C, a\in\A.$$
We next define the weakly$^*$--continuous function $f_\sigma:\S(\A)\to\S(\A)$ taking 
$\phi$ to $\phi\circ\sigma$ and we check that $f_\sigma$ leaves $K$ invariant.
First, any $\phi\in K$ restricts to $\mu$ on $\C$, thus, being $\C$ and $\mu$ $\sigma$--invariant,
$f_\sigma(\phi)$ restricts to $\mu$ on $\C$, as well. We are left to show
that for $\phi\in K$, $\C$ is in the centralizer of $f_\sigma(\phi)$.
We compute, for $c\in\C$, $a\in\A$, 
$$f_\sigma(\phi)(ca)=\phi(\sigma(ca))=\phi(\sigma(c)\sigma(a))=$$
$$\phi(\sigma(a)\sigma(c))=\phi(\sigma(ac))=f_\sigma(\phi)(ac).$$
Thus, applying again the Schauder--Tychonov fixed point Theorem,
we find a fixed point $\phi$ of $f_\sigma$, which is the desired 
extension of $\mu$.
The last assertion now follows from the previous lemma.\medskip

We now collect all the results we have obtained, in form of a Theorem. \medskip

\noindent{\bf 8.5. Theorem} {\sl 
Let $\A$ be a unital $C^*$--algebra,  
endowed with  a faithful ucp map $\sigma$, 
and let $\gamma: M\to\A$ be a unital $^*$--monomorphism
from a commutative finite--dimensional $C^*$--algebra $M$.
Assume that the smallest $\sigma$--stable $C^*$--subalgebra  $\C$ of $\A$
containing $\gamma(M)$ 
is commutative and that 
$$\sigma(ca)=\sigma(c)\sigma(a),\quad c\in\C, a\in\A.$$ 
Let $\mu$ be a faithful $\sigma$--invariant state on $\C$ extended to 
a $\sigma$--invariant state $\phi$ on $\A$ centralized by $\C$, this
being possible by Prop. 8.4.
 Then 
$$h_\phi(\sigma)\geq h_{\mu}(\sigma\upharpoonright_\C),$$
where the r.h.s. denotes the classical m.t. entropy of the epimorphism
of the spectrum of $\C$ defined by the restriction of $\sigma$ to $\C$.}\medskip

We now go back to the situation where 
 $\sigma=\sigma_{\{x_j\}}$ is the ucp map
implemented by  a finite subset $\{x_j\}$ of $\A$ such that $\sum_j
x_j{x_j}^*=I$.\medskip

\noindent{\bf 8.6. Theorem} {\sl Let $\A$ be a unital $C^*$--algebra, and let
$\{x_j\}$ by a finite subset constituted by $d$ nonzero partial isometries
satisfying
$$\sum_jx_j{x_j}^*=I,$$
$$\sum_j{x_j}^*x_j\ \text{is invertible},$$
$${x_i}^*x_j=0, i\neq j,$$
$$[{x_{i_1}\dots x_{i_r}}({x_{i_1}\dots x_{i_r}})^*, {x_j}^*x_j]=0, 
\quad j, i_1,\dots,i_r=1,\dots,d, r\in{\Bbb N}.$$
Let $\sigma_{\{x_j\}}$ be the  ucp map
implemented by the multiplet $\{x_j\}$. 
Let $$\Phi:{\C(\Lambda_{\{x_j\}}}_+)\to\A$$ be the natural
$^*$--monomorphism
 defined in Corollary 7.2. Induce a shift--invariant measure $\mu_+$
on ${\Lambda_{\{x_j\}}}_+$ from a shift--invariant measure $\mu$
on ${\Lambda_{\{x_j\}}}$ and then extend $\mu_+$ to a  
$\sigma$--invariant state $\phi$ on $\A$
centralized
by $\Phi(\C({\Lambda_{\{x_j\}}}_+)),$ by Prop. 8.4. Then
$$h_\phi(\sigma_{\{x_j\}})\geq
h_\mu(\sigma\upharpoonright_{\Lambda_{\{x_j\}}}).$$}\medskip

\noindent{\it Proof} Consider the finite 
dimensional commutative $C^*$--algebra $M$ of 
$\C({\Lambda_{\{x_j\}}}_+)$ generated by the characteristic functions of
the cylinder
sets $[i], i: x_i\neq0$.  $\C({\Lambda_{\{x_j\}}}_+)$ is naturally
embedded in $\A$
via the $^*$--monomorphism $\Phi$ defined in Cor. 7.2. Clearly 
$\Phi(\C({\Lambda_{\{x_j\}}}_+))$
is
generated by the ranges of
${\sigma_{\{x_j\}}}^i\circ \Phi(M)$. Also, $\sigma$ is faithful
as $\sum_{j}{x_j}^*x_j$ is invertible.
 The commutation relations between the domain projections
and the
range projections of the iterated products of the $x_j$'s easily show 
that $\sigma_{\{x_j\}}(ca)=\sigma_{\{x_j\}}(c)\sigma_{\{x_j\}}(a)$, 
$c\in\Phi(\C({\Lambda_{\{x_j\}}}_+))$,
$a\in\A$. In particular, $\sigma_{\{x_j\}}$ is a $^*$--monomorphism on
$\Phi(\C({\Lambda_{\{x_j\}}}_+))$. Thus the previous theorem applies.
\medskip

If for example $\A=\O_A$ or, more generally, $\A=\O_\Lambda$, the
assumptions of the previous theorem hold true.

Note that, with the notation and assumptions of the previous result,  we
know that, using also Corollary 7.6,
$$h_\mu(\Lambda_{\{x_j\}})\leq h_\phi(\sigma_{\{x_j\}})\leq 
\text{ht}(\sigma_{\{x_j\}})\leq h_{\text{top}}
(\Lambda_{\{x_j\}})+\lim_\alpha\text{ht}_2(\{\phi_{x_i,
x_j}\},\omega_\alpha),\eqno(8.1)$$
where $\omega_\alpha$ is any net of finite subsets of $\A^0$
with total union. 
The middle inequality is due to Voiculescu \cite{V}.
 
In classical ergodic theory, a probability measure $m$ on a dynamical system
$(X, T)$ such that $m\circ T^*=m$ is called an equilibrium measure, or a measure
with maximal entropy, if
it maximizes the entropy, i.e. if $h_m(X)=h_{\text{top}}(X)$.
It is well known that  dynamical systems  arising from  subshifts
admit equilibrium measures, see 
\cite{DGS}.
 Applyling this fact to our subshift $\Lambda_{\{x_j\}}$,
we see that 
  there exists a shift--invariant measure $\mu$
on  $\Lambda_{\{x_j\}}$ with $$h_\mu(\Lambda_{\{x_j\}})=h_{\text{top}}(\Lambda_{\{x_j\}}).$$
Combining with the previous inequality, we obtain, under the simplifying
assumption that the second summand in $(8.1)$ vanishes,
an existence theorem of equilibrium states in the noncommutative situation
above considered. \medskip

\noindent{\bf 8.7. Corollary} {\sl Consider the same situation as in 
 Theorem 8.6.
Let $\mu$ be a shift--invariant measure on $\Lambda_{\{x_j\}}$ with
maximal entropy,
and let us extend it to a $\sigma_{\{x_j\}}$--invariant state $\phi$ on $\A$ centralized 
by $\Phi(\C({\Lambda_{\{x_j\}}}_+))$. Assume furthermore that for 
a net $\omega_\alpha$ of finite subsets of $\A^0$ with total  union,
$\text{ht}_2(\{\phi_{x_j, x_j}\}, \omega_\alpha)=0$.
Then
$$h_\mu(\Lambda_{\{x_j\}})= h_\phi(\sigma_{\{x_j\}})=
\text{ht}(\sigma_{\{x_j\}})=h_{\text{top}}(\Lambda_{\{x_j\}}).$$}\medskip

By virtue of the remark following Theorem 7.10, we obtain the following
result.\medskip

\noindent{\bf 8.8. Corollary} {\sl 
Let $\Lambda$ be a subshift of one of the following types:
\roster 
\item
Markov shift,
\item
sofic subshift
\item
$\beta$--shift.
\endroster
Let us extend an invariant measure $\mu$ on $\Lambda$ with
maximal entropy to a state $\phi$ on $\O_\Lambda$ centralized
by the canonical commutative subalgebra $\C(\Lambda_+)\subset\O_\Lambda.$
Then
$$h_\mu(\Lambda)=h_\phi(\sigma_{\{S_i\}})=\text{ht}(\sigma_{\{S_i\}})=
h_{\text{top}}(\Lambda),$$
where $\{S_i\}$ is the canonical set of generating partial isometries
of the Matsumoto algebra $\O_\Lambda$.}\medskip

\heading{9.
From KMS states to equilibrium states}
\endheading

In this section we  make an attempt to show a  closer connection
between KMS states on full periodic
$C^*$--dynamical systems
studied in sections 1--6 and equilibrium states considered in sections 7--8.
To motivate the result of this section, we consider the classical situation
of a topological dynamical system $(X, T)$ over a compact space $X$. A
Borel probability measure
$m$ on $X$ is called conformal if $m\circ T^*$ is equivalent to $m$. The study
of conformal measures is of particular importance as it leads to equilibrium states
of the system \cite{DU}. Now in our noncommutative setting, where we replace
$X$ by a unital $C^*$--algebra $\A$ endowed with a full action of the circle,
and $T$ by the ucp map $\sigma_{\{x_j\}}$, KMS states provide a natural
class of states on $\A$ which play the role of conformal measures. Indeed we
have the following immediate result.\medskip

\noindent{\bf 9.1. Proposition} {\sl Let $(\A, \gamma)$ be a full periodic $C^*$--dynamical
system over a unital $C^*$--algebra $\A$, and let $\{x_j\}$ be a finite subset of $\A^1$
such that $\sum_j x_j{x_j}^*=I$ and $\sum_j{x_j}^*x_j$ is invertible. 
If $\omega$ is a KMS state at inverse temperature $\beta$ then
$$\omega\circ\sigma_{\{x_j\}}=\omega(a\cdot)$$
where $a=e^{-\beta}\sum_j {x_j}^*x_j$ is obviously a positive and invertible
element of $\A^0$.}\medskip

We show how to produce $\sigma$--invariant states on $\A$ from KMS states
of  the system $(\A, \gamma)$.
Consider the completely positive map $S_{\{x_j\}}: a\in\A\to\sum_j{x_j}^*ax_j$ 
already considered in section 2.
Let $\omega$ be a faithful KMS state for $(\A, \gamma)$ at maximal inverse temperature 
$\beta_{\text{max}}=\log(\lambda_{\text{max}})$.
Then, for
$t> \lambda_{\text{max}}$, consider the series
$$\sum_{k=0}^{+\infty} \frac{{S_{\{x_j\}}}^k(I)}{t^{k+1}},$$
 which we claim to be  Cauchy  for every seminorm  $p_T$, $T\in\A$, where $p_T(a)=
|\omega(aT)|$,
$a\in\A$. We show the claim.
$$\omega(\sum_{n}^m\frac{{S_{\{x_j\}}}^k(I)}{t^{k+1}}T)=
{\lambda_{\text{max}}}^{-1}\sum_n^m
(\frac{\lambda_{\text{max}}}{t})^{k+1}\omega({\sigma_{\{x_j\}}}^k(T)),$$
Now the r.h.s. is converging to $0$ as $m, n\to\infty$, since $\sum
\frac{\sigma^k}{\mu^{k+1}}$ is
norm
 converging for $\mu>1$.
Being $\omega$ faithful, 
$b_t:=\sum_{k=0}^{+\infty} \frac{{S_{\{x_j\}}}^k(I)}{t^{k+1}}$ lies in the enveloping von Neumann
algebra of $\A$.
Set $$a_t=(t-\lambda_{\text{max}})b_t,\eqno(9.1)$$
so $\omega(a_t)=1$,
and define $$\omega_t:=\omega(a_t\ \cdot\ ).$$
Then for any $T\in\A$,
$$\omega_t(T)-\omega_t(\sigma_{\{x_j\}}(T))=\omega(a_t T)-{\lambda_{\text{max}}}^{-1}
\omega(S_{\{x_j\}}(a_t) T)=$$
$${\lambda_{\text{max}}}^{-1}(t-\lambda_{\text{max}})\omega((\lambda_{\text{max}}-
S_{\{x_j\}})(b_t)T)\to 0$$
for $t\to\lambda_{\text{max}}$.
So any weak$^*$--limit point $\phi$ of $\omega_t$ for $t\to\lambda_{\text{max}}$ is a 
$\sigma_{\{x_j\}}$--invariant state on $\A$. Note that if every
${x_\mu}^*x_\mu$ commutes with any $x_\nu{x_\nu}^*$, for all multiindices $\mu$ and $\nu$,
then the Banach space generated by the 
$x_\nu{x_\nu}^*$ is in the centralizer of any such $\phi$.

The next result can be regarded as an  example where the construction of
equilibrium states out of KMS states is explicit. This 
is the noncommutative analogue of the  well known relationship between
the Perron--Frobenius Theorem and equilibrium states for  Markov subshifts, see
Proposition 17.14 in \cite{DGS}. \medskip

\noindent{\bf 9.2. Theorem} {\sl Let $(\A, \gamma, {\Bbb T})$ be a
unital, full, periodic
$C^*$--dynamical system, and let $\O_A\subset\A$
be a unital ${\Bbb Z}$--graded inclusion of the Cuntz--Krieger  algebra
associated to an irreducible matrix $A$,
in  $\A$. If $\omega$ is a faithful KMS state of $\A$ at maximal
inverse temperature $\log(\lambda_{\text{max}})$, then 
\roster
\item $\lambda_{\text{max}}=r(A)$,
\item 
 $\omega_t$ is norm convergent, for $t\to{\lambda_{\text{max}}}^+$,
to a $\sigma_{\{S_j\}}$--invariant state $\phi$
 centralized by $\C({\Lambda_A}_+)$, where $\{S_j\}$ is the canonical set
of
generators of $\O_A$.
\item $\phi$
restricts on $\C({\Lambda_A}_+)$ to the unique probability measure  $\mu$
for which
$$h_\mu(\Lambda_A)=h_{\text{top}}(\Lambda_A).$$
\endroster
In particular, if for a net $\omega_\alpha$ of finite subsets of $\A^0$
with total union 
$$\text{ht}_2(\{\phi_{x_i, x_j}\}, \omega_\alpha)=0$$ then 
$$h_\mu(\Lambda_A)=h_\phi(\sigma_{\{x_j\}})=
h_{\text{top}}(\sigma_{\{x_j\}})=h_{\text{top}}(\Lambda_A)=\log(r(A)).$$}\medskip

\noindent{\it Proof} 
It is known  that Markov subshifts defined by irreducible matrices
have a unique maximal measure, see Theorem 19.14 in \cite{DGS}.
The elements $a_t$, $t\geq\lambda_{\text{max}}$, defined as in $(9.1)$
belong to the finite--dimensional $C^*$--subalgebra of $\C({\Lambda_A}_+)$
generated by 
$S_i{S_i}^*$,
the characteristic functions 
of the cylinders $[i], i=1,\dots d$. Since $\omega(a_t)=1$ and $\omega$ is
faithful,
there exists a norm--limit point $a$ of $a_t$, for
$t\to{\lambda_{\text{max}}}^+$. Inspection
shows that $a$ is an eigenvector of $S_{\{S_i\}}$ with eigenvalue
$\lambda_{\text{max}}$,
and therefore it corresponds to a left eigenvalue $(v_i)$ of $A$,
normalized so that $\omega(a)=1$. In particular, $a_t$ is convergent.
Since $\omega$ is a KMS state of $\A$, and hence
 of $\O_\Lambda$ w.r.t. the gauge action,
evaluating $\omega$ on $S_i{S_i}^*$ gives the unique, up to a scalar,
 positive right  eigenvector $(u_j)$ of $A$.
The normalization $\omega(a)=1$ yields $\sum_i u_iv_i=1$.
Evaluating $\phi$ on ${S_{i_1}\dots S_{i_r}}({S_{i_1}\dots S_{i_r}})^*$
gives
$$\phi({S_{i_1}\dots S_{i_r}}({S_{i_1}\dots
S_{i_r}})^*)=\omega(a{S_{i_1}\dots S_{i_r}}({S_{i_1}\dots S_{i_r}})^*)=$$
$$v_{i_1}\omega({S_{i_1}\dots S_{i_r}}({S_{i_1}\dots S_{i_r}})^*)=
\frac{v_{i_1}}{\lambda^r}\omega(({S_{i_1}\dots S_{i_r}})^*{S_{i_1}\dots
S_{i_r}})=$$
$$\frac{v_{i_1}}{\lambda^r}a_{i_1,i_2}\dots a_{i_{r-1},i_r}\sum_i
\sum_j\omega(a_{i_r,j}S_j{S_j}^*)=
\frac{v_{i_1}u_{i_r}}{\lambda^{r-1}}a_{i_1,i_2}\dots a_{i_{r-1},i_r}.$$
If we now compare with the formula given in Prop. 17.14 in \cite{DGS},
we see that $\mu=\phi\upharpoonright_{\C({\Lambda_A}_+)}$ restricts  
precisely to the
unique measure on ${\Lambda_A}_+$ with maximal entropy.

\bigskip

\noindent{\bf Acknowledgments}
Part of this paper was written during a visit of C.P. at the Mathematics 
Department of the
University of Orleans. She wishes to thank C. Anantharaman--Delaroche 
for invitation and 
for drawing attention to the Furstenberg's example, and J. Renault for 
many fruitful discussions.  
We are also indebted to H. Matui for pointing out an error in section 7
of a previous version of this paper.

\heading
References
\endheading 

\item{[Bl]} Blanchard, F.: $\beta$--Expansions and symbolic dynamics.
Theor. Computer Sci. {\bf 65}, 131--141 (1989).

\item{[BDR]}  Blackadar, B.,   Dadarlat, M.,   R\o rdam, M.: The real
rank
of inductive limit $C^*$-algebras.  Math. Scand. {\bf 69}, 
211--216 (1991). 

\item{[BG]}  Boca, F.,  Goldstein, P.: Topological entropy for the
canonical 
endomorphism
of the Cuntz--Krieger algebras. Preprint 1999. math.OA/9906210.

\item{[BEH]}  Bratteli, O.,  Elliott, G.,  Herman, R.: On the possible 
temperatures
of a dynamical system.  Comm. Math. Phys. {\bf 74},
281--295 
(1980).

\item{[BEK]}  Bratteli, O.,  Elliott, E.,  Kishimoto, A.: The temperature 
state space
of a $C^*$--dynamical system I.  Yokohama Math. J.
{\bf 28}, 125--167 (1980).

\item{[BJ]}  Bratteli, O.,  J\o rgensen, P.E.T.: Isometries , shifts, 
Cuntz
algebras and multiresolusion wavelet analysis of scale $N$.
Integral 
Equations Operator Theory {\bf 28}, 382--443 (1997).

\item{[B]}  Brown, N.P.: Topological entropy in exact $C^*$--algebras.
 Math. Ann. {\bf 314}, 347--367 (1999).

\item{[Ch]}  Choda, M.: Endomorphisms of shift type (entropy for 
endomorphisms
of Cuntz algebras). In: S. Doplicher, R. Longo, J.E. Roberts,
L. Zsido (eds){\it Op\-era\-tor Algebras and Quantum Field
Theory\/}. Proceedings, Rome 1996, pp. 469--475. Cambridge, MA: 
International Press.

\item{[CNT]}  Connes, A.,  Narnhofer, H.,  Thirring, W.: Dynamical
entropy of $C^*$--algebras and von Neumann algebras.  Comm. Math. 
Phys. {\bf 112}, 691--719 (1987).

\item{[Co]}  Connes, A.: Compact metric spaces: Fredholm modules, and
hyperfiniteness.  Ergod. Th. \& Dynam. Sys.  {\bf 9},
207--220
(1989).

\item{[C]}  Cuntz, J.: Simple $C^*$--algebras generated by isometries.
 Comm. Math. Phys. {\bf 57}, 173--185 (1977).

\item{[CK]}  Cuntz, J.,   Krieger, W.: A class of $C^*$--algebras
and 
topological Ma
rkov chains.
 Invent. Math. {\bf 56}, 251--268 (1980).

\item{[DGS]}  Denker, M.,  Grillenberger, C.,  Sigmund, K.: {\it Ergodic
theory 
on compact spaces\/}.
LNM 527, Springer--Verlag 1976.

\item{[DU]}  Denker, M.,   Urba\'nski, M.: On the existence of conformal
measures.
 Tans. Amer. Math. Soc. {\bf 328}, 563--587 (1991).

\item{[E]}  Evans, D.: Gauge actions on $\O_A$.  J. Operator Theory
 {\bf 7}, 79--100 (1982).

\item{[EFW]}  Enomoto, M.,  Fuji, M.,  Watatani, Y.: KMS states for
gauge 
actions on
$\O_A$.  Math. Japon. {\bf 29}, 607--619 (1984).

\item{[G]}  Gantmacher, F.R.: {\it The theory of matrices\/}, vol. 2.
Chelsea 
Publishing Company,
1964.

\item{[GP]}  Goldstein, P.,  Pinzari, C.: Work in progress.

\item{[GHJ]}  Goodman, F.M.,   de la Harpe, P.,   Jones,
V.F.R.: {\it Coxeter 
graphs and 
towers of algebras\/}. New York: Springer--Verlag, 1989.

\item{[H]}  Hutchinson, J.E.: Fractals and self--similarity.
 Indiana Univ. Math. J. {\bf 30}, 713--747 (1981).

\item{[JP]} J\o rgensen, P.E.T.,   Pedersen, S.: Harmonic analysis of
fractal measures.  Constr. Approx. {\bf 12}, 1--30
(1996) 

\item{[KPW]}  Kajiwara, T., Pinzari, C.,  Watatani, Y.:  Ideal structure 
and simplicity of the $C^*$--algebras generated by Hilbert bimodules.
 J. Funct. Anal.  {\bf 159}, 295--322 (1998).

\item{[K]}  Katayama, Y.: Generalized Cuntz algebras $\O_N^M$.  RIMS 
Kokyuroku {\bf 858}, 87--90 (1994).

\item{[KMW]} Katayama, Y., Matsumoto, K., Watatani Y.: Simple
$C^*$--algebras arising from $\beta$--expansion of real numbers.
Ergod Th.\& Dynam. Sys. {\bf 18}, 937--962 (1998).

\item{[Ma]}  Ma\~ne, R.: {\it Ergodic theory and differentiable
dynamics\/}.
Berlin Heidelberg: 
Springer--Verlag, 1987.

\item{[MP]}  Martin, M., Pasnicu, C.: Some comparability results in 
inductive limit $C^*$-algebras.  J. Operator Th.
{\bf 30}, 137--147 (1993).

\item{[M]} Matsumoto, K.: On $C^*$-algebras associated with subshifts.
 Internat.
J. Math. {\bf 8}, 357--374 (1997). 

\item{[M2]}  Matsumoto, K.: Dimension groups for subshifts and
simplicity 
of the corresponding $C^*$-algebras. To appear in  J. Math. 
Soc. Japan.

\item{[MWY]}  Matsumoto, K.,  Watatani, Y.,  Yoshida, M.: KMS states for
gauge 
actions on $C^*$--algebras
associated with subshifts.  Math.--Z. {\bf 3}, 489--509 
(1998).

\item{[OPI]}  Olesen, D.,  Pedersen, G.K.: Some $C^*$--dynamical systems 
with a single KMS state.  Math. Scand. {\bf 42}, 111-118
(1978).

\item{[OP]}  Olesen, D., Pedersen, G.K.: Applications of the Connes 
spectrum to $C^*$--dynamical
systems III.  J. Funct. Anal.  {\bf 45}, 357--390 (1982).

\item{[Par]} Parry, W.: On the $\beta$--expansions of real numbers. 
Acta Math. Acad. Sci. Hung. {\bf 11}, 401--416 (1960).

\item{[P]} Pimsner, M.: A class of $C^*$--algebras generalizing both 
Cuntz--Krieger algebras and crossed products by ${\Bbb Z}$. In:
Voiculescu, D. (ed.)
{\it Free probability theory\/}, AMS, 1997.

\item{[Re]} Renyi, A.: Representations of real numbers and their ergodic
properies. Acta Math. Acad. Sci. Hung. {\bf 8}, 477--493 (1957).

\item{[R]}  R\o rdam, M.: Classification of certain infinite simple
 $C^*$--algebras.  J. Funct. Anal. {\bf 131}, 415--458 
(1995).

\item{[RiI]}  Rieffel, M.: Metrics on states from actions of compact
groups.   Doc. Math.
{\bf 3}, 215--229 (1998);  math.OA/9807084.

\item{[RiII]} Rieffel, M.: Metrics on state spaces. Doc. Math. {\bf 4},
559--600 (1999); math.OA/9906151.

\item{[T]}  Takesaki, M.: {\it Theory of operator algebras I\/}. New York: 
Springer--Verlag,
1979.

\item{[V]}  Voiculescu, D.: Dynamical approximation entropies and 
topological entropy in operator algebras.  Comm. Math. Phys.
 {\bf 170}, 249--281 (1995).

\item{[W]}  Wassermann, A.: Exact $C^*$--algebras and related topics.
Lecure
Notes Series no. 19, GARC, Seoul National University, 1994.
\enddocument